\documentclass{article}
\usepackage{graphicx}
\usepackage{float}
\usepackage{nameref}
\usepackage{amsthm}
\usepackage{amsmath}
\usepackage{amsfonts}
\usepackage{multicol}
\usepackage{xassoccnt}
\usepackage{hyperref}
\hypersetup{colorlinks, linkcolor = {red}, citecolor = {green}, urlcolor = {blue}}
\usepackage{cleveref}
\usepackage{mathrsfs}
\usepackage{stmaryrd}
\usepackage[shortlabels]{enumitem}
\usepackage{amssymb}
\usepackage[T1]{fontenc}
\usepackage{pxfonts}
\DeclareMathAlphabet{\mathpzc}{OMS}{cmsy}{m}{n}
\usepackage{tikz}
\usepackage{tikz-cd}
\usetikzlibrary{automata, arrows.meta, positioning}
\usepackage[english]{babel}
\usepackage{comment}
\usepackage[margin=1.1in]{geometry}
\usepackage{pifont}
\newcommand{\cmark}{\ding{51}}%
\newcommand{\xmark}{\ding{55}}%

\newtheorem{theorem}{Theorem}[section]
\crefname{theorem}{theorem}{theorems}
\newtheorem{proposition}{Proposition}[section]
\crefname{proposition}{proposition}{propositions}
\newtheorem{lemma}{Lemma}[section]
\crefname{lemma}{lemma}{lemmas}
\newtheorem{corollary}{Corollary}[section]
\crefname{corollary}{Corollary}{corollaries}
\newtheorem{conjecture}{Conjecture}[section]
\crefname{conjecture}{Conjecture}{conjectures}
\newtheorem{remark}{Remark}[section]
\crefname{remark}{remark}{remarks}
\newtheorem{example}{Example}[section]
\crefname{example}{example}{examples}
\newtheorem{definition}{Definition}[section]
\crefname{definition}{definition}{definitions}

\DeclareCoupledCountersGroup{theorems}
\DeclareCoupledCounters[name=theorems]{theorem,proposition,lemma,remark,example,corollary,definition,conjecture}

\title{Smooth sequences over two letters of same parity}
\author{Julien Cassaigne\footnote{CNRS, Aix Marseille Univ, I2M, 3 place Victor Hugo, Case 19, 13331 Marseille Cedex 3, France} ~~~~~~~~~~~~ Raphaël Henry\footnote{Université de Bordeaux, IMB, 351 cours de la Libération, F 33405 Talence, France}}
\date{}

\begin{document}

\maketitle

\pagenumbering{arabic}
\setcounter{page}{1}

\begin{abstract}
	Smooth sequences over an alphabet of integers $\{a,b\}$ are sequences that are infinitely derivable, the emblematic example being the Oldenburger-Kolakoski sequence over $\{1,2\}$. In this paper we study smooth sequences when $a$ and $b$ have same parity: we give a well-behaved $S$-adic representation and we prove uniform or linear recurrence, linear factor complexity and uniform frequencies for almost all these sequences. In the process we show that their language is anti-complement-invariant and we describe precisely the factor complexity in many cases. As an application we compute the frequencies of morphic smooth sequences over any odd alphabet.
\end{abstract}

\newcommand{\kola}{
    \begin{gathered}
    \tikzpicture[every node/.style={anchor=south west}]
        \node[minimum width=1cm,minimum height=0.4cm] at (0,0) {$\kappa = \underbrace{22}_2 \underbrace{11}_2 \underbrace{2}_1 \underbrace{1}_1 \underbrace{22}_2 \underbrace{1}_1 \underbrace{22}_2 \underbrace{11}_2 \underbrace{2}_1 \underbrace{11}_2 \underbrace{22}_2 \underbrace{1}_1 \underbrace{2}_1 \underbrace{22}_2 \underbrace{1}_1 ...$};
        \node[minimum width=1cm,minimum height=0.4cm] at (14.15,0.05) {$... = \kappa$};
    \endtikzpicture
    \end{gathered}
}

\newcommand{\kolabis}{
    \begin{gathered}
    \tikzpicture[every node/.style={anchor=south west}]
        \node[minimum width=1cm,minimum height=0.4cm] at (0,0) {$\kappa_{3,1} = \underbrace{333}_3 \underbrace{111}_3 \underbrace{333}_3 \underbrace{1}_1 \underbrace{3}_1 \underbrace{1}_1 \underbrace{333}_3 \underbrace{111}_3 \underbrace{333}_3 \underbrace{1}_1 \underbrace{333}_3 \underbrace{1}_1 \underbrace{333}_3 \underbrace{111}_3 ...$};
        \node[minimum width=1cm,minimum height=0.4cm] at (13.575,-0.05) {$... = \kappa_{3,1}$};
    \endtikzpicture
    \end{gathered}
}

\section*{Introduction}

The emblematic Oldenburger-Kolakoski sequence is a fixed point of the \textit{run-length encoding} over the alphabet $\{1,2\}$:
$$\kola$$
Initially studied by Oldenburger \cite{Oldenburger}, it was later popularized by Kolakoski \cite{Kolakoski} and is now registered as A000002 in the OEIS \cite{OEIS}. Early on, the community of combinatorics on words asked the following questions: Does $\kappa$ belong to a known class of sequences? Is it eventually periodic? (uniformly) recurrent? Does it have factor frequencies? What is its factor complexity? Kolakoski's original question had already been answered in \cite{Oldenburger}: $\kappa$ is not eventually periodic. Also, in \cite{CKL} the authors show that $\kappa$ is a 2-alternating D0L sequence but not a purely morphic sequence. The other questions remain unsolved.

Facing the difficulties raised by $\kappa$, researchers started to investigate its variations in the hope of finding better answers. The first idea is to change the alphabet: if $a$ and $b$ are two different integers, let $\kappa_{a,b}$ denote the unique fixed point of the run-length starting with $a$. With this notation, the Oldenburger-Kolakoski sequence is denoted by $\kappa_{2,1}$ and we observe that $\kappa_{1,2} = 1 \kappa_{2,1}$. The next sequence to bring attention, registered as A064353 in the OEIS \cite{OEIS}, is
$$\kolabis$$
In this paper we shall consider a family of sequences originally introduced over $\{1,2\}$ by Dekking \cite[Section 4]{Dekking01}, which contains the sequences $\kappa_{a,b}$. Over any binary alphabet $\mathpzc{A} = \{a,b\}$, if the run-length encoding is \textit{derivative}, then a sequence in $\mathpzc{A}^\mathbb{N}$ is \textit{smooth} if all its derivatives are in $\mathpzc{A}^\mathbb{N}$. We write $\mathcal{C}^\infty$ the set of smooth sequences. Over any binary alphabet, we also consider the computable language of \textit{smooth words} denoted $\mathcal{C}_f^\infty$, which is closed by an operation called \textit{finite derivative} that simulates locally the run-length encoding. The connection between $ \mathcal{C}_f^\infty$ and $\mathcal{C}^\infty$ is that $\mathscr{L}(\mathcal{C}^\infty) = \mathcal{C}_f^\infty$ (we have $\mathscr{L}(\mathcal{C}^\infty) \subseteq \mathcal{C}_f^\infty$ by construction, and we have proven the equality in \cite{CH}).

Over \textit{mixed} alphabets, i.e., $\mathpzc{A} = \{a,b\}$ where $a+b$ is odd, very few things are known about smooth sequences and in particular about the Oldenburger-Kolakoski sequence. The first main conjecture is that $\mathcal{C}_f^\infty$ is uniformly recurrent, which would imply that every smooth sequence $x$ satisfies $\mathscr{L}(x) = \mathcal{C}_f^\infty$ and is uniformly recurrent. Note that we have proven in \cite{CH} that $p_{\mathcal{C}_f^\infty}(n)$ grows at least like $n^\rho$ where $\rho = \frac{\log(a+b)}{\log\left(\frac{a+b}{2}\right)}$. The second main conjecture is that  $\mathcal{C}_f^\infty$ has uniform frequencies, which generalizes Keane's conjecture that $\kappa$ has uniform letter frequencies.
	
Over \textit{even} and \textit{odd} alphabets, i.e., $\mathpzc{A} = \{a,b\}$ where $a$ and $b$ are both even (res. both odd), smooth sequences seem to have different properties. In \cite{Dekking79,Dekking01} Dekking noticed that $\kappa_{3,1}$ is morphic, then in \cite{Sing02,BS,Sing10} Sing showed the same for every fixed point $\kappa_{a,b}$ using other substitutions that are always primitive. It immediatly follows that these smooth sequences are linearly recurrent, have linear factor complexity and uniform frequencies, in particular they have computable algebraic letter frequencies (the frequency of the letter 3 in $\kappa_{3,1}$ is approximately 0.6028). The substitutive nature of smooth sequences extends beyond the fixed points: in \cite{JMPR} the authors provide an $S$--adic representation of all smooth sequences over $\{1,3\}$, using a coding and two substitutions over four letters that generalize Dekking's construction. The authors deduce that these smooth sequences have uniform frequencies and are almost all uniformly recurrent, and they characterize their letter frequencies. More generally, in \cite{BJP} the authors study smooth sequences over any even or odd alphabet: they show that all these smooth sequences are recurrent, they characterize the minimal and maximal smooth sequences with respect to the lexicographic order and they compute the letter frequencies of these two sequences.

In this paper we combine the approaches of \cite{Sing02,BS,Sing10} about morphic fixed points, of \cite{JMPR} about the $S$-adic representation and of \cite{BJP} that considers any even or odd alphabet. 
\begin{itemize}
	\item With \Cref{smoothSadicSystem} we use the substitutions introduced by Sing in \cite{Sing02,BS,Sing10} to give a general $S$-adic representation of smooth sequences over any even and odd alphabet.
\end{itemize}
We use these substitutions because they behave better than the ones used in \cite{JMPR}: they allow to apply general results about $S$-adic sequences to a majority of smooth sequences. Combined with arguments of our own, we improve the previous results and obtain new ones about smooth sequences.
\begin{itemize}
	\item With \Cref{antiInvariance} we show that the language of every smooth sequence over an even and odd alphabet is anti-complement-invariant.

	\item With \Cref{URLR} we show that smooth sequences over an even or odd alphabet $\{a,b\}$ are all uniformly recurrent and are linearly recurrent if and only if they belong to the class $\mathscr{M}$ (cf. \Cref{propertyM}). This class covers all smooth sequences if $a > 1$ and a minority of smooth sequences if $a = 1$.
	
	\item With \Cref{sameComplexity} we show that, over an even alphabet or an odd alphabet of the form $\{a,a+2\}$, all smooth sequences have the same factor complexity $p(n)$ and we characterize $\liminf p(n)/n$ and $\limsup p(n)/n$. We deduce with \Cref{linearCorollary} that the smooth sequences in $\mathscr{M}'$ (cf. \Cref{propertyM2}) have linear factor complexity, where the class $\mathscr{M}'$ adds all smooth sequences over $\{1,3\}$ to $\mathscr{M}$.
	
	\item With \Cref{frequenciesCorollary} we show that the smooth sequences in $\mathscr{M}'$ also have uniform frequencies. Furthermore, as an application of \Cref{smoothSadicSystem}, with \Cref{smoothFreq} we show that the letter frequencies of all smooth sequences over an even alphabet are $1/2$ and we compute the letter frequencies of the morphic smooth sequences over odd alphabets.
\end{itemize}

We summarize the new properties of smooth sequences in \Cref{summary}.

\subsubsection*{Outline of the paper}

In \Cref{prelim} we provide the basic definitions of our paper and we define smooth sequences and smooth words over binary alphabets. Then we detail the difference between smooth sequences over mixed alphabets and over even or odd alphabets by recalling the existing conjectures and results, and finally we formulate the results of this paper.
\newline
In \Cref{structure} we prove \Cref{smoothSadicSystem} and we characterize the morphic smooth sequences.
\newline
In \Cref{occurrences} we prove \Cref{antiInvariance} among other technical results that will serve as a toolbox in the rest of the paper. 
\newline
In \Cref{recurrence} we prove \Cref{URLR}.
\newline
In \Cref{complexity} we prove \Cref{sameComplexity} and \Cref{linearCorollary}.
\newline
In \Cref{frequencies} we prove \Cref{frequenciesCorollary} and \Cref{smoothFreq}.

\section{Preliminaries}\label{prelim}

\subsection{Notation and basic concepts}

\subsubsection{Words and sequences}

\begin{definition}
	A binary alphabet is a set $\mathpzc{A} = \{a,b\}$ where $a$ and $b$ are integers such that $1 \leq a < b$. If $a$ and $b$ are both even (resp. both odd), we say that $\mathpzc{A}$ is \textup{even} (resp. \textup{odd}). Otherwise, i.e., if $a+b$ is odd, we say that $\mathpzc{A}$ is \textup{mixed}.
\end{definition}

Given a finite alphabet $\mathpzc{A}$, $\mathpzc{A}^n$ (resp. $\mathpzc{A}^*$) denotes the set of finite words of length $n$ (resp. the set of all finite words) over $\mathpzc{A}$, and $\lvert u \rvert$ denotes the length of the word $u \in \mathpzc{A}^*$. In particular $\varepsilon$ is the empty word and we define $\mathpzc{A}^+ := \mathpzc{A}^*\backslash\{\varepsilon\}$. If words $u,v,w,w' \in \mathpzc{A}^*$ are such that $v = wuw'$, we say that $u$ is a \textit{factor} of $v$ and we write $u \sqsubset v$. Moreover, if $w = \varepsilon$ (resp. $w' = \varepsilon$), we say that $u$ is a \textit{prefix} (resp. \textit{suffix}) of $v$ and we write $u \sqsubset_p v$. Over a binary alphabet $\mathpzc{A} = \{a,b\}$, the \textit{complement} of a word $u \in \mathpzc{A}^*$ (resp. $x \in \mathpzc{A}^\mathbb{N}$) is the word $\overline{u}$ (resp. $\overline{x}$) where every letter $a$ is replaced with the letter $b$ and vice versa. If $u = u_1 ... u_n \in \mathpzc{A}^n$ and $c \in \mathpzc{A}$, we define $\lvert u \rvert_c := \left\lvert \left\{ i \in \left\llbracket 1,n \right\rrbracket ~\middle|~ u_i = c \right\} \right\rvert$.

Let $\mathpzc{A}^\mathbb{N}$ be the set of all (right-infinite) sequences over $\mathpzc{A}$, where $0 \in \mathbb{N}$. The \textit{shift} over $\mathpzc{A}^\mathbb{N}$ is the map $S : (x_n)_{n \geq 0} \mapsto (x_{n+1})_{n \geq 0}$. The factors of a sequence $x \in \mathpzc{A}^\mathbb{N}$ are the finite words $u$ of the form $x_{[i,j)} := x_i x_{i+1} ... x_{j-1}$ for $0 \leq i \leq j$, with the convention that $x_{[i,i)} := \varepsilon$, and we write $u \sqsubset x$. If $u$ is of the form $x_{[0,j)}$, we say that $u$ is a \textit{prefix} of $x$ and we write $u \sqsubset_p x$. The \textit{language} of a sequence $x \in \mathpzc{A}^\mathbb{N}$ is the set of its factors ${\mathscr{L}(x) := \{u \in \mathpzc{A}^* \mid u \sqsubset x\}}$. Similarly, the language of a set $X \subseteq \mathpzc{A}^\mathbb{N}$ is the set $\mathscr{L}(X) := \displaystyle \bigcup_{x \in X} \mathscr{L}(x)$.

A sequence $x \in \mathpzc{A}^\mathbb{N}$ is \textit{recurrent} if each of its factors occurs infinitely often in $x$. A word $u \in \mathpzc{A}^*$ \textit{occurs with bounded gaps} in $x$ if there exists $L \geq 1$ such that, for all $v \in \mathscr{L}(x)$, if $\lvert v \rvert \geq L$ then $u \sqsubset v$. Stronger than being recurrent, $x$ is \textit{uniformly recurrent} if all its factors occur with bounded gaps in it. Even stronger than being uniformly recurrent, $x$ is \textit{linearly recurrent} if there exists a constant $L \geq 1$ such that, for all $u,v \in \mathscr{L}(x)$, if $\lvert v \rvert \geq L \lvert u \rvert$ then $u \sqsubset v$.

\subsubsection{Morphic and $S$-adic sequences}

Let $\mathpzc{A}$ and $\mathpzc{B}$ be two finite alphabets. A map $\sigma : \mathpzc{A}^* \rightarrow \mathpzc{B}^*$ is a \textit{substitution} if $\sigma(uv) = \sigma(u)\sigma(v)$ for all $u,v \in \mathpzc{A}^*$ and $\sigma(a) \neq \varepsilon$ for all $a \in \mathpzc{A}$, and we write $\sigma : \mathpzc{A} \rightarrow \mathpzc{B}^+$. If $\sigma : \mathpzc{A} \rightarrow \mathpzc{A}^+$ and there is a letter $a \in \mathpzc{A}$ and a word $u \in \mathpzc{A}^+$ such that $\sigma(a) = au$, then iterating $\sigma$ on the letter $a$ produces prefixes of increasing length of an sequence in $\mathpzc{A}^\mathbb{N}$. The produced sequence is called a \textit{purely morphic sequence} and is denoted $\sigma^\infty(a)$. If we apply another substitution $\tau : \mathpzc{A} \rightarrow \mathpzc{B}^+$, then the sequence $\tau(\sigma^\infty(a)) \in \mathpzc{B}^\mathbb{N}$ is called a \textit{morphic sequence}.

Let $\mathpzc{A}$ be an alphabet, let $a$ be a letter of $\mathpzc{A}$, let $S$ be a finite set of substitutions $\sigma$ from $\mathpzc{A}(\sigma) \subset \mathpzc{A}$ to $\mathpzc{A}^+$ and let $\left(\sigma_n : \mathpzc{A}_{n+1} \rightarrow \mathpzc{A}_n^+\right)_{n \geq 0}$ be a sequence of $S^\mathbb{N}$ such that $\left(\sigma_0 \sigma_1 ...  \sigma_n (a^\omega)\right)_{n \geq 0}$ converges in $\mathpzc{A}^\mathbb{N}$ to $x$. We will say that $x$ is an \textit{S-adic sequence} on $\mathpzc{A}$ (generated by $\left(\sigma_i\right)_{i \geq 0} \in S^\mathbb{N}$ and $a$). If there exists an integer $s_0$ such that for all $r \geq 0$, all $b \in \mathpzc{A}_r$ and all $c \in \mathpzc{A}_{r+s_0+1}$, the letter $b$ has an occurrence in $\sigma_{r+1} \sigma_{r+2} ... \sigma_{r+s_0} (c)$, then we say that $x$ is a \textit{primitive $S$-adic sequence} (with constant $s_0$). We say that a substitution $\sigma : \mathpzc{A} \rightarrow \mathpzc{B}^+$ is \textit{proper} if there exist two letters $\ell,r \in \mathpzc{B}$ such that for all $a \in \mathpzc{A}$ the first (resp. last) letter of $\sigma(a)$ is $\ell$ (resp. $r$). We say that a sequence $x \in \mathpzc{A}^\mathbb{N}$ is a \textit{proper $S$-adic sequence} if it is an $S$-adic sequence and the substitutions in \textit{S} are proper.

\subsubsection{Factor complexity and bispecial words}

A language $\mathscr{L} \subset \mathpzc{A}^*$ is \textit{factorial} if $u \sqsubset v \in \mathscr{L}$ implies $u \in \mathscr{L}$, and it is left-extendable (resp. right-extendable) if $u \in \mathscr{L}$ implies that there exists $a \in \mathpzc{A}$ such that $au \in \mathscr{L}$ (resp. $ua \in \mathscr{L}$). In particular, the language $\mathscr{L}(x)$ of a word $x \in \mathpzc{A}^\mathbb{N}$ is factorial and right-extendable, and if in addition $x$ is recurrent then it is also left-extendable.

\begin{definition}
	We define the following complexity functions $p : \mathbb{N} \rightarrow \mathbb{N}$.
	\begin{itemize}
		\item The \textup{complexity} of a language $\mathscr{L} \subset \mathpzc{A}^*$ is
	${p_\mathscr{L}(n) := \left\lvert \mathscr{L} \cap \mathpzc{A}^n \right\rvert}$,
	
		\item The \textup{factor complexity} of a word $x \in \mathpzc{A}^\mathbb{N}$ is $p_x(n) := p_{\mathscr{L}(x)}(n)$,
		
		\item The \textup{factor complexity} of a set of words $X \subset \mathpzc{A}^\mathbb{N}$ is $p_X(n) := p_{\mathscr{L}(X)}(n)$.
	\end{itemize}
\end{definition}

In order to describe the asymptotics of a complexity function, we use the following notations.

\begin{definition}
	Given a language $\mathscr{L} \subset \mathpzc{A}^*$ and a function $f : \mathbb{N} \rightarrow \mathbb{N}$, we write
	\begin{itemize}
		\item $p_\mathscr{L}(n) = \mathcal{O}(f(n))$ if there exists $C > 0$ and $N \geq 0$ such that ${p_\mathscr{L}(n) \leq C f(n)}$ for all $n \geq N$,
		\item $p_\mathscr{L}(n) = \Omega(f(n))$ if there exists $C > 0$ and $N \geq 0$ such that ${p_\mathscr{L}(n) \geq C f(n)}$ for all $n \geq N$,
		\item ${p_\mathscr{L}(n) = \Theta(f(n))}$ if ${p_\mathscr{L}(n) = \mathcal{O}(f(n))}$ and ${p_\mathscr{L}(n) = \Omega(f(n))}$.
	\end{itemize}
\end{definition}

Let us recall how to compute the complexity of a language over a binary alphabet from the bispecial factors (we refer to \cite[Section 3.3]{Cassaigne}). First, a language $\mathscr{L} \subseteq \mathpzc{A}^*$ is \textit{factorial} if $u \sqsubset v \in \mathscr{L}$ implies $u \in \mathscr{L}$, and it is \textit{left-extendable} (resp. \textit{right-extendable}) if, for all $u \in \mathscr{L}$, there exists $a \in \mathpzc{A}$ such that $au \in \mathscr{L}$ (resp. $ua \in \mathscr{L}$).

\begin{definition}
	Let $\mathpzc{A} = \{a,b\}$ be a binary alphabet and let $\mathscr{L} \subseteq \mathpzc{A}^*$ be a factorial, left- and right-extendable language. A word $u \in \mathscr{L}$ is \textup{bispecial} if ${\{au,bu,ua,ub\} \subseteq \mathscr{L}}$. In that case, the set $\{aua,aub,bua,bub\} \cap \mathscr{L}$ has 2, 3 or 4 elements and its cardinal determines the \textup{type} of $u$: we say that $u$ is \textup{weak} if it has 2 elements, \textup{neutral} if it has 3 and \textup{strong} if it has 4.
\end{definition}

Let $s(n)$ be the first finite difference of the factor complexity, i.e., $s(n) := p_\mathscr{L}(n+1) - p_\mathscr{L}(n)$, and let $b(n)$ be the second finite difference, i.e., $b(n) := s(n+1) - s(n)$. The seminal result that allows to compute the complexity of the language $\mathscr{L}$ is the following.

\begin{theorem}[{\cite[Proposition 3.2]{Cassaigne}}]\label{pfrombsw}
	Let $\mathscr{L}$ be a factorial, left- and right-extendable language over a binary alphabet, and let $bs(n)$ (resp. $bw(n)$) denote the number of strong (resp. weak) bispecial words of $\mathscr{L}$. Then for all $n \geq 0$ we have
	$$b(n) = bs(n) - bw(n).$$
\end{theorem}

\subsection{Smooth sequences and smooth words}

\subsubsection{Smooth sequences}

Over a binary alphabet $\mathpzc{A} = \{a,b\}$, any sequence $x \in \mathpzc{A}^\mathbb{N}$ where the letters $a$ and $b$ occur infinitely often is canonically factorized as ${x_0}^{p_0} \overline{x_0}^{p_1} {x_0}^{p_2} \overline{x_0}^{p_3} ...$ where $x_0$ is the first letter of $x$ and $p_i \geq 1$ for all $i \geq 0$.

\begin{definition}
	Over a binary alphabet $\mathpzc{A}$, we say that a sequence $x \in \mathpzc{A}^\mathbb{N}$ is \textup{derivable} if its canonical factorization is ${x_0}^{p_0} \overline{x_0}^{p_1} {x_0}^{p_2} \overline{x_0}^{p_3}...$ with $p_i \in \mathpzc{A}$ for all $i \geq 0$. We write $\mathcal{C}^1$ the set of derivable sequences and we define the \textup{derivative} as the map
	\begin{align*}
		\mathcal{D}~:~~~~~~~~~~~~~~~~~~ \mathcal{C}^1 & \longrightarrow \mathpzc{A}^\mathbb{N} \\
    	{x_0}^{p_0} \overline{x_0}^{p_1} {x_0}^{p_2} ... & \longmapsto p_0 p_1 p_2...
	\end{align*}
\end{definition}

\begin{example}
	Over $\{1,2\}$, we have $(221)^{\omega} \in \mathcal{C}^1$ and $\mathcal{D}\left((221)^{\omega}\right) = (21)^{\omega} \in \mathcal{C}^1$ but $\mathcal{D}\left((21)^{\omega}\right) = 1^{\omega} \notin \mathcal{C}^1$.
\end{example}

\begin{definition}
	For $n \geq 2$, we define by induction the set of sequences that are derivable $n$ times
	\begin{equation*}
		\mathcal{C}^n := \left\{x \in \mathcal{C}^1 ~\middle |~ \mathcal{D}(x) \in \mathcal{C}^{n-1} \right\}.
	\end{equation*}
	We then define the set of \textup{smooth sequences}, that is the set of sequences that are infinitely derivable
	\begin{equation*}
		\mathcal{C}^\infty := \displaystyle \bigcap_{n \geq 1} \mathcal{C}^n.
	\end{equation*}
\end{definition}

\begin{example}\label{kab}
	As a primary example of smooth sequences, over any binary alphabet $\mathpzc{A} = \{a,b\}$, the derivative $\mathcal{D}$ has two fixed points: one that starts with $a$, denoted by $\kappa_{a,b}$; the other that starts with $b$, denoted by $\kappa_{b,a}$. Here are notable examples, among which is the Oldenburger-Kolakoski word $\kappa_{2,1}$:
	\begin{align*}
		\kappa_{2,1} = &~2211212212211211221211212211211212212211212212112112212211212212211211212212112... \\
		\kappa_{3,1} = &~3331113331313331113331333133311133313133311133313331113331333111333131333111333... \\
		\kappa_{2,4} = &~2244222244442244224422224444222244442244222244442244222244442244224422224444222... \\
		\kappa_{2,5} = &~2255222225555522552255225555522222555552222255555225522222555552255222225555522...
	\end{align*}
	Also, we notice that $\kappa_{1,b} = 1 \kappa_{b,1}$ for all $b \geq 2$.
\end{example}

Oldenburger originally proved in \cite[Theorem 2]{Oldenburger} that $\kappa_{2,1}$ is not eventually periodic. It was proven again by Üçoluk in \cite{Ucoluk} and by Dekking in \cite[Example 4]{Dekking79}, and the proof is easily adaptable to any smooth sequence.

\begin{proposition}\label{aperiodic}
	No smooth sequence is eventually periodic.
\end{proposition}
\begin{proof}
	Suppose that a word $x \in \mathpzc{A}^\mathbb{N}$ is smooth, eventually periodic, and that its period $u \in \mathpzc{A}^+$ is the smallest period among all smooth sequences. Then $\mathcal{D}(x)$ is also smooth and eventually periodic with a strictly smaller period than $x$, unless $\mathpzc{A} = \{1,b\}$ and $u \in \{1b,b1\}$. In the first case it contradicts the minimality of $u$, and in the second case $\mathcal{D}(x)$ is not smooth.
\end{proof}

At last, we must warn the reader that the set $\mathcal{C}^\infty$ is not a subshift.

\begin{remark}
	The set $\mathcal{C}^\infty$ is closed but it is not shift-invariant: if $\mathpzc{A} \neq \{1,2\}$, we can shift any smooth sequence until it is not derivable. Over $\{1,2\}$, we have for example $S(\kappa_{2,1}) \notin \mathcal{C}^\infty$ since
	$$\mathcal{D}^4(\kappa_{2,1}) = \mathcal{D}^4(211212212211211...) = \mathcal{D}^3(1211212212...) = \mathcal{D}^2(1121121...) = \mathcal{D}(2121...) = 111... \notin \mathcal{C}^1.$$
\end{remark}

\subsubsection{Smooth words}

Over a binary alphabet $\mathpzc{A}$, any word $u \in \mathpzc{A}^*$ is canonically factorized as $a_1^{p_1} ... a_n^{p_n}$ where $a_i \in \mathpzc{A}$, $a_{i+1} = \overline{a_i}$ and $p_i \geq 1$, with the convention that $n = 0$ if $u = \varepsilon$.

\begin{definition}\label{Df}
	Over a binary alphabet $\mathpzc{A}$, we say that a word $u \in \mathpzc{A}^*$ is \textup{derivable} if its canonical factorization is $a_1^{p_1} ... a_n^{p_n}$ with $p_i \in \mathpzc{A}$ for all $i \in \llbracket 2,n-1 \rrbracket$ and ${p_1,p_n \in \llbracket 1,b \rrbracket}$. We write $\mathcal{C}^1_f$ the set of derivable words and we define the \textup{finite derivative} as the map
	\begin{align*}
		\mathcal{D}_f~:~~~~~~~~~ \mathcal{C}^1_f & \longrightarrow \mathpzc{A}^* \\
    	a_1^{p_1}...a_n^{p_n} & \longmapsto cut(p_1)~p_2~...~p_{n-1}~cut(p_n)
	\end{align*}
	where $cut(p) = \begin{cases} \varepsilon \textrm{ if } 1 \leq p \leq a \\ b \textrm{ if } a < p \leq b \end{cases}$, $\mathcal{D}_f(\varepsilon) = \varepsilon$ and ${\mathcal{D}_f({a_1}^{p_1}) = cut(p_1)}$.
\end{definition}

\begin{example}
	Over $\{1,2\}$, $\mathcal{C}^1_f$ is the set of words that do not have $111$ or $222$ as factors, and we have ${\mathcal{D}_f(\varepsilon) = \mathcal{D}_f(1) = \mathcal{D}_f(21) = \varepsilon}$ and $\mathcal{D}_f(2211) = \mathcal{D}_f(122112) = 22$. Over $\{1,3\}$, we have $\mathcal{D}_f(331113) = 33$.
\end{example}

\begin{definition}
	For $n \geq 2$, we define by induction the set of words that are derivable $n$ times
	\begin{equation*}
		\mathcal{C}_f^n := \left\{u \in \mathcal{C}^1_f ~\middle |~ \mathcal{D}_f(u) \in \mathcal{C}_f^{n-1} \right\}.
	\end{equation*}
	We then define the set of \textup{smooth words}, that is the set of words that are infinitely derivable
	\begin{equation*}
		\mathcal{C}_f^\infty := \displaystyle \bigcap_{n \geq 1} \mathcal{C}_f^n.
	\end{equation*}
\end{definition}

We notice that the map $\mathcal{D}_f$ is \textit{contracting}, i.e., $\left\lvert \mathcal{D}_f(u) \right\rvert < \lvert u \rvert$ if ${u \in \mathcal{C}^1_f \backslash \{\varepsilon\}}$. This implies that a word is smooth if and only if it can be derived down to $\varepsilon$.

\begin{definition}
	The \textup{height} of a smooth word $u$ is the first integer $n \geq 0$ such that $\mathcal{D}_f^n(u) = \varepsilon$.
\end{definition}

\begin{example}
	Over $\{1,2\}$, ${\mathcal{D}_f^4(221121221) = \mathcal{D}_f^3(22112) = \mathcal{D}_f^2(22) = \mathcal{D}_f(2) = \varepsilon}$ so $221121221 \in \mathcal{C}_f^\infty$ with height $4$, but $\mathcal{D}_f(12121) = 111 \notin \mathcal{C}^1_f$ so $12121 \notin \mathcal{C}_f^\infty$.
\end{example}

From the definitions of smooth sequences and smooth words, it is clear that $\mathscr{L}(\mathcal{C}^\infty) \subseteq \mathcal{C}_f^\infty$, which yields $p_x(n) \leq p_{\mathcal{C}^\infty}(n) \leq p_{\mathcal{C}_f^\infty}(n)$ for all $x \in \mathcal{C}^\infty$. In \cite{CH} we proved the equality of these languages.

\begin{theorem}[{\cite[Theorems 1.30]{CH}}]\label{CH1}
	Over a binary alphabet $\mathpzc{A} = \{a,b\}$, we have ${\mathscr{L}(\mathcal{C}^\infty) = \mathcal{C}_f^\infty}$, and in particular $p_{\mathcal{C}^\infty}(n) = p_{\mathcal{C}_f^\infty}(n)$.
\end{theorem}

In the same paper we studied the following conjecture.

\begin{conjecture}[{\cite{Sing02}}]\label{CfTheta}
	Over a binary alphabet $\mathpzc{A} = \{a,b\}$, we have $p_{\mathcal{C}_f^\infty}(n) = \Theta(n^\rho)$ where $\rho := \frac{\log(a+b)}{\log\left(\frac{a+b}{2}\right)}$.
\end{conjecture}

\begin{theorem}[{\cite[Theorems 1.31 and 1.32 (ii)]{CH}}]\label{CH2}
	Let $\mathpzc{A} = \{a,b\}$ be a binary alphabet and let $\rho := \frac{\log(a+b)}{\log\left(\frac{a+b}{2}\right)}$.
	
	\textup{(i)} We have $p_{\mathcal{C}_f^\infty}(n) = \Omega(n^\rho)$.
	
	\textup{(ii)} If $\mathpzc{A}$ is even, we have $p_{\mathcal{C}_f^\infty}(n) = \mathcal{O}(n^\rho)$.
\end{theorem}

\subsection{Properties of smooth sequences}

The study of smooth sequences differs greatly between mixed alphabets and even or odd alphabets, this shows on two levels:
\begin{enumerate}
	\item Smooth sequences seem to have different properties between the two cases,
	
	\item The study of smooth sequences has been more fruitful over  even and odd alphabets than over mixed alphabets.
\end{enumerate}

However, over any binary alphabet, it seems that all smooth sequences have similar properties. Let us detail the conjectures and results in the two cases.

\subsubsection{Smooth sequences over mixed alphabets}

Over mixed alphabets, the study of smooth sequences is governed by two main conjectures.

\begin{conjecture}\label{mixedCf}
	Let $\mathpzc{A} = \{a,b\}$ be a mixed alphabet.
	
	\textup{(i)} The language $\mathcal{C}_f^\infty$ is uniformly recurrent.
	
	\textup{(ii)} The language $\mathcal{C}_f^\infty$ is uniquely ergodic, and in particular it has uniform letter frequencies.
\end{conjecture}

One derives several sub-conjectures that can be considered independantly from one another.
	
\begin{conjecture}\label{subConjectures}
	Let $x$ be a smooth sequence over a mixed alphabet $\mathpzc{A} = \{a,b\}$.
	
	\textup{(i)} $\mathscr{L}(x) = \mathcal{C}_f^\infty$.
	
	\textup{(ii)} $x$ is unifomly recurrent.
	
	\textup{(iii)} $x$ has uniform frequencies, and it particular it has uniform letter frequencies.
	
	\textup{(iv)} $\mathscr{L}(x)$ is complement-invariant.
	
	\textup{(v)} $\mathscr{L}(x)$ is mirror-invariant.
	
	\textup{(vi)} $p_x(n) = \Theta(n^\rho)$ where $\rho := \frac{\log(a+b)}{\log\left(\frac{a+b}{2}\right)}$.
	
	\textup{(vii)} $x$ is not linearly recurrent.
	
	\textup{(viii)} $x$ is not morphic.
\end{conjecture}

These conjectures have been mostly looked at over $\{1,2\}$ to study the Oldenburger-Kolakoski sequence $\kappa$, the existing results are sparse and conditional.

\begin{proposition}[\cite{Rao,Nilsson}]\label{freq12}
	Over $\mathpzc{A} = \{1,2\}$, we have $\displaystyle\liminf_{u \in \mathcal{C}_f^\infty} \frac{\lvert u \rvert_2}{\lvert u \rvert} > \frac{1}{2} - 0.00008$.
\end{proposition}

\begin{proposition}[{\cite[Propositions 1 and 2]{Dekking81}}]\label{caracCfk}
	Let $\mathpzc{A} = \{1,2\}$.
	
	\textup{(i)} If $\mathscr{L}(\kappa)$ is complement-invariant, then $\kappa$ is recurrent.

	\textup{(ii)} $\mathscr{L}(\kappa)$ is complement-invariant if and only if $\mathscr{L}(\kappa) = \mathcal{C}_f^\infty$.
\end{proposition}

\begin{theorem}[{\cite[Theorem 1.40]{CH}}]\label{dichotomy}
	Let $\mathpzc{A} = \{1,2\}$.
	
	\textup{(i)} If $\mathscr{L}(\kappa) = \mathcal{C}_f^\infty$, then $\kappa$ is recurrent but not linearly recurrent and not morphic, and we have $p_{\kappa}(n) = \Omega\left(n^\rho\right)$ and $p_{\kappa}(n) = \mathcal{O}\left(n^{\rho+0.00023}\right)$ where ${\rho := \frac{\log(3)}{\log(3/2)} \approx 2.70951}$.
	
	\textup{(ii)} If $\mathscr{L}(\kappa) \subsetneq \mathcal{C}_f^\infty$ and $\kappa$ is uniformly recurrent, then $\kappa$ is linearly recurrent and morphic, and we have $p_{\kappa}(n) = \Theta(n)$.
\end{theorem}
	
\subsubsection{Smooth sequences over even or odd alphabets}

Over even or odd alphabets, the study of smooth sequences started with the fixed points $\kappa_{a,b}$. In \cite{Dekking79,Dekking01} Dekking noticed that $\kappa_{3,1}$ is morphic, then in \cite{Sing02,BS,Sing10} Sing showed the same for every fixed point $\kappa_{a,b}$ using other substitutions.
\begin{proposition}
	If $a$ and $b$ are two even (resp. odd) integers, then the smooth sequence $\kappa_{a,b}$ is morphic and generated by a primitive substitution.
\end{proposition}

Several properties immediatly follow.

\begin{corollary}\label{primitiveMorphicProperties}
	Let $a$ and $b$ be two even (resp. odd) integers.
	
	\textup{(i)} $\kappa_{a,b}$ is linearly recurrent.
	
	\textup{(ii)} $p_{\kappa_{a,b}} = \Theta(n)$.
	
	\textup{(iii)} $\kappa_{a,b}$ has uniform frequencies.
	
	\textup{(iv)} $\kappa_{a,b}$ has computable algebraic letter frequencies.
	
	\textup{(v)} $\mathscr{L}(\kappa_{a,b}) \subsetneq \mathcal{C}_f^\infty$.
\end{corollary}

\begin{example}\label{k31}
	We compare the constructions of Dekking and Sing for the sequence $\kappa_{3,1}$.
	\begin{itemize}
		\item Dekking's idea is to mark the parity of the letters with the alphabet $\{1,1',3,3'\}$. This yields $\kappa_{3,1} = \tau(\varsigma^\infty(3))$ where
	\begin{align*}
		\varsigma :~~~ & 1 \longmapsto 3, & & 1' \longmapsto 1', & & 3 \longmapsto 33'3, & & 3' \longmapsto 1'11' \\
		\tau :~~~ & 1 \longmapsto 1, & & 1' \longmapsto 1, & & 3 \longmapsto 3, & & 3' \longmapsto 3.
	\end{align*}
	
		\item Sing's idea is to code the sequence by three blocks of length two. This yields $\kappa_{3,1} = \pi(\sigma^\infty(3))$ where
		\begin{align*}
		\sigma :~~~ & A \longmapsto ABC, & & B \longmapsto AB, & & C \longmapsto B \\
		\pi :~~~ & A \longmapsto 33, & & B \longmapsto 31, & & C \longmapsto 11.
	\end{align*}
	\end{itemize}
	The substitution $\sigma$ has better properties than $\varsigma$: it is primitive whereas $\varsigma$ is not, and it is requires fewer letters. With a simple computation involving the incidence matrix of $\varsigma$ or $\sigma$, the frequency of the letter $3$ in $\kappa_{3,1}$ is approximately 0.6028.
\end{example}

Beyond the fixed points $\kappa_{a,b}$, the papers \cite{BJP} and \cite{JMPR} study all smooth sequences over an even or odd alphabet. In the first one, the authors obtain results about recurrence, mirror- and complement-invariance.

\begin{theorem}[\cite{BJP}]\label{BJP}
	Let $\mathpzc{A}$ be an even or odd alphabet.
	
	\textup{(i)} Every smooth sequence is recurrent.
	
	\textup{(ii)} If $\mathpzc{A}$ is even, then the language of the minimal smooth sequence with respect to the lexicographic order is not mirror-invariant. If $\mathpzc{A}$ is odd, then the language of every smooth sequence is mirror-invariant.
	
	\textup{(iii)} The language of the minimal smooth sequence with respect to the lexicographic order is not complement-invariant.
\end{theorem}
	
In particular, the sequences such that $\mathscr{L}(x)$ is not mirror-invariant or complement-invariant satisfy $\mathscr{L}(x) \subsetneq \mathcal{C}_f^\infty$. In the same paper, the authors manage to compute the letter frequency of a smooth sequence that is not the fixed point $\kappa_{a,b}$.
	
\begin{proposition}[\cite{BJP}]\label{minLexFreq}
	Over the odd alphabet $\mathpzc{A} = \{1,b\}$, the frequency of $b$'s in the minimal smooth sequence with respect to the lexicographic order is $\frac{1}{\sqrt{2b-1}+1}$.
\end{proposition}
	
In \cite{JMPR} the authors prove stronger results over $\{1,3\}$. Extending Dekking's idea of marking the parity of letters (cf. \Cref{k31}), they give an $S$-adic representation of every (bi-infinite) smooth sequence. They derive the following results.
	
\begin{theorem}[\cite{JMPR}]\label{JMPR}
	Let $\mathpzc{A} = \{1,3\}$.
	
	\textup{(i)} Every (bi-infinite) smooth sequence has uniform frequencies.
	
	\textup{(ii)} Almost every (bi-infinite) smooth sequence is uniformly recurrent.
\end{theorem}

They also show that these smooth sequences have diverse letter frequencies.

\subsection{New results}

In this paper we provide results about the $S$-adic structure, the complement-invariance, the recurrence, the factor complexity and the frequencies of smooth sequences over even and odd alphabets. In \Cref{summary} we summarize the main properties we prove based on a partition of smooth sequences (cf. \Cref{propertyM,propertyM2}). 

\begin{figure}[H]
	\begin{center}
	\begin{tabular}{|c|c|c|c|c|c|}
		\hline
  		$x \in \mathcal{C}^\infty$ & $\substack{\textup{Complement-invariant} \\ \textup{(\Cref{notInvariant} (i))}}$ & $\substack{\textup{Uniformly recurrent} \\ \textup{(\Cref{URLR} (i))}}$ & $\substack{\textup{Linearly recurrent} \\ \textup{(\Cref{URLR} (ii))}}$ & $\substack{\textup{Linear factor complexity} \\ \textup{(\Cref{linearCorollary})}}$ & $\substack{\textup{Uniform frequencies} \\ \textup{(\Cref{frequenciesCorollary})}}$ \\
		\hline
		\hline
		$x \in \mathscr{M}$ & \xmark & \cmark & \cmark & \cmark & \cmark \\
		\hline
		$x \in \mathscr{M}' \backslash \mathscr{M}$ & \xmark & \cmark & \xmark & \cmark & \cmark \\
		\hline
		$x \notin \mathscr{M}'$ & \xmark & \cmark & \xmark & \textbf{?} & \textbf{?} \\
		\hline
	\end{tabular}
	\end{center}
	\caption{The properties of smooth sequences over even and odd alphabets.}
	\label{summary}
\end{figure}

\subsubsection{$S$-adic structure of smooth sequences}

In \cite{Dekking01}, Dekking showed that a central object to understand smooth sequences is the following map.

\begin{definition}\label{trace}
	Over $\mathpzc{A} = \{a,b\}$, the \textup{trace} is the map
	\begin{align*}
		\mathcal{T} : \mathcal{C}^\infty & \longrightarrow \mathpzc{A}^\mathbb{N} \\
    	x & \longmapsto \displaystyle \left(\mathcal{D}^n(x)_0\right)_{n \geq 0}.
	\end{align*}
\end{definition}

In \Cref{structure} we use the substitutions introduced by Sing in \cite{BS} (illustrated in \Cref{k31}, properly introduced in \Cref{substitutions}) to exhibit the $S$-adic structure of smooth sequences.
\newpage
\begin{proposition}\label{smoothSadicSystem}
	Over an even or odd alphabet $\mathpzc{A} = \{a,b\}$, there exists a finite alphabet $\Sigma$, four substitutions $\sigma_a, \sigma_b : \Sigma \rightarrow \Sigma^+$ and $\pi_a, \pi_b : \Sigma \rightarrow \mathpzc{A}^2$ and a letter $A \in \Sigma$ such that the set of smooth sequences $\mathcal{C}^\infty$ is the $S$-adic system
	$$\left\{\displaystyle\lim_{n \rightarrow \infty} \pi_{t_0} \sigma_{t_1} \sigma_{t_2} ... \sigma_{t_n}(A^\omega) ~\middle|~ t \in \mathpzc{A}^\mathbb{N}\right\}.$$
	Moreover, the directive sequence of a smooth sequence $x$ is $t = \mathcal{T}(x)$.
\end{proposition}

\subsubsection{Complement-invariance of smooth sequences}

In \Cref{occurrences} we show that the language of any smooth sequence over an even or odd alphabet is anti-complement-variant.

\begin{proposition}\label{antiInvariance}
	Let $\mathpzc{A} = \{a,b\}$ be an even or odd alphabet. Then there exists a finite set $\mathfrak{S} \subset \mathpzc{A}^*$ such that, if $u \in \mathscr{L}(x)$ and $\mathcal{D}_f(u) \notin \mathfrak{S}$, then $\overline{u} \notin \mathscr{L}(x)$.
\end{proposition}

We deduce that all smooth sequences over even or odd alphabets are far from containing all smooth words.

\begin{corollary}\label{notInvariant}
	Let $x$ be a smooth sequence over an even or odd alphabet.
	
	\textup{(i)} The language $\mathscr{L}(x)$ is not complement-invariant.
	
	\textup{(ii)} We have $\mathscr{L}(x) \subsetneq \mathcal{C}_f^\infty$.
\end{corollary}

\subsubsection{Recurrence of smooth sequences}

In \Cref{recurrence} we determine how recurrent smooth sequences over even or odd alphabets are, based on the following definition.

\begin{definition}\label{propertyM}
	Let $x$ be a smooth sequence over an odd alphabet $\mathpzc{A} = \{1,b\}$. 
	\begin{itemize}
		\item We say that $x$ has property $\mathfrak{P}1$ if there exists $n \geq 1$ such that $(1b)^n \nsqsubset \mathcal{T}(x)$.
		
		\item We say that $x$ has property $\mathfrak{P}2$ if there exists $n \geq 0$ such that $S^n(\mathcal{T}(x)) = (1b)^\omega$.
	\end{itemize}
	We define the class of smooth sequences
	$$\mathscr{M} := \left\{x \in \mathcal{C}^\infty \textit{ over } \mathpzc{A} = \{a,b\} ~\middle|~ \mathpzc{A} \textit{ is even or odd, and } x \textit{ has property } \mathfrak{P}1 \textit{ or } \mathfrak{P}2 \textit{ if } a = 1 \right\}.$$
\end{definition}

\begin{theorem}\label{URLR}
	Let $x$ be a smooth sequence over an even or odd alphabet.
	
	\textup{(i)} $x$ is uniformly recurrent.
	
	\textup{(ii)} $x$ is linearly recurrent if and only if $x \in \mathscr{M}$.
\end{theorem}

The proofs of (i) and (ii) are independant: (i) relies on an induction on the height of factors and on technical lemmas while (ii) relies on the main result of \cite{Durand03}. In particular, this improves \Cref{JMPR} \textup{(ii)}.

\subsubsection{Factor complexity of smooth sequences}

A corollary of \Cref{URLR} \textup{(ii)} is that the smooth sequences in $\mathscr{M}$ have linear factor complexity. In \Cref{complexity} we go further by precisely  describing the bispecial factors of smooth sequences. We deduce in many cases a strong homogeneity of smooth sequences in terms of factor complexity.

\begin{theorem}\label{sameComplexity}
	Let $\mathpzc{A} = \{a,b\}$ be an even alphabet, or an odd alphabet such that $b = a+2$.
	
	\textup{(i)} If $a > 1$ (resp. if $a = 1$), then all smooth sequences (resp. all smooth sequences that do not have property $\mathfrak{P}2$) have the same factor complexity $p(n)$.
	
	\textup{(ii)} If $b = a+2$, then we have $p(n) = \begin{cases} 1 \textrm{ if } n = 0 \\ 2n \textrm{ if } n \in \llbracket 1,a+1 \rrbracket \\ 2n+2 \textrm{ if } n \geq a+2 \end{cases}$.
	
	\textup{(iii)} If $\mathpzc{A}$ is even and $b \in \left\llbracket a+4, \frac{1}{2}a^2 + 2a \right\rrbracket$, then we have $2 < \displaystyle\liminf_{n \rightarrow \infty} \frac{p(n)}{n} < \displaystyle\limsup_{n \rightarrow \infty} \frac{p(n)}{n} < 3$.
	
	\textup{(iv)} If $\mathpzc{A}$ is even and $b \in \left\llbracket \frac{1}{2}a^2 + 2a + 2, \frac{1}{2}a^2 + 3a \right\rrbracket$, then we have $2 < \displaystyle\liminf_{n \rightarrow \infty} \frac{p(n)}{n} < 3 < \displaystyle\limsup_{n \rightarrow \infty} \frac{p(n)}{n} < 4$.
	
	\textup{(v)} If $\mathpzc{A}$ is even and $b \in \left\llbracket \frac{1}{2}a^2 + 3a + 2, +\infty \right)$, then we have $3 < \displaystyle\liminf_{n \rightarrow \infty} \frac{p(n)}{n} < \displaystyle\limsup_{n \rightarrow \infty} \frac{p(n)}{n} < 4$.
\end{theorem}

We deduce linear factor complexity for a class of smooth sequences larger than $\mathscr{M}$.

\begin{definition}\label{propertyM2}
	We define the class of smooth sequences $\mathscr{M}' := \mathscr{M} \cup \left\{x \in \mathcal{C}^\infty \textit{ over } \{1,3\} \right\}$.
\end{definition}

\begin{corollary}\label{linearCorollary}
	A smooth sequence in $\mathscr{M}'$ has linear factor complexity.
\end{corollary}

\begin{remark}
	 We deduce another proof of \Cref{notInvariant} \textup{(ii)} for smooth sequences in $\mathscr{M}'$: \Cref{linearCorollary} states that these sequences have linear factor complexity whereas \Cref{CH2} \textup{(i)} implies that the complexity of $\mathcal{C}_f^\infty$ grows more than linearly.
\end{remark}

\subsubsection{Frequencies in smooth sequences}

In \Cref{frequencies} we show that the smooth sequences in $\mathscr{M}'$ also have uniform frequencies.

\begin{corollary}\label{frequenciesCorollary}
	A smooth sequence in $\mathscr{M}'$ has uniform frequencies.
\end{corollary}

This improves \Cref{JMPR} \textup{(i)}, with a different proof. Then we show how to compute the letter frequencies in some smooth sequences.

\begin{theorem}\label{smoothFreq}
	Let $x$ be a smooth sequence over an even or odd alphabet $\mathpzc{A} = \{a,b\}$.
	
	\textup{(i)} If $\mathpzc{A}$ is even, then the frequency of $a$'s and $b$'s in $x$ are $1/2$.
	
	\textup{(ii)} If $\mathpzc{A}$ is odd, suppose that $\mathcal{T}(x) = c u^\omega$ for some $c \in \mathpzc{A}$ and some $u \in \mathpzc{A}^n$ ($n \geq 1$), and that $u$ is not of the form $(1b)^n$ or $(b1)^n$ if $a = 1$. Then the frequency of $c$'s in $x$ is the algebraic number $\alpha + \frac{1}{2}\beta$ where $\left(\begin{matrix}\alpha \\ \beta \\ \gamma \end{matrix}\right)$ is the unitary right Perron eigenvector of the primitive matrix $\displaystyle \prod_{i=1}^n M_{u_i}$ and
$M_a := \left(\begin{matrix} \frac{a-1}{2} & \frac{a-1}{2} & \frac{b-1}{2} \\ 1 & 1 & 1 \\ \frac{a-1}{2} & \frac{b-1}{2} & \frac{b-1}{2} \end{matrix}\right)$ and $M_b := \left(\begin{matrix} \frac{b-1}{2} & \frac{b-1}{2} & \frac{a-1}{2} \\ 1 & 1 & 1 \\ \frac{b-1}{2} & \frac{a-1}{2} & \frac{a-1}{2} \end{matrix}\right)$.
	
	\textup{(iii)} If $\mathpzc{A} = \{1,b\}$ is odd and $\mathcal{T}(x) \in \{(1b)^\omega,(b1)^\omega\}$, then the frequency of $b$'s in $x$ is $\frac{1}{\sqrt{2b-1}+1}$.
\end{theorem}

\begin{remark}
	\Cref{smoothFreq} unifies the existing results regarding the value of letter frequencies in smooth sequences over even or odd alphabets.
	\begin{itemize}
		\item \Cref{smoothFreq} \textup{(i)} was known to have an easy proof but it did not have a proper one yet.
		
		\item Over $\{1,3\}$, we recover the frequency of 3's in $\kappa_{3,1}$ ($\approx 0.6028$, cf. \Cref{k31}). More generally, we find smooth sequences with various frequencies, as detailed in \cite{JMPR}. For example $\mathcal{E}((331)^\omega)$ (resp. $\mathcal{E}((31131)^\omega)$) has a frequency of $3$'s approximately equal to $0.6624$ (to $0.68$).
	
	\item Over an odd alphabet $\mathpzc{A} = \{1,b\}$, in \cite{BJP} the authors show that the minimum smooth sequence $x$ with respect to the lexicographic order satisfies $\mathcal{T}(x) = (1b)^\omega$. Therefore \Cref{smoothFreq} \textup{(iii)} is a reformulation of \Cref{minLexFreq} with a simpler proof.
	\end{itemize}
\end{remark}

\section{How to generate smooth sequences}\label{structure}

In this section we characterize the smooth sequences over any binary alphabet based on the trace map $\mathcal{T} : x \mapsto \left(\mathcal{D}^n(x)_ 0\right)_{n \geq 0}$, this is \Cref{conjugacy}. This result implies that there are infinitely many smooth sequences over any binary alphabet, and it allows to compute arbitrarily large prefixes of any smooth sequence from its trace.

Then we show that every smooth sequence over an even or odd alphabet is the image by a substitution of an $S$-adic sequence generated by two substitutions, this is \Cref{smoothSadic}. Combined with \Cref{conjugacy}, it yields \Cref{smoothSadicSystem}.

Lastly, among the smooth sequences over even or odd alphabets that we now know are $S$-adic, we focus on the sequences that are morphic. In particular, in \Cref{1bb1morphic} we study two smooth sequences that behave differently from all the other smooth sequences and will have to be handled carefully in the rest of the paper.

\subsection{The Kolakoski shift}

We begin by constructing the reciprocal of the trace $\mathcal{T} : x \mapsto \left(\mathcal{D}^n(x)_ 0\right)_{n \geq 0}$.

\begin{definition}\label{primitive}
	Over $\mathpzc{A} = \{a,b\}$, the \textup{primitives} are the continuous maps
	\begin{align*}
		\mathcal{P}_a~:~\mathpzc{A}^\mathbb{N} & \longrightarrow \mathcal{C}^1 & \mathcal{P}_b~:~ \mathpzc{A}^\mathbb{N} & \longrightarrow \mathcal{C}^1 \\
    	x & \longmapsto a^{x_0}b^{x_1}a^{x_2}b^{x_3}... & x & \longmapsto \overline{\mathcal{P}_a(x)}
	\end{align*}
\end{definition}

The primitives $\mathcal{P}_a$ and $\mathcal{P}_b$ are partial reciprocals of the derivative $\mathcal{D}$, in the sense that
$$\mathcal{P}_{x_0} \circ \mathcal{D} (x) = \mathcal{D} \circ \mathcal{P}_a (x) = \mathcal{D} \circ \mathcal{P}_b (x) = x$$
for all $x \in \mathpzc{A}^\mathbb{N}$. This easily generalizes.

\begin{lemma}\label{PPDD}
	Let $x \in \mathcal{C}^\infty$ over $\mathpzc{A} = \{a,b\}$ and let $t = \mathcal{T}(x)$. Then, for all $n \geq 1$, we have
	$$x = \mathcal{P}_{t_0} \mathcal{P}_{t_1} ... \mathcal{P}_{t_{n-1}} (\mathcal{D}^n(x)).$$
\end{lemma}

Now we show that composing primitives produces smooth sequences.

\begin{lemma}\label{smoothConvergence}
	Let $\mathpzc{A} = \{a,b\}$ be a binary alphabet. If $t \in \mathpzc{A}^\mathbb{N}$, then the sequence $\left(\mathcal{P}_{t_0} \mathcal{P}_{t_1} ... \mathcal{P}_{t_{n-1}}(t_n^\omega) \right)_{n \geq 0}$ converges in $\mathpzc{A}^\mathbb{N}$ and $\displaystyle\lim_{n \rightarrow \infty} \mathcal{P}_{t_0} \mathcal{P}_{t_1} ... \mathcal{P}_{t_{n-1}}(t_n^\omega) \in \mathcal{C}^\infty$.
\end{lemma}
\begin{proof}
	The first step is to find prefixes of the infinite word $\mathcal{P}_{t_0} \mathcal{P}_{t_1} ... \mathcal{P}_{t_{n-1}}(t_n^\omega)$ for each $n \geq 0$. To do so, consider the maps
	\begin{align*}
		\mathfrak{p}_a~:~~~~~~~~\mathpzc{A}^* & \longrightarrow \mathpzc{A}^+ & \mathfrak{p}_b~: \mathpzc{A}^* & \longrightarrow \mathpzc{A}^+ \\
    	u_1...u_n & \longmapsto \begin{cases} a^{u_1}b^{u_2}...b^{u_n} a^a \textrm{ if } n \textrm{ is even,} \\
   		a^{u_1}b^{u_2}...a^{u_n} b^a \textrm{ otherwise.} \end{cases} & u & \longmapsto \overline{\mathfrak{p}_a(u)}
	\end{align*}
	We observe that, if $c \in \mathpzc{A}$ and $u \sqsubset_p x \in \mathpzc{A}^\mathbb{N}$, then $\mathfrak{p}_c(u) \sqsubset_p \mathcal{P}_c(x)$. We directly deduce that
	\begin{equation}\label{prefixpP}
		\mathfrak{p}_{t_0} \mathfrak{p}_{t_1} ... \mathfrak{p}_{t_{n-1}}(t_n) \sqsubset_p \mathcal{P}_{t_0} \mathcal{P}_{t_1} ... \mathcal{P}_{t_{n-1}}(t_n^\omega)
	\end{equation}
	for all $n \geq 0$.
	
	The second step is to show that the sequence of prefixes we exhibited converges in $\mathpzc{A}^\mathbb{N}$. To do so, we show that, if $c \in \mathpzc{A}$ and $u <_p v \in \mathpzc{A}^*$, then
	\begin{equation}\label{prefixpp}
		\mathfrak{p}_c(u) <_p \mathfrak{p}_c(v).
	\end{equation}
	By writting $v = v_1...v_n$ with $n \geq 0$, there exists $m \in \llbracket 0,n \llbracket$ such that ${u = v_1...v_m}$ and we have
	$$\mathfrak{p}_c(u) = a_1^{v_1} ... a_m^{v_m} a_{m+1}^a <_p a_1^{v_1} ... a_m^{v_m} a_{m+1}^{v_{m+1}} ... a_n^{v_n} a_{n+1}^a = \mathfrak{p}_c(v)$$
	where $a_1 = c$ and $a_{i+1} = \overline{a_i}$. This concludes the proof of \Cref{prefixpp}. Then, with \Cref{prefixpp}, a quick induction on $n$ shows that
	\begin{equation}\label{prefixppp}
		\mathfrak{p}_{t_0} \mathfrak{p}_{t_1} ... \mathfrak{p}_{t_{n-1}}(t_n) <_p \mathfrak{p}_{t_0} \mathfrak{p}_{t_1} ... \mathfrak{p}_{t_n}(t_{n+1})
	\end{equation}
	for all $n \geq 0$. Now \Cref{prefixppp} ensures that the sequence $\left(\mathfrak{p}_{t_0} \mathfrak{p}_{t_1} ... \mathfrak{p}_{t_{n-1}}(t_n) \right)_{n \geq 0}$ converges in $\mathpzc{A}^\mathbb{N}$. With \Cref{prefixpP}, this implies that the sequence $\left(\mathcal{P}_{t_0} \mathcal{P}_{t_1} ... \mathcal{P}_{t_{n-1}}(t_n^\omega)\right)_{n \geq 0}$ converges in $\mathpzc{A}^\mathbb{N}$.
	
	The last step is to show that the word $x := \displaystyle\lim_{n \rightarrow \infty} \mathcal{P}_{t_0} \mathcal{P}_{t_1} ... \mathcal{P}_{t_{n-1}}(t_n^\omega) \in \mathpzc{A}^\mathbb{N}$ is smooth. For all $m \geq 1$, the previous result applied to the sequence $(t_n)_{n \geq m}$ yields the convergence of the sequence $\left(\mathcal{P}_{t_m} \mathcal{P}_{t_{m+1}} ... \mathcal{P}_{t_{m+n-1}}({t_{m+n}}^\omega)\right)_{n \geq 0}$. Then, by continuity of $\mathcal{P}_a$ and $\mathcal{P}_b$, we obtain
	\begin{equation}\label{PlimP}
		x = \mathcal{P}_{t_0} ... \mathcal{P}_{t_{m-1}} \left(\displaystyle\lim_{n \rightarrow \infty} \mathcal{P}_{t_m} \mathcal{P}_{t_{m+1}} ... \mathcal{P}_{t_{m+n-1}}(t_{m+n}^\omega)\right) \in \mathcal{C}^m
	\end{equation}
	for all $m \geq 1$, which means that $x \in \mathcal{C}^\infty$.
\end{proof}

\Cref{smoothConvergence} yields the following definition.

\begin{definition}
	Over $\mathpzc{A} = \{a,b\}$, the \textup{expansion} is the map
	\begin{align*}
		\mathcal{E} : \mathpzc{A}^\mathbb{N} & \longrightarrow \mathcal{C}^\infty \\
    	t & \longmapsto \displaystyle\lim_{n \rightarrow \infty} \mathcal{P}_{t_0} \mathcal{P}_{t_1} ... \mathcal{P}_{t_{n-1}}(t_n^\omega).
	\end{align*}
\end{definition}

\begin{remark}
	Over $\mathpzc{A} = \{a,b\}$, if $c \in \mathpzc{A}$ and $t \in \mathpzc{A}^\mathbb{N}$, we have $\mathcal{E}(ct) = \overline{\mathcal{E}(\overline{c}t)}$.
\end{remark}

We show that the map $\mathcal{E} : \mathpzc{A}^\mathbb{N} \rightarrow \mathcal{C}^\infty$ is the reciprocal of the map $\mathcal{T} : \mathcal{C}^\infty\rightarrow \mathpzc{A}^\mathbb{N}$, which conjugates the action of the derivation $\mathcal{D}$ on $\mathcal{C}^\infty$ to the action of the shift $S$ on $\mathpzc{A}^\mathbb{N}$.

\begin{proposition}\label{conjugacy}
	Let $\mathpzc{A} = \{a,b\}$ be a binary alphabet.

	\textup{(i)} The map $\mathcal{T} : \mathcal{C}^\infty \rightarrow \mathpzc{A}^\mathbb{N}$ is bijective and its reciprocal is the map $\mathcal{E} : \mathpzc{A}^\mathbb{N} \rightarrow \mathcal{C}^\infty$.
	
	\textup{(ii)} $(\mathcal{C}^\infty, \mathcal{D})$ is conjugated to the full-shift $(\mathpzc{A}^\mathbb{N}, S)$, which provides the commuting diagram
		\[\begin{tikzcd}
			\mathcal{C}^\infty \arrow{r}{\mathcal{D}} \arrow[swap, d, shift left, "\mathcal{T}"'] \arrow[swap, d, leftarrow, shift right, "\mathcal{E}"] & \mathcal{C}^\infty \arrow[swap, d, shift left, "\mathcal{T}"'] \arrow[swap, d, leftarrow, shift right, "\mathcal{E}"] \\%
\mathpzc{A}^\mathbb{N} \arrow{r}{S}& \mathpzc{A}^\mathbb{N}
		\end{tikzcd}\]
\end{proposition}
\begin{proof}
	(i) Let $x \in \mathcal{C}^\infty$ and let $t = \mathcal{T}(x)$. Recall from the proof of \Cref{smoothConvergence} that $\mathcal{E}(t) = \displaystyle\lim_{n \rightarrow \infty} \mathfrak{p}_{t_0} \mathfrak{p}_{t_1} ... \mathfrak{p}_{t_{n-1}} (t_n)$, and that, if $c \in \mathpzc{A}$ and $u \sqsubset_p x \in \mathpzc{A}^\mathbb{N}$, then $\mathfrak{p}_c(u) \sqsubset_p \mathcal{P}_c(x)$. For all $n \geq 0$, since $t_n = \mathcal{D}^n(x)_0 \sqsubset_p \mathcal{D}^n(x)$, we directly deduce that 
	$$\mathfrak{p}_{t_0} \mathfrak{p}_{t_1} ... \mathfrak{p}_{t_{n-1}} (t_n) \sqsubset_p \mathcal{P}_{t_0} \mathcal{P}_{t_1} ... \mathcal{P}_{t_{n-1}}(\mathcal{D}^n(x)) = x$$
	thanks to \Cref{PPDD}. This implies that
	$$\mathcal{E} \circ \mathcal{T} (x) = \mathcal{E}(t) = \displaystyle\lim_{n \rightarrow \infty} \mathfrak{p}_{t_0} \mathfrak{p}_{t_1} ... \mathfrak{p}_{t_{n-1}} (t_n) = x.$$
		
	Reciprocally, let $t \in \mathpzc{A}^\mathbb{N}$ and let $x = \mathcal{E}(t) = \displaystyle\lim_{n \rightarrow \infty} \mathcal{P}_{t_0} \mathcal{P}_{t_1} ... \mathcal{P}_{t_{n-1}}(t_n^\omega)$. For all $m \geq 1$, \Cref{PlimP} yields $\mathcal{T}(x)_m = \mathcal{D}^m(x)_0 = t_m$, and we also have $\mathcal{T}(x)_0 = x_0 = t_0$, therefore $t = \mathcal{T} \circ \mathcal{E}(t)$.
	
	(ii) Let $x \in \mathcal{C}^\infty$ and let $t = \mathcal{T}(x)$. Then we have
	\begin{align*}
		\mathcal{D}(x) & = \mathcal{D} \circ \mathcal{E}(t) ~~~~~~~~~~~~\textup{thanks to (i)} \\
		& = \mathcal{D} \circ \mathcal{P}_{t_0} \left(\displaystyle\lim_{n \rightarrow \infty} \mathcal{P}_{t_1} \mathcal{P}_{t_2} ... \mathcal{P}_{t_n}(t_{n+1}^\omega)\right) ~~~~~~\textrm{thanks to \Cref{PlimP}} \\
		& = \displaystyle\lim_{n \rightarrow \infty} \mathcal{P}_{t_1} \mathcal{P}_{t_2} ... \mathcal{P}_{t_n}(t_{n+1}^\omega) \\
		& = \mathcal{E} \left(\left(t_n\right)_{n \geq 1}\right) \\
		& = \mathcal{E} \circ S \circ \mathcal{T} (x).
	\end{align*}
\end{proof}

With \Cref{conjugacy} we can now express any smooth sequence as $\mathcal{E}(t)$ for some $t \in \mathpzc{A}^\mathbb{N}$. This is a contructive expression because we exhibited in the proof of \Cref{smoothConvergence} the computable sequence of prefixes $\left(\mathfrak{p}_{t_0} \mathfrak{p}_{t_1} ... \mathfrak{p}_{t_{n-1}}(t_n)\right)_{n \geq 0}$ that satisfies $\displaystyle \lim_{n \rightarrow \infty} \mathfrak{p}_{t_0} \mathfrak{p}_{t_1} ... \mathfrak{p}_{t_{n-1}}(t_n) = \mathcal{E}(t)$.

\begin{example}\label{computeE}
	Over $\{1,2\}$, starting from the first terms of a trace $t = 22112...$, we compute prefixes of the sequence $\mathcal{E}(t)$ like so:
	\begin{align*}
		\mathfrak{p}_2(2) = &~ 221 \\
		\mathfrak{p}_2 \mathfrak{p}_2(1) = \mathfrak{p}_2 (21) = &~ 2212 \\
		\mathfrak{p}_2 \mathfrak{p}_2 \mathfrak{p}_1(1) = \mathfrak{p}_2 \mathfrak{p}_2 (12) = \mathfrak{p}_2 (2112) = &~ 2212112 \\
		\mathfrak{p}_2 \mathfrak{p}_2 \mathfrak{p}_1 \mathfrak{p}_1(2) = \mathfrak{p}_2 \mathfrak{p}_2 \mathfrak{p}_1(112) = \mathfrak{p}_2 \mathfrak{p}_2 (12112) = \mathfrak{p}_2 (21121221) = &~ 2212112112212
	\end{align*}
\end{example}

\begin{example}\label{Eab}
	Over $\{1,2\}$, \Cref{conjugacy} \textup{(ii)} yields that the smooth sequences
	\begin{align*}
		\mathcal{E}\left((21)^\omega\right) =	21221211221221121121221121122121121122122121121221121121221221121221211211... \\
		\mathcal{E}\left((12)^\omega\right) = 11211221221211221221121221211211221211212211211212212112212211211212212211...
	\end{align*}
	are the derivative of each other.
\end{example}

\begin{remark}\label{fixedPoints}
	It follows from \Cref{conjugacy} that, over any binary alphabet $\mathpzc{A} = \{a,b\}$, the fixed points of $\mathcal{D}^n$ ($n \geq 1$) are the smooth sequences of the form $\mathcal{E}(u^\omega)$ where $u \in \mathpzc{A}^n$. In particular, the two fixed points of $\mathcal{D}$ introduced in \Cref{kab} are $\kappa_{a,b} = \mathcal{E}(a^\omega)$ and $\kappa_{b,a} = \mathcal{E}(b^\omega)$.
\end{remark}

\subsection{An $S$-adic representation}

We first introduce the substitutions that will produce the desired $S$-adic system.

\begin{definition}\label{substitutions}
	Over an even or odd alphabet $\mathpzc{A} = \{a,b\}$, we define an alphabet $\Sigma$ and four substitutions $\sigma_a, \sigma_b : \Sigma \rightarrow \Sigma^+$ and $\pi_a, \pi_b : \Sigma \rightarrow \mathpzc{A}^+$ as follows.
	
	(i) If $\mathpzc{A}$ is even, define the substitutions over the alphabet $\Sigma := \{A,B\}$
	\begin{align*}
		\sigma_a :~~~ & A \longmapsto A^\frac{a}{2}~B^\frac{a}{2}, & & B \longmapsto A^\frac{b}{2}~B^\frac{b}{2}, \\
		\sigma_b :~~~ & A \longmapsto A^\frac{b}{2}~B^\frac{b}{2}, & & B \longmapsto A^\frac{a}{2}~B^\frac{a}{2}, \\
		\pi_a :~~~ & A \longmapsto aa, & & B \longmapsto bb, \\
		\pi_b :~~~ & A \longmapsto bb, & & B \longmapsto aa.
	\end{align*}
	(ii) If $\mathpzc{A}$ is odd, define the substitutions over the alphabet $\Sigma := \{A,B,C\}$
	\begin{align*}
		\sigma_a :~~~ & A \longmapsto A^\frac{a-1}{2}~B~C^\frac{a-1}{2}, & & B \longmapsto A^\frac{a-1}{2}~B~C^\frac{b-1}{2}, & & C \longmapsto A^\frac{b-1}{2}~B~C^\frac{b-1}{2}, \\
		\sigma_b :~~~ & A \longmapsto A^\frac{b-1}{2}~B~C^\frac{b-1}{2}, & & B \longmapsto A^\frac{b-1}{2}~B~C^\frac{a-1}{2}, & & C \longmapsto A^\frac{a-1}{2}~B~C^\frac{a-1}{2} \\
		\pi_a :~~~ & A \longmapsto aa, & & B \longmapsto ab, & & C \longmapsto bb \\
		\pi_b :~~~ & A \longmapsto bb, & & B \longmapsto ba, & & C \longmapsto aa.
	\end{align*}
\end{definition}

Over an even or odd alphabet $\mathpzc{A} = \{a,b\}$, the maps $\pi_a : \Sigma^\mathbb{N} \rightarrow \mathpzc{A}^\mathbb{N}$ and $\pi_b : \Sigma^\mathbb{N} \rightarrow \mathpzc{A}^\mathbb{N}$ are clearly injective so we can define the maps $\pi_a^{-1} : \pi_a(\Sigma^\mathbb{N}) \rightarrow \Sigma^\mathbb{N}$ and $\pi_b^{-1} : \pi_b(\Sigma^\mathbb{N}) \rightarrow \Sigma^\mathbb{N}$. In particular, for all $x \in \Sigma^\mathbb{N}$ and all $c \in \mathpzc{A}$, we have $\pi_c^{-1} \pi_c(x) = x$.

\begin{proposition}\label{smoothSadic}
	Let $x$ be a smooth sequence over an even or odd alphabet $\mathpzc{A} = \{a,b\}$ and let $t = \mathcal{T}(x)$. Then we have
	$$x = \displaystyle\lim_{n \rightarrow \infty} \pi_{t_0} \sigma_{t_1} \sigma_{t_2} ... \sigma_{t_n}(A^\omega).$$
\end{proposition}
\begin{proof}
	The proof relies on the following lemma: if $y = d^\omega$ for some $d \in \mathpzc{A}$ or $y \in \mathcal{C}^1$, then $y \in \pi_{y_0}(\Sigma^\mathbb{N})$ and, for all $c \in \mathpzc{A}$, we have
	\begin{equation}\label{SP}
		\mathcal{P}_c(y) = \pi_c \sigma_{y_0} \pi_{y_0}^{-1}(y).
	\end{equation}
	Suppose that $\mathpzc{A}$ is even. If $y = d^\omega$, we have $y = \pi_d(A^\omega) \in \pi_d(\Sigma^\mathbb{N})$ and 
	$$\pi_c \sigma_d \pi_d^{-1}(y) = \pi_c \sigma_d (A^\omega) = \pi_c \left( \left(A^\frac{d}{2} B^\frac{d}{2}\right)^\omega \right) = \left(c^d \overline{c}^d\right)^\omega = \mathcal{P}_c(x).$$
	If $y \in \mathcal{C}^1$, by setting $(p_n)_{n \geq 0} := \left(\frac{\mathcal{D}(y)_n}{2}\right)_{n \geq 0} \in \mathbb{N}^\mathbb{N}$, we get
	$$y = \mathcal{P}_{y_0}(\mathcal{D}(y)) = {y_0}^{2p_0}~\overline{y_0}^{2p_1}~{y_0}^{2p_2}~\overline{y_0}^{2p_3}~... = \pi_{y_0}(A^{p_0}~B^{p_1}~A^{p_2}~B^{p_3}~...) \in \pi_{y_0}(\Sigma^\mathbb{N})$$
	and
	\begin{align*}
		\pi_c \sigma_{y_0} \pi_{y_0}^{-1}(y) & = \pi_c \sigma_{y_0} \left(A^{p_0}~B^{p_1}~A^{p_2}~B^{p_3}~...\right) \\
		& = \pi_c \left( \left(A^\frac{y_0}{2}~B^\frac{y_0}{2}\right)^{p_0}~\left(A^\frac{\overline{y_0}}{2}~B^\frac{\overline{y_0}}{2}\right)^{p_1}~\left(A^\frac{y_0}{2}~B^\frac{y_0}{2}\right)^{p_2}~\left(A^\frac{\overline{y_0}}{2}~B^\frac{\overline{y_0}}{2}\right)^{p_3}~... \right) \\
		& = \left(c^{y_0} \overline{c}^{y_0}\right)^{p_0}~\left(c^{\overline{y_0}} \overline{c}^{\overline{y_0}}\right)^{p_1}~\left(c^{y_0} \overline{c}^{y_0}\right)^{p_2}~\left(c^{\overline{y_0}} \overline{c}^{\overline{y_0}}\right)^{p_3}~... \\
		& = \mathcal{P}_c(y).
	\end{align*}
	Suppose that $\mathpzc{A}$ is odd. If $y = d^\omega$, we have $y = \pi_d(A^\omega) \in \pi_d(\Sigma^\mathbb{N})$ and
	$$\pi_c \sigma_d \pi_d^{-1}(y) = \pi_c \sigma_d (A^\omega) = \pi_c \left( \left(A^\frac{d-1}{2} B C^\frac{d-1}{2}\right)^\omega \right) = \left(c^d \overline{c}^d\right)^\omega = \mathcal{P}_c(x).$$
	If $y \in \mathcal{C}^1$, by setting $(p_n)_{n \geq 0} := \left(\frac{\mathcal{D}(y)_n-1}{2}\right)_{n \geq 0} \in \mathbb{N}^\mathbb{N}$, we get
	$$y = \mathcal{P}_{y_0}(\mathcal{D}(y)) = {y_0}^{2p_0+1}~\overline{y_0}^{2p_1+1}~{y_0}^{2p_2+1}~\overline{y_0}^{2p_3+1}~... = \pi_{y_0}(A^{p_0}~B~C^{p_1}~A^{p_2}~B~C^{p_3}~...) \in \pi_{y_0}(\Sigma^\mathbb{N})$$
	and
	\begin{align*}
		\pi_c \sigma_{y_0} \pi_{y_0}^{-1}(y) & = \pi_c \sigma_{y_0} \left(A^{p_0}~B~C^{p_1}~A^{p_2}~B~C^{p_3}~...\right) \\
		& = \pi_c \left( \left(A^\frac{y_0-1}{2}~B~C^\frac{y_0-1}{2}\right)^{p_0}~A^\frac{y_0-1}{2}~B~C^\frac{\overline{y_0}-1}{2}~\left(A^\frac{\overline{y_0}-1}{2}~B~C^\frac{\overline{y_0}-1}{2}\right)^{p_1}~...\right) \\
		& = \left(c^{y_0} \overline{c}^{y_0}\right)^{p_0} c^{y_0} \overline{c}^{y_1} \left(c^{y_1} \overline{c}^{y_1}\right)^{p_1} \left(c^{y_2} \overline{c}^{y_2}\right)^{p_2} c^{y_2} \overline{c}^{y_3} \left(c^{y_3} \overline{c}^{y_3}\right)^{p_3} ... \\
		& = \mathcal{P}_c(y).
	\end{align*}
	This concludes the proof of \Cref{SP}. Now, for all $n \geq 1$, we get
	\begin{align*}
		\mathcal{P}_{t_0} \mathcal{P}_{t_1} ... \mathcal{P}_{t_{n-1}} (t_n^\omega) & = \pi_{t_0} \sigma_{t_1} \pi_{t_1}^{-1} \pi_{t_1} \sigma_{t_2} \pi_{t_2}^{-1} \pi_{t_2} ... \sigma_{t_{n-1}} \pi_{t_{n-1}}^{-1} \pi_{t_{n-1}} \sigma_{t_n} \pi_{t_n}^{-1} (t_n^\omega) \textup{ with \Cref{SP}} \\
		& = \pi_{t_0} \sigma_{t_1} \sigma_{t_2} ... \sigma_{t_n}(A^\omega) \textup{~~~~~~because } \pi_{t_n}^{-1} (t_n^\omega) = A^\omega.
	\end{align*}
	Then, with \Cref{conjugacy} (i), we deduce that
	$$x = \mathcal{E}(t) = \displaystyle\lim_{n \rightarrow \infty} \mathcal{P}_{t_0} \mathcal{P}_{t_1} ... \mathcal{P}_{t_{n-1}}(t_n^\omega) = \lim_{n \rightarrow \infty} \pi_{t_0} \sigma_{t_1} \sigma_{t_2} ... \sigma_{t_n}(A^\omega).$$
\end{proof}

Then \Cref{conjugacy,smoothSadic} prove \Cref{smoothSadicSystem}. To ease statements and proofs in the rest of the paper, we introduce the following notation.

\begin{definition}
	Over an even or odd alphabet $\mathpzc{A} = \{a,b\}$, we define the map
	\begin{align*}
		\mathcal{F} : \mathpzc{A}^\mathbb{N} & \longrightarrow \Sigma^\mathbb{N} \\
		t & \longmapsto \lim_{n \rightarrow \infty} \sigma_{t_0} \sigma_{t_1} ... \sigma_{t_n}(A^\omega).
	\end{align*}
	In particular, for all $x \in \mathcal{C}^\infty$ we have $x = \pi_{x_0} \circ \mathcal{F} \circ S(t)$ where $t := \mathcal{T}(x)$.
\end{definition}

\begin{remark}\label{F}
	Let $x$ be a smooth sequence over an even or odd alphabet $\mathpzc{A} = \{a,b\}$ and let $t := \mathcal{T}(x)$.
	
	\textup{(i)} We have $\mathcal{D}^k(x) = \pi_{t_k} \circ \mathcal{F} \circ S^{k+1} (t)$ for all $k \geq 0$.
	
	\textup{(ii)} We have $\mathcal{F}(t) = \sigma_{t_0} \sigma_{t_1} ... \sigma_{t_{k-1}} \circ \mathcal{F} \circ S^k (t)$ for all $k \geq 0$.
\end{remark}

\subsection{A morphic representation}

Over even and odd alphabets, \Cref{smoothSadic} provides a morphic representation of the smooth sequences $x$ such that $S(\mathcal{T}(x))$ is periodic.

\begin{lemma}\label{smoothMorphic}
	Let $\mathpzc{A} = \{a,b\}$ be an even or odd alphabet. If $u \in \mathpzc{A}^n$ ($n \geq 1$) and $c \in \mathpzc{A}$, we have
	$$\mathcal{F}(u^\omega) = \sigma^\infty(D) ~~~~~~ \textit{and} ~~~~~~ \mathcal{E}(c u^\omega) = \pi_c \left(\sigma^\infty(D)\right)$$
	where $\sigma := \sigma_{u_1} \sigma_{u_2} ... \sigma_{u_n}$ and $D := \begin{cases} B \textit{ if } a = 1 \textit{ and } u_1 = 1, \\ A \textit{ otherwise} \end{cases}$.
\end{lemma}
\begin{proof}
	The proof is left to the reader.
\end{proof}

\begin{example}
	Over an even or odd alphabet, combining \Cref{fixedPoints} and \Cref{smoothMorphic} yields
	$$\kappa_{a,b} = \begin{cases} \pi_a\left(\sigma_a^\infty(B)\right) \textit{ if } a = 1, \\ \pi_a\left(\sigma_a^\infty(A)\right) \textit{ otherwise } \end{cases} ~~~~~~ \textit{and} ~~~~~~ \kappa_{b,a} = \pi_b\left(\sigma_b^\infty(A)\right).$$
	We recover the construction of Sing in \Cref{k31}.
\end{example}

There are two specific smooth sequences that we need to treat separately.

\begin{lemma}\label{1bb1morphic}
	Let $\mathpzc{A} = \{1,b\}$ be an odd alphabet.
	
	\textup{(i)} We have $\mathcal{F}\left((b1)^\omega\right) = \sigma^\infty(A)$ and $\mathcal{E}\left((1b)^\omega\right) = \pi_1\left(\sigma^\infty(A)\right)$ where $\sigma : A \mapsto A^\frac{b-1}{2} B, B \mapsto A^\frac{b-1}{2} B^\frac{b+1}{2}$.
	
	\textup{(ii)} We have $\mathcal{F}\left((1b)^\omega\right) = \sigma^\infty(B)$ and $\mathcal{E}\left((b1)^\omega\right) = \pi_b\left(\sigma^\infty(B)\right)$ where
$\sigma : B \mapsto B^\frac{b+1}{2} C^\frac{b-1}{2}, C \mapsto B C^\frac{b-1}{2}$.
\end{lemma}
\begin{proof}
	(i) \Cref{smoothMorphic} states that $\mathcal{F}((b1)^\omega) = (\sigma_b \sigma_1)^\infty(A)$ where
	$$\sigma_b \sigma_1 : A \mapsto A^\frac{b-1}{2} B,~~ B \mapsto A^\frac{b-1}{2} B^\frac{b+1}{2},~~ C \mapsto \left[A^\frac{b-1}{2}~B~C^\frac{b-1}{2}\right]^\frac{b-1}{2} A^\frac{b-1}{2} B^\frac{b+1}{2}.$$
	We observe that $(\sigma_b \sigma_1)^n (A) \in \{A,B\}^*$ for all $n \geq 1$ so we deduce that $\mathcal{F}((b1)^\omega) = \sigma^\infty(A)$ where ${\sigma := (\sigma_b \sigma_1)|_{\{A,B\}}}$. Finally we recall that $\mathcal{E}\left((1b)^\omega\right) = \pi_1 \circ \mathcal{F}((b1)^\omega)$.
	
	(ii) Similarly to (i), \Cref{smoothMorphic} states that $\mathcal{F}((1b)^\omega) = (\sigma_1 \sigma_b)^\infty(A)$ where
	$$\sigma_1 \sigma_b : A \mapsto B^\frac{b+1}{2} C^\frac{b-1}{2} \left[A^\frac{b-1}{2} B C^\frac{b-1}{2}\right]^\frac{b-1}{2},~~ B \mapsto B^\frac{b+1}{2} C^\frac{b-1}{2},~~ C \mapsto B C^\frac{b-1}{2}.$$
	We observe that $(\sigma_1 \sigma_b)^n (B) \in \{B,C\}^*$ for all $n \geq 1$ so we deduce that $\mathcal{F}((1b)^\omega) = \sigma^\infty(B)$ where $\sigma := (\sigma_1 \sigma_b)|_{\{B,C\}}$. Finally we recall that $\mathcal{E}\left(b1)^\omega\right) = \pi_b \circ \mathcal{F}((1b)^\omega)$.
\end{proof}

\begin{corollary}
	Over an odd alphabet $\mathpzc{A} = \{1,b\}$, we have $bb \nsqsubset \mathcal{E}\left((1b)^\omega\right)$ and $bb \nsqsubset \mathcal{E}\left((b1)^\omega\right)$.
\end{corollary}

\section{Occurrences in smooth sequences}\label{occurrences}

In this section we provide technical lemmas to be able to tell if a word occurs in a smooth sequence or not. In the process we show that the language of any smooth sequence over an even or odd alphabet is anti-complement-invariant in \Cref{antiInvariance2}, which is the detailed statement of \Cref{antiInvariance}. In particular, the following definition will be crucial.

\begin{definition}\label{shortWords}
	Over a binary alphabet $\mathpzc{A} = \{a,b\}$, we say that a word of the form $c^n$ where $c \in \mathpzc{A}$ and $n \in \llbracket 0,b-1 \rrbracket$ is \textup{short}, and a word that is not short is \textup{long}.
\end{definition}

\subsection{Anti-complement-invariance}

\begin{lemma}\label{abba}
	Let $x$ be a smooth sequence over a binary alphabet $\mathpzc{A} = \{a,b\}$.
	
	\textup{(i)} If $\mathpzc{A}$ is even, then the words $x_0 \overline{x_0}$ and $\overline{x_0} x_0$ occur only at odd positions in $x$.
	
	\textup{(ii)} If $\mathpzc{A}$ is odd, then the word $x_0 \overline{x_0}$ occurs only at even positions and the word $\overline{x_0} x_0$ occurs only at odd positions in $x$.
\end{lemma}
\begin{proof}
	It suffices to consider the canonical factorization $x = {x_0}^{p_0} \overline{x_0}^{p_1} {x_0}^{p_2} \overline{x_0}^{p_3} ...$ and to observe that $p_i + p_{i+1}$ is always even.
\end{proof}

\begin{proposition}\label{antiInvariance2}
	Let $x$ be a smooth sequence over an even or odd alphabet $\mathpzc{A} = \{a,b\}$. If $u \in \mathscr{L}(x)$ and $\mathcal{D}_f(u)$ is long, then $\overline{u} \notin \mathscr{L}(x)$.
\end{proposition}
\begin{proof}
	Let $a_1^{p_1} ... a_n^{p_n}$ be the canonical factorization of $u$ and suppose that $\overline{u} \in \mathscr{L}(x)$. We are going to show that $\mathcal{D}_f(u)$ occurs both at an even and at an odd position in $\mathcal{D}(x)$. First, since $\mathcal{D}_f(u)$ is long, we have $u \neq \varepsilon$ and $n \geq 1$. If $p_1 \leq a$, let $v$ (resp. $v'$) be the element of $\{u, \overline{u}\}$ that starts with $x_0$ (resp. with $\overline{x_0}$). By assumption, we have $v \sqsubset x$ so there exists $i \geq 0$ such that $x = {x_0}^{y_0} \overline{x_0}^{y_1} ... {x_0}^{y_{2i}} \overline{x_0}^{p_2} x_0^{p_3} ...$ so
	$$\mathcal{D}(x) = y_0 y_1 ...y_{2i} p_2 p_3 ... = y_0 y_1 ...y_{2i} \mathcal{D}_f(u) ...,$$
	which means that $\mathcal{D}_f(u)$ occurs at an odd position in $\mathcal{D}(x)$. Also by assumption, we have $v' \sqsubset x$ so there exists $i \geq 0$ such that $x = {x_0}^{y_0} \overline{x_0}^{y_1}...\overline{x_0}^{y_{2i+1}} {x_0}^{p_2} \overline{x_0}^{p_3} ...$ so
	$$\mathcal{D}(x) = y_0 y_1 ... y_{2i+1} p_2 p_3 ... = y_0 y_1 ...y_{2i+1} \mathcal{D}_f(u) ...,$$
	which means that $\mathcal{D}_f(u)$ occurs at an even position in $\mathcal{D}(x)$. If $p_1 > a$, similar arguments yield the same result. Now, if $\mathcal{D}_f(u)$ is of the form $c^n$ where $c \in \mathpzc{A}$ and $n \geq 0$, it is long so $n \geq b$ and it is smooth so $n \leq b$, therefore $n = b$. Then every occurrence of $c^b$ is followed by $\overline{c}$ in $\mathcal{D}(x)$, and the word $c \overline{c}$ is a factor of $\mathcal{D}_f(u) \overline{c}$, so $c \overline{c}$ occurs both at an even and at an odd position in $\mathcal{D}(x)$, which contradicts \Cref{abba}. Similarly, if $\mathcal{D}_f(u)$ contains some $a$'s and some $b$'s, then either $ab$ or $ba$ is a factor of $\mathcal{D}_f(u)$ so it occurs both at an even and at an odd position in $\mathcal{D}(x)$, which contradicts \Cref{abba}.
\end{proof}

\subsection{Bounded gaps}

\begin{definition}
	Let $x$ be a sequence over a finite alphabet $\mathpzc{A}$. We say that a word $u \in \mathpzc{A}^*$ \textup{occurs with bounded gaps} in $x$ if there exists $L > 0$ such that every $v \in \mathscr{L}(x)$ such that $\lvert v \rvert \geq L$ satisfies $u \sqsubset x$, and in that case we write $u \prec x$. Furthermore, if the gaps between two consecutive occurrences at even positions and the gaps between two consecutive occurrence at odd positions are bounded, then we write $u \ll x$.
\end{definition}

\begin{lemma}\label{DDuu}
	Let $x$ be a smooth sequence over an even or odd alphabet $\mathpzc{A} = \{a,b\}$ and let $u \in \mathpzc{C}_f^1$.
	
	\textup{(i)} If $\mathcal{D}_f(u) \sqsubset \mathcal{D}(x)$, then we have $u \sqsubset x$ or $\overline{u} \sqsubset x$.
	
	\textup{(ii)} If $\mathcal{D}_f(u) \prec \mathcal{D}(x)$, then we have $u \prec x$ or $\overline{u} \prec x$.
	
	\textup{(iii)} If $\mathcal{D}_f(u) \ll \mathcal{D}(x)$, then we have $u \prec x$ and $\overline{u} \prec x$.
\end{lemma}
\begin{proof}
	If $u = \varepsilon$ this is trivial. Otherwise let $u = a_1^{p_1} ... a_n^{p_n}$ with $n \geq 1$ and suppose that $p_1 \leq a$ and $p_n \leq a$. Since $x = {x_0}^{y_0} \overline{x_0}^{y_1} {x_0}^{y_2} \overline{x_0}^{y_3} ...$ where $y := \mathcal{D}(x)$, we observe that each occurrence of $\mathcal{D}_f(u) = p_2 ... p_{n-1}$ at an even position in $S(y)$ gives rise to a different occurrence of the word $v := b_2^{p_2} ... b_{n-1}^{p_{n-1}}$ where $b_2 = \overline{x_0}$ in $x$. Furthermore, each occurrence of $v$ in $x$ is included in a different occurrence of the word $w := b_1^{p_1} b_2^{p_2} ... b_{n-1}^{p_{n-1}} b_n^{p_n}$ where $b_1 = x_0$ and $b_n = \overline{b_{n-1}}$ in $x$. Similarly, each occurrence of $\mathcal{D}_f(u)$ at an odd position in $S(y)$ gives rise to a different occurrence of the word $\overline{v}$ in $x$, that is included in a different occurrence of the word $\overline{w}$ in $x$. Finally, we observe that $\{w,\overline{w}\} = \{u,\overline{u}\}$, and that the gap between two occurrences of $w$ (resp. $\overline{w}$) in $x$ is at most $b$ times the gap between the two associated occurrences of $\mathcal{D}_f(u)$ in $S(y)$, which yields (ii) and (iii). To prove (i), \Cref{BJP} (i) states that $y$ is recurrent so, if $\mathcal{D}_f(u) \sqsubset y$, then $\mathcal{D}_f(u) \sqsubset S(y)$ and the previous reasoning yields $u \sqsubset x$ or $\overline{u} \sqsubset x$.
	
	The proof is similar in the three other cases ($p_1 \leq a$ and $p_n > a$; $p_1 > a$ and $p_n \leq a$; $p_1 > a$ and $p_n > a$).
\end{proof}
\newpage
\begin{lemma}\label{ACbgF}
	Let $\mathpzc{A} = \{1,b\}$ be an odd alphabet. If $t \in \mathpzc{A}^\mathbb{N}$ is such that $t \notin \left\{(1b)^\omega,(b1)^\omega\right\}$, then we have $A,C \prec \mathcal{F}(t)$.
\end{lemma}
\begin{proof}
	The fact that $t \notin \left\{(1b)^\omega,(b1)^\omega\right\}$ provides $k \geq 0$ such that $t_k = t_{k+1}$, and \Cref{F} (ii) yields $\mathcal{F} \circ S^k(t) = {\sigma_{t_k}}^2 \sigma_{t_{k+2}} \circ \mathcal{F} \circ S^{k+3}(t)$. Now we observe that $B$ occurs in the image of every letter of $\Sigma$ by $\sigma_1$ and by $\sigma_b$ so we deduce that $B \prec \sigma_{k+2} \circ \mathcal{F} \circ S^{k+3}(t)$. Then, we observe that $A,C \sqsubset {\sigma_1}^2(B), {\sigma_b}^2(B)$ so we deduce that $A,C \prec {\sigma_{t_k}}^2 \sigma_{t_{k+2}} \circ \mathcal{F} \circ S^{k+3}(t) = \mathcal{F} \circ S^k(t)$. Finally, we observe that $A,C \sqsubset \sigma_1(C), \sigma_b(A)$ so we deduce that $A,C \prec \sigma_{t_0} \sigma_{t_1} ... \sigma_{t_{k-1}} \circ \mathcal{F} \circ S^k (t) = \mathcal{F}(t)$ thanks to \Cref{F} (ii).
\end{proof}

\begin{lemma}\label{limitABC}
	Let $\mathpzc{A} = \{1,b\}$ be an odd alphabet and let $t \in \mathpzc{A}^\mathbb{N}$. Then every factor of $\mathcal{F}(t)$ of length $\frac{b+1}{2}$ either contains $B$ or contains $A$ and $C$.
\end{lemma}
\begin{proof}
	An equivalent statement is that the words $A^\frac{b+1}{2}$ and $C^\frac{b+1}{2}$ do not occur in $\mathcal{F}(t)$. We observe that $\pi_1(A^\frac{b+1}{2}) = 1^{b+1} \notin \mathcal{C}_f^\infty$ (resp. $\pi_1(C^\frac{b+1}{2}) = b^{b+1} \notin \mathcal{C}_f^\infty$), and we deduce that $\pi_1(A^\frac{b+1}{2}) \nsqsubset \mathcal{E}(1t) = \pi_1 \circ \mathcal{F}(t)$ (resp. $\pi_1(C^\frac{b+1}{2}) \nsqsubset \pi_1 \circ \mathcal{F}(t)$) and that $A^\frac{b+1}{2} \nsqsubset \mathcal{F}(t)$ (resp. $C^\frac{b+1}{2} \nsqsubset \mathcal{F}(t)$).
\end{proof}

\begin{definition}
	Let $x$ be a smooth sequence over an odd alphabet $\mathpzc{A} = \{1,b\}$. We say that $x$ has property $\mathfrak{P}3$ if $S^2(\mathcal{T}(x)) \notin \left\{(1b)^\omega,(b1)^\omega\right\}$.
\end{definition}

\begin{lemma}\label{Sfactors}
	Let $x$ be a smooth sequence over an even or odd alphabet $\mathpzc{A} = \{a,b\}$ that has property $\mathfrak{P}3$ if $a = 1$. If $u \in \mathpzc{C}_f^1$ and $\mathcal{D}_f(u)$ is short, then $u \prec x$.
\end{lemma}
\begin{proof}
	Let $\mathcal{D}_f(u) = c^n$ with $c \in \mathpzc{A}$ and $n \in \llbracket 0,b-1 \rrbracket$. With \Cref{DDuu} (iii) it suffices to show that $c^n \ll \mathcal{D}(x)$, and since $c^n$ occurs with both parities in $c^b$ it suffices to show that $c^b \prec \mathcal{D}(x)$.
	
	If $a > 1$, every occurence of $b$ in $\mathcal{D}^2(x)$ is preceded or followed by another $b$ so we have $b \ll \mathcal{D}^2(x)$ and \Cref{DDuu} (iii) yields $c^b \prec \mathcal{D}(x)$.

	  If $a = 1$, we have $S^2(t) \notin \left\{(1b)^\omega,(b1)^\omega\right\}$ so \Cref{ACbgF} states that $A,C \prec \mathcal{F} \circ S^2(t)$ where $t := \mathcal{T}(x)$. Since ${t_1}^2 = \pi_{t_1}(A)$ and $\overline{t_1}^2 = \pi_{t_1}(C)$, we deduce that ${t_1}^2,\overline{t_1}^2 \prec \pi_{t_1} \circ \mathcal{F} \circ S^2(t) = \mathcal{D}(x)$ thanks to \Cref{F}. Since $a = 1$, every occurrence of ${t_1}^2$ (resp. of $\overline{t_1}^2$) in $\mathcal{D}(x)$ is included in an occurrence of ${t_1}^b$ (resp. of $\overline{t_1}^b$) in $\mathcal{D}(x)$ so we deduce that $c^b \prec \mathcal{D}(x)$.
\end{proof}

\begin{remark}\label{pathological}
	Over an odd alphabet $\mathpzc{A} = \{1,b\}$, the eight smooth sequences do not have property $\mathfrak{P}3$ are $\mathcal{E}(11(1b)^\omega)$, $\mathcal{E}((1b)^\omega)$, $\mathcal{E}(b1(1b)^\omega)$, $\mathcal{E}(bb(1b)^\omega)$, $\mathcal{E}(11(b1)^\omega)$, $\mathcal{E}(1b(b1)^\omega)$, $\mathcal{E}((b1)^\omega)$ and $\mathcal{E}(bb(b1)^\omega)$. One can check that \Cref{Sfactors} does not hold for them, for example in \Cref{1bb1morphic} we show that $\mathcal{E}((1b)^\omega)$ and $\mathcal{E}((b1)^\omega)$ do not contain the word $bb$.
\end{remark}

\subsection{Forbidden and mandatory factors}

\subsubsection{Over even alphabets}

\begin{lemma}\label{evenForbidden}
	If $x$ is a smooth sequence over an even alphabet $\mathpzc{A} = \{a,b\}$, define the words
	\begin{align*}
		& f_1 := {x_0}^{a+1} \overline{x_0}^a x_0, & & f_2 := \overline{x_0} {x_0}^a \overline{x_0}^{a+1}, \\
		& f_3 := {x_0} \left[\overline{x_0}^a {x_0}^a\right]^\frac{b}{2} \overline{x_0}, & & f_4 := \overline{x_0}^{a+1} \left[{x_0}^b \overline{x_0}^b\right]^\frac{b-2}{2} {x_0}^{a+1}.
	\end{align*}
	If $u \in \mathpzc{A}^*$ satisfies $f_i \sqsubset u$ for some $i \in \{1,2,3,4\}$, then we have $u \not\sqsubset x$.
\end{lemma}
\begin{proof}
	It suffices to show that $f_i \not\sqsubset x$ for all $i \in \{1,2,3,4\}$.
	
	If $f_1 \sqsubset x$, this provides $i \geq 0$ such that $x = {x_0}^{p_0} \overline{x_0}^{p_1} ... \overline{x_0}^{p_{2i-1}} {x_0}^b \overline{x_0}^a {x_0}^{p_{2i+2}} ...$ and $\mathcal{D}(x) = p_0 p_1 ... p_{2i-1} b a ...$ which contradicts \Cref{abba} (i).
	
	If $f_2 \sqsubset x$, this provides $i \geq 0$ such that $x = {x_0}^{p_0} \overline{x_0}^{p_1} ... \overline{x_0}^{p_{2i-1}} {x_0}^a \overline{x_0}^b {x_0}^{p_{2i+2}} ...$ and $\mathcal{D}(x) = p_0 p_1 ... p_{2i-1} a b ...$ which contradicts \Cref{abba} (i).
	
	If $f_3 \sqsubset x$, this provides $i \geq 0$ such that $x = {x_0}^{p_0} \overline{x_0}^{p_1} ... {x_0}^{p_{2i}} \left[\overline{x_0}^a {x_0}^a\right]^\frac{b}{2} \overline{x_0}^{p_{2i+b+1}} ...$ Then, if $p_{2i} = b$ this contradicts \Cref{abba} (i), and if $p_{2i} = a$ we have $\mathcal{D}(x) = p_0 p_1 ... p_{2i-1} a^{b+1} p_{2i+b+1} ... \notin \mathpzc{C}^1$.
	
	If $f_4 \sqsubset x$, this provides $i \geq 0$ such that $x = {x_0}^{p_0} \overline{x_0}^{p_1} ... {x_0}^{p_{2i}} \left[\overline{x_0}^b {x_0}^b\right]^\frac{b}{2} \overline{x_0}^{p_{2i+b+1}} ...$ Then, if $p_{2i} = a$ this contradicts \Cref{abba} (i), and if $p_{2i} = b$ we have $\mathcal{D}(x) = p_0 p_1 ... p_{2i-1} b^{b+1} p_{2i+b+1} ... \notin \mathpzc{C}^1$.
\end{proof}

From the forbidden factors given in \Cref{evenForbidden}, we obtain a convenient lemma to find some mandatory factors.

\begin{lemma}\label{tinyEvenMandatory}
	Let $x$ be a smooth sequence over an even alphabet $\mathpzc{A}$. If $u \in \mathcal{C}_f^2$ is such that $\left(\mathcal{D}_f\right)^2(u)$ is short and that $f_i \sqsubset \overline{u}$ for some $i \in \{1,2,3,4\}$, then we have $u \prec x$. 
\end{lemma}
\begin{proof}
	\Cref{Sfactors} states that $\mathcal{D}_f(u) \prec \mathcal{D}(x)$ so \Cref{DDuu} (ii) yields $u \prec x$ or $\overline{u} \prec x$, but \Cref{evenForbidden} states $\overline{u} \nsqsubset x$.
\end{proof}

\begin{lemma}\label{evenMandatory}
	If $x$ is a smooth sequence over an even alphabet $\mathpzc{A} = \{a,b\}$, define the words
	\begin{align*}
		& m_1 := \overline{x_0}^{a+1} \left[{x_0}^a \overline{x_0}^a\right]^\frac{b}{2} {x_0}^{a+1}, & & m_2 := \overline{x_0}^{a+1} \left[{x_0}^a \overline{x_0}^a\right]^\frac{a}{2} {x_0}^{a+1}, \\
		& m_3 := x_0 \overline{x_0}^a \left[{x_0}^b \overline{x_0}^b\right]^\frac{b}{2} {x_0}^a \overline{x_0}, & & m_4 := x_0 \overline{x_0}^a \left[{x_0}^b \overline{x_0}^b\right]^\frac{a}{2} {x_0}^a \overline{x_0}.
	\end{align*}
	If $u \in \mathpzc{A}^*$ and there exists $i \in \{1,2,3,4\}$ such that $u \sqsubset m_i$, then we have $u \sqsubset x$.
\end{lemma}
\begin{proof}
	It suffices to show that $m_i \sqsubset x$ for all $i \in \{1,2,3,4\}$. We observe that $\left(\mathcal{D}_f\right)^2(m_1) = \mathcal{D}_f(b a^b b) = b$ is short and that $f_1 \sqsubset_p \overline{m_1}$; $\left(\mathcal{D}_f\right)^2(m_2) = \mathcal{D}_f(b a^a b) = a$ and $f_1 \sqsubset_p \overline{m_2}$; $\left(\mathcal{D}_f\right)^2(m_3) = \mathcal{D}_f(a b^b a) = b$ and $f_2 \sqsubset_p \overline{m_3}$; $\left(\mathcal{D}_f\right)^2(m_4) = \mathcal{D}_f(a b^a a) = a$ and $f_2 \sqsubset_p \overline{m_4}$; so \Cref{tinyEvenMandatory} yields $m_i \sqsubset x$ for all $i \in \{1,2,3,4\}$.
\end{proof}

\subsubsection{Over odd alphabets}

\begin{lemma}\label{oddForbidden}
	If $x$ is a smooth sequence over an odd alphabet $\mathpzc{A} = \{a,b\}$, define the forbidden words
	\begin{align*}
		& f_1 := x_0 \overline{x_0}^a {x_0}^{a+1}, & & f_2 := {x_0}^{a+1} \overline{x_0}^a x_0, \\
		& f_3 := x_0 \left[\overline{x_0}^a {x_0}^a\right]^\frac{b-1}{2} \overline{x_0}^a x_0, & & f_4 := {x_0}^{a+1} \left[\overline{x_0}^b {x_0}^b\right]^\frac{b-3}{2} \overline{x_0}^b {x_0}^{a+1}.
	\end{align*}
	If $\mathcal{D}(x)_0 = a$ (resp. if $\mathcal{D}(x)_0 = b$), $u \in \mathpzc{A}^*$ and $f_i \sqsubset u$ (resp. $\overline{f_i} \sqsubset u$) for some $i \in \{1,2,3,4\}$, then we have $u \not\sqsubset x$.
\end{lemma}
\begin{proof}
	Suppose that $\mathcal{D}(x)_0 = a$, so that \Cref{abba} (ii) states that $ab$ occurs only at even positions in $\mathcal{D}(x)$ and $ba$ occurs only at odd positions in $\mathcal{D}(x)$. Then it suffices to show that $f_i \not\sqsubset x$ for all $i \in \{1,2,3,4\}$.
	
	If $f_1 \sqsubset x$, this provides $i \geq 0$ such that $x = {x_0}^{p_0} \overline{x_0}^{p_1} ... {x_0}^{p_{2i}} \overline{x_0}^a {x_0}^b ...$ and $\mathcal{D}(x) = p_0 p_1 ... p_{2i} a b ...$ which contradicts \Cref{abba} (ii).
	
	If $f_2 \sqsubset x$, this provides $i \geq 0$ such that $x = {x_0}^{p_0} \overline{x_0}^{p_1} ... \overline{x_0}^{p_{2i-1}} {x_0}^b \overline{x_0}^a {x_0}^{p_{2i+2}} ...$ and $\mathcal{D}(x) = p_0 p_1 ... p_{2i-1} b a ...$ which contradicts \Cref{abba} (ii).
	
	If $f_3 \sqsubset x$, this provides $i \geq 0$ such that $x = {x_0}^{p_0} \overline{x_0}^{p_1} ... {x_0}^{p_{2i}} \left[\overline{x_0}^a {x_0}^a\right]^\frac{b-1}{2} \overline{x_0}^a {x_0}^{p_{2i+b+1}} ...$ and $\mathcal{D}(x) = p_0 ... p_{2i}~a^b~p_{2i+b+1} ...$ Then, if $p_{2i} = b$ this contradicts \Cref{abba} (ii), and if $p_{2i} = a$ we have $\mathcal{D}(x) \notin \mathpzc{C}^1$.
	
	If $f_4 \sqsubset x$, this provides $i \geq 0$ such that $x = {x_0}^{p_0} \overline{x_0}^{p_1} ... \overline{x_0}^{p_{2i-1}} \left[{x_0}^b \overline{x_0}^b\right]^\frac{b-1}{2} {x_0}^b \overline{x_0}^{p_{2i+b}}...$ and $\mathcal{D}(x) = p_0 ... p_{2i-1} b^b p_{2i+b} ...$ Then, if $p_{2i-1} = a$ this contradicts \Cref{abba} (ii), and if $p_{2i} = b$ we have $\mathcal{D}(x) \notin \mathpzc{C}^1$.
	
	The proof is similar when $\mathcal{D}(x)_0 = b$.
\end{proof}

From the forbidden factors given in \Cref{oddForbidden}, we obtain a convenient lemma to find some mandatory factors.

\begin{lemma}\label{tinyOddMandatory}
	Let $x$ be a smooth sequence over an odd alphabet $\mathpzc{A} = \{a,b\}$ that has property $\mathfrak{P}3$ if $a = 1$. If $\mathcal{D}(x)_0 = a$ (resp. $\mathcal{D}(x)_0 = b$), $u \in \mathcal{C}_f^2$ is such that $\mathcal{D}^2_f(u)$ is short and there exists $i \in \{1,2,3,4\}$ such that $f_i \sqsubset \overline{u}$ (resp. $\overline{f_i} \sqsubset \overline{u}$), then we have $u \sqsubset x$.
\end{lemma}
\begin{proof}
	\Cref{Sfactors} states that $\mathcal{D}_f(u) \prec \mathcal{D}(x)$ so \Cref{DDuu} (ii) yields $u \prec x$ or $\overline{u} \prec x$, but \Cref{oddForbidden} states $\overline{u} \nsqsubset x$.
\end{proof}

\begin{lemma}\label{oddMandatory}
	If $x$ is a smooth sequence over an odd alphabet $\mathpzc{A} = \{a,b\}$ that has property $\mathfrak{P}3$ if $a = 1$, define the mandatory words
	\begin{align*}
		& m_1 := \overline{x_0}^{a+1} \left[{x_0}^a \overline{x_0}^a\right]^\frac{b-1}{2} {x_0}^a \overline{x_0}^{a+1}, & & m_2 := \overline{x_0}^{a+1} \left[{x_0}^a \overline{x_0}^a\right]^\frac{a-1}{2} {x_0}^a \overline{x_0}^{a+1}, \\
		& m_3 := \overline{x_0} {x_0}^a \left[\overline{x_0}^b {x_0}^b\right]^\frac{b-1}{2} \overline{x_0}^b {x_0}^a \overline{x_0}, & & m_4 := \overline{x_0} {x_0}^a \left[\overline{x_0}^b {x_0}^b\right]^\frac{a-1}{2} \overline{x_0}^b {x_0}^a \overline{x_0}.
	\end{align*}
	If $\mathcal{D}(x)_0 = a$ (resp. if $\mathcal{D}(x)_0 = b$), $u \in \mathpzc{A}^*$ and $u \sqsubset m_i$ (resp. $u \sqsubset \overline{m_i}$) for some $i \in \{1,2,3,4\}$, then we have $u \sqsubset x$.
\end{lemma}
\begin{proof}
	Suppose that $\mathcal{D}(x)_0 = a$. Then it suffices to show that $m_i \sqsubset x$ for all $i \in \{1,2,3,4\}$. We observe that $\left(\mathcal{D}_f\right)^2(m_1) = \mathcal{D}_f(b a^b b) = b$ is short and that $f_1 \sqsubset_p \overline{m_1}$; $\left(\mathcal{D}_f\right)^2(m_2) = \mathcal{D}_f(b a^a b) = a$ and $f_1 \sqsubset_p \overline{m_2}$; $\left(\mathcal{D}_f\right)^2(m_3) = \mathcal{D}_f(a b^b a) = b$ and $f_2 \sqsubset_p \overline{m_3}$; $\left(\mathcal{D}_f\right)^2(m_4) = \mathcal{D}_f(a b^a a) = a$ and $f_2 \sqsubset_p \overline{m_4}$; so \Cref{tinyOddMandatory} yields $m_i \sqsubset x$ for all $i \in \{1,2,3,4\}$.
	
	The proof is similar when $\mathcal{D}(x)_0 = b$.
\end{proof}

\section{Recurrence of smooth sequences}\label{recurrence}

In this subsection we show \Cref{URLR}. Since a smooth sequence over an even or odd alphabet $\mathpzc{A}$ has the form $\pi_c \circ \mathcal{F} \circ S(t)$ where $c \in \mathpzc{A}$ and $t \in \mathpzc{A}^\mathbb{N}$, we investigate the uniform (resp. linear) recurrence of the sequence $\mathcal{F} \circ S(t)$.

\subsection{Uniform recurrence}

To prove uniform recurrence, we first handle a special case.

\begin{lemma}\label{1bLR}
	Let $\mathpzc{A} = \{1,b\}$ be an odd alphabet and let $t \in \mathpzc{A}^\mathbb{N}$. If there exists $n \geq 0$ such that $S^n(t) = (1b)^\omega$, then the sequence $\mathcal{F}(t)$ is linearly recurrent.
\end{lemma}
\begin{proof}
	We have $\mathcal{F}(t) = \sigma_{t_0} \sigma_{t_1} ... \sigma_{t_{n-1}} \circ \mathcal{F} \circ S^n(t)$, and \Cref{1bb1morphic} (ii) states that $\mathcal{F} \circ S^n(t) = \mathcal{F}\left((1b)^\omega\right) = \sigma^\infty(B)$ where $\sigma$ is a primitive substitution. It is known that a purely morphic sequence generated by a primitive substitution is linearly recurrent (see \cite{Queffelec}), so we deduce that $\mathcal{F} \circ S^n(t)$ and $\mathcal{F}(t)$ are linearly recurrent.
\end{proof}

Now we can show uniform recurrence everywhere, this is \Cref{URLR} (i).

\begin{proof}[Proof of \Cref*{URLR} \textup{(i)}]
	Let $\mathpzc{A} = \{a,b\}$ be an even or odd alphabet. We are going to show that, for all $x \in \mathpzc{C}^\infty$, every word $u \in \mathscr{L}(x)$ satisfies $u \prec x$. We proceed by induction on the height of $u$: first, if $u$ has height $0$, i.e., $u = \varepsilon$, this is trivial. Now let $h \geq 0$ be such that, for all $x \in \mathpzc{C}^\infty$, every word of $u \in \mathscr{L}(x)$ of height $h$ satisfies $u \prec x$. Let $x \in \mathpzc{C}^\infty$, let $t := \mathcal{T}(x)$ and let $u \in \mathscr{L}(x)$ be of height $h+1$. 
	
	If $\mathcal{D}_f(u)$ is short and $S^2(t) \notin \left\{(1b)^\omega,(b1)^\omega\right\}$ if $a = 1$, \Cref{Sfactors} directly states that $u \prec x$.
	
	If $\mathcal{D}_f(u)$ is short, $a = 1$ and $S^2(t) \in \left\{(1b)^\omega,(b1)^\omega\right\}$, we have $S^2(t) = (1b)^\omega$ or $S^3(t) = (1b)^\omega$ so \Cref{1bLR} states that $\mathcal{F} \circ S (t)$ is linearly recurrent. We deduce that $x = \pi_{x_0} \circ \mathcal{F} \circ S (t)$ is also linearly recurrent, and a fortiori we have $u \prec x$.
	
	If $\mathcal{D}_f(u)$ is long, consider the canonical factorization $u = a_1^{p_1} ... a_n^{p_n}$ with $n \geq 1$. Then a quick inspection of four cases (wether $p_1$ and $p_n$ are $\leq a$ or $>a$) shows that $\mathcal{D}_f(u) \in \mathscr{L}(\mathcal{D}(x))$. Since the word $\mathcal{D}_f(u)$ has height $h$, by hypothesis, we have $\mathcal{D}_f(u) \prec \mathcal{D}(x)$. Now \Cref{DDuu} (ii) ensures that $u \prec x$ or $\overline{u} \prec x$, but \Cref{antiInvariance2} ensures that $\overline{u} \not\sqsubset x$.
\end{proof}

\subsection{Linear recurrence}

In \Cref{1bLR} we proved linear recurrence in a special case, now we prove it for the class of smooth sequences $\mathscr{M}$ (cf. \Cref{propertyM}).

\begin{lemma}\label{FLR}
	Let $\mathpzc{A} = \{a,b\}$ be an even or odd alphabet and let $t \in \mathpzc{A}^\mathbb{N}$. If $a > 1$, or if $a = 1$ and there exists $n \geq 1$ such that $(1b)^n \nsqsubset t$, then $\mathcal{F}(t)$ is a proper primitive $S$-adic sequence and it is linearly recurrent.
\end{lemma}
\begin{proof}
	Recall that $\mathcal{F}(t) := \displaystyle \lim_{n \rightarrow \infty} \sigma_{t_0} \sigma_{t_1} ... \sigma_{t_n}(A^\omega)$ is an $S$-adic sequence.
	
	If $a > 1$, we immediatly observe that the substitutions $\sigma_a$ and $\sigma_b$ are proper and have positive entries so $\mathcal{F}(t)$ is a proper primitive $S$-adic sequence.
	
	If $a = 1$ and there exists $n \geq 0$ such that $(1b)^n \nsqsubset t$, we have the substitutions
	\begin{align*}
		\sigma_1 :~~~ & A \longmapsto B, & & B \longmapsto B~C^\frac{b-1}{2}, & & C \longmapsto A^\frac{b-1}{2}~B~C^\frac{b-1}{2}, \\
		\sigma_b :~~~ & A \longmapsto A^\frac{b-1}{2}~B~C^\frac{b-1}{2}, & & B \longmapsto A^\frac{b-1}{2}~B, & & C \longmapsto B.
	\end{align*}
	We observe that the substitutions $\sigma_1^2$, $\sigma_1 \sigma_b$, $\sigma_b\sigma_1$ and $\sigma_b^2$ are proper (every image starts with $B$ and ends with $C$). Also, the fact that $(1b)^n \nsqsubset t$ implies that one of the words $111$, $11b$, $bbb$ or $bb1$ occurs with bounded gaps in $t$, and one can check that the images of every letter of $\Sigma$ by the substitutions $\sigma_1^3$, $\sigma_1^2 \sigma_b$, $\sigma_b^3$ and $\sigma_b^2 \sigma_1$ all contain every letter of $\Sigma$. We deduce that $\mathcal{F}(t)$ is a proper and primitive $S$-adic sequence.
	
	To conclude the proof, \cite[Proposition 1.1]{Durand03} states that a proper primitive $S$-adic sequence is linearly recurrent.
\end{proof}

\begin{proposition}\label{LR}
	A smooth sequence in $\mathscr{M}$ is linearly recurrent.
\end{proposition}
\begin{proof}
	Let $x$ be a smooth sequence in $\mathscr{M}$ over an even or odd alphabet $\mathpzc{A} = \{a,b\}$ and let $t := \mathcal{T}(x)$. Since $x = \pi_{x_0} \circ \mathcal{F} \circ S(t)$, it suffices to show that $\mathcal{F} \circ S(t)$ is linearly recurrent.
	
	If $a > 1$, or if $a = 1$ and $x$ has property $\mathfrak{P}2$ (i.e., there exists $n \geq 1$ such that $(1b)^n \nsqsubset t$), then we have $(1b)^n \nsqsubset S(t)$ so \Cref{FLR} states that $\mathcal{F} \circ S(t)$.
	
	If $a = 1$ and $x$ has property $\mathfrak{P}1$ (i.e., there exists $n \geq 0$ such that $S^n(t) = (1b)^\omega$), then we have $S^{n+1} (S(t)) = (1b)^\omega$ so \Cref{1bLR} states that $\mathcal{F} \circ S(t)$ is linearly recurrent.
\end{proof}

\subsection{Non-linear recurrence}

To disprove linear recurrence in the remaining cases, we first show a recognizability property in our $S$-adic system.

\begin{lemma}\label{recognizability}
	Let $\mathpzc{A} = \{1,b\}$ be an odd alphabet and let $t \in \mathpzc{A}^\mathbb{N}$.
	
	\textup{(i)} If there exists $c \in \mathpzc{A}$, $u,v \sqsubset \Sigma^*$ and $D,E \in \Sigma$ such that $\sigma_c(DuE) \sqsubset \sigma_c(v)$, then we have $u \sqsubset v$.
	
	\textup{(ii)} If there exists $k \geq 0$, $u,v \sqsubset \Sigma^*$ and $w,w' \in \Sigma^\frac{b+3}{2}$ such that $wuw', v \sqsubset \mathcal{F} \circ S^{k+1}(t)$ and $\sigma(wuw') \sqsubset \sigma(v)$ where $\sigma := \sigma_{t_1} \sigma_{t_2} ... \sigma_{t_k}$ ($\sigma := Id$ if $k = 0$), then we have $u \sqsubset v$.
\end{lemma}
\begin{proof}
	(i) Recall the substitutions
	\begin{align*}
		\sigma_1 :~~~ & A \longmapsto B, & & B \longmapsto B~C^\frac{b-1}{2}, & & C \longmapsto A^\frac{b-1}{2}~B~C^\frac{b-1}{2}, \\
		\sigma_b :~~~ & A \longmapsto A^\frac{b-1}{2}~B~C^\frac{b-1}{2}, & & B \longmapsto A^\frac{b-1}{2}~B, & & C \longmapsto B.
	\end{align*}
	We observe that the image of every letter of $\Sigma$ by $\sigma_1$ (by $\sigma_b$) is of the form $A^p B C^q$ where $(p,q)$ is unique to each letter. In particular, one can uniquely desubstitute the word $\sigma_1(u)$ (resp. $\sigma_b(u)$) up to the first and last letters of $u$ because $\sigma_1(A) \sqsubset_p \sigma_1(B) \sqsubset_s \sigma_1(C)$ (resp. because $\sigma_b(C) \sqsubset_s \sigma_b(B) \sqsubset_p \sigma_b(A)$).  Therefore, by adding the letters $D,E \in \Sigma$, we deduce that $\sigma_c(DuE) \sqsubset \sigma_c(v)$ implies $u \sqsubset v$.

	(ii) We proceed by induction on $k$. If $k = 0$ this is trivial, now let $k \geq 0$ be such that (ii) holds. Let $u,v \sqsubset \Sigma^*$ and $w,w' \in \Sigma^\frac{b+3}{2}$ such that $wuw', v \sqsubset \mathcal{F} \circ S^{k+2}(t)$ and $\sigma(wuw') \sqsubset \sigma(v)$ where $\sigma := \sigma_{t_1} \sigma_{t_2} ... \sigma_{t_{k+1}}$. Let $D := w_\frac{b+3}{2}$ and let $E := w'_1$, so that
	\begin{equation}\label{nextStep}
		\sigma_{t_{k+1}} \left( w_{\left[1,\frac{b+1}{2}\right]}~D~u~E~w'_{\left[2,\frac{b+3}{2}\right]} \right) = \sigma_{t_{k+1}}(wuw') \sqsubset \sigma_{t_{k+1}} \circ \mathcal{F} \circ S^{k+2}(t) = \mathcal{F} \circ S^{k+1}(t)
	\end{equation}
	and
	$$\sigma_{t_{k+1}}(v) \sqsubset \sigma_{t_{k+1}} \circ \mathcal{F} \circ S^{k+2}(t) = \mathcal{F} \circ S^{k+1}(t).$$
	
	Now \Cref{limitABC} states that the words $w_{\left[1,\frac{b+1}{2}\right]}$ and $w'_{\left[2,\frac{b+3}{2}\right]}$ either contain $B$ or contain $A$ and $C$, and we observe that $\lvert \sigma_{t_{k+1}}(B) \rvert \geq 2$ and $\lvert \sigma_{t_{k+1}}(C) \rvert \geq 2$ if $t_{k+1} = 1$ (resp. $\lvert \sigma_{t_{k+1}}(A) \rvert \geq 2$ if $t_{k+1} = b$). We deduce that $\ell := \left\lvert \sigma_{t_{k+1}}\left(w_{\left[1,\frac{b+1}{2}\right]}\right) \right\rvert \geq \frac{b+3}{2}$ and $\ell' := \left\lvert \sigma_{t_{k+1}}\left(w'_{\left[2,\frac{b+3}{2}\right]}\right) \right\rvert \geq \frac{b+3}{2}$. By setting $\sigma' := \sigma_{t_1} \sigma_{t_2} ... \sigma_{t_k}$, we get
	\begin{align*}
		\sigma'\left(\sigma_{t_{k+1}}\left(w_{\left[1,\frac{b+1}{2}\right]}\right)_{\left[\ell-\frac{b+1}{3},\ell\right]} ~ \sigma_{t_{k+1}}(DuE) ~ \sigma_{t_{k+1}}\left(w'_{\left[2,\frac{b+3}{2}\right]}\right)_{\left[1,\frac{b+3}{2}\right]} \right) &~ \sqsubset \sigma' \sigma_{t_{k+1}} \left( w_{\left[1,\frac{b+1}{2}\right]}~D~u~E~w'_{\left[2,\frac{b+3}{2}\right]} \right) \\
		&~\sqsubset \sigma' \sigma_{t_{k+1}} (wuw') = \sigma(wuw') ~~~\textup{ thanks to \Cref{nextStep}} \\
		&~\sqsubset \sigma(v) = \sigma' \sigma_{t_{k+1}} (v) ~~~~~~~~\textup{ by hypothesis}.
	\end{align*}
	Then the induction hypothesis yields $\sigma_{t_{k+1}}(DuE) \sqsubset \sigma_{t_{k+1}}(v)$ and (i) yields $u \sqsubset v$.
\end{proof}

\begin{lemma}\label{FnotLR}
	Let $\mathpzc{A} = \{1,b\}$ be an odd alphabet and let $t \in \mathpzc{A}^\mathbb{N}$ such that 
$(1b)^n \sqsubset t$ for all $n \geq 1$ and $S^n(t) \neq (1b)^\omega$ for all $n \geq 0$. Then the sequence $\mathcal{F}(t)$ is not linearly recurrent.
\end{lemma}
\begin{proof}
	The condition on $t$ means that one of the following properties holds:
	
	(i) $1(1b)^n \sqsubset t$ for all $n \geq 0$,
	
	(ii) $b(b1)^n \sqsubset t$ for all $n \geq 0$.
	
	Suppose that (i) holds, let $n \geq 0$ and let $u := (\sigma_1 \sigma_b)^n (B)$. Then there exists $k \geq 0$ such that $t_{[k,k+2n]} = 1(1b)^n$ and it is clear that $B \sqsubset \mathcal{F} \circ S^{k+2n+1} (t)$ so we deduce that $u \sqsubset (\sigma_1 \sigma_b)^n \circ \mathcal{F} \circ S^{k+2n+1} (t) = \mathcal{F} \circ S^{k+1} (t)$. Furthermore, we have $\sigma_1 \sigma_b (B) = B^\frac{b+1}{2} C^\frac{b-1}{2}$ and $\sigma_1 \sigma_b (C) = B C^\frac{b-1}{2}$ so we deduce that $u \in \{B,C\}^*$ and that $\lvert u \rvert_C \geq n$. Also, \Cref{ACbgF} states that $A \prec \mathcal{F} \circ S^{k+2n+1}(t)$, which provides $w,w' \in \Sigma^\frac{b+3}{2}$ such that $wAw' \sqsubset \mathcal{F} \circ S^{k+1}(t)$. Now set $\sigma := \sigma_{t_0} ... \sigma_{t_{k-1}}$, so that $\sigma \sigma_1 (u), \sigma \sigma_1(wAw') \sqsubset \sigma \sigma_1 \circ \mathcal{F} \circ S^{k+1} (t) = \mathcal{F}(t)$ since $t_k = 1$. We also have
	\begin{equation}\label{notFactor}
		\sigma \sigma_1(wAw') \nsqsubset \sigma \sigma_1(u),
	\end{equation}
	otherwise \Cref{recognizability} (ii) would yield $A \sqsubset u$. Furthermore, we define $V(v) := \left(\begin{matrix} \left\lvert v \right\rvert_A \\ \left\lvert v \right\rvert_B \\ \left\lvert v \right\rvert_C \end{matrix}\right)$ for all $v \in \Sigma^*$, so that $V(\sigma(v)) = M_\sigma V(v)$. We get
	$$V(\sigma_1(wAw')) = V(\sigma_1(w)) +  V(\sigma_1(A)) + V(\sigma_1(w')) \leq \frac{b+3}{2} \left(\begin{matrix} \frac{b-1}{2} \\ 1 \\ \frac{b-1}{2} \end{matrix}\right) + \left(\begin{matrix} 0 \\ 1 \\ 0 \end{matrix}\right) + \frac{b+3}{2} \left(\begin{matrix} \frac{b-1}{2} \\ 1 \\ \frac{b-1}{2} \end{matrix}\right) \leq (b+4) \left(\begin{matrix} \frac{b-1}{2} \\ 1 \\ \frac{b-1}{2} \end{matrix}\right)$$
	and
	$$V(\sigma_1(u)) \geq \lvert u \rvert_C V(\sigma_1(C)) \geq n \left(\begin{matrix} \frac{b-1}{2} \\ 1 \\ \frac{b-1}{2} \end{matrix}\right) \geq \frac{n}{b+4} V(\sigma_1(wAw')).$$
	We deduce that
	$$V(\sigma \sigma_1(u)) = M_\sigma V(\sigma_1(u)) \geq \frac{n}{b+4} M_\sigma V(\sigma_1(wAw')) = \frac{n}{b+4}V(\sigma \sigma_1(wAw'))$$
	so
	\begin{equation}\label{largeRatio}
		\left\lvert \sigma \sigma_1(u) \right\rvert \geq \frac{n}{b+4} \left\lvert \sigma \sigma_1(wAw') \right\rvert
	\end{equation}
	Finally, \Cref{notFactor,largeRatio} mean that $\mathcal{F} (t)$ is not linearly recurrent.
	
	If (ii) holds, the same proof works by exchanging $1$ and $b$ and exchanging $A$ and $C$.
\end{proof}

\begin{proposition}\label{notLR}
	Over an even or odd alphabet, a smooth sequence not in $\mathscr{M}$ is not linearly recurrent.
\end{proposition}
\begin{proof}
	Let $x$ be a smooth sequence over an even or odd alphabet $\mathpzc{A} = \{1,b\}$ not in $\mathscr{M}$ and let $t := \mathcal{T}(x)$. Since $x = \pi_{x_0} \circ \mathcal{F} \circ S(t)$, it suffices to show that $\mathcal{F} \circ S(t)$ is linearly not recurrent: the fact that $x \notin \mathscr{M}$ yields $(1b)^n \sqsubset S(t)$ for all $n \geq 1$ and $S^{n+1}(t) \neq (1b)^\omega$ for all $n \geq 0$ so \Cref{FnotLR} yields the result.
\end{proof}

Then \Cref{LR,notLR} prove \Cref{URLR} (ii).

\section{Factor complexity of smooth sequences}\label{complexity}

In this section we prove \Cref{sameComplexity} and \Cref{linearCorollary}, based on a precise description of the bispecial factors of a given smooth sequence.

\subsection{Description of the bispecial factors}

In \cite{CH} we described the bispecial smooth words and their type (with respect to the language $\mathcal{C}_f^\infty$) using \textit{finite primitives}. Since the bispecial factors of a smooth sequence are also bispecial smooth words, we begin by recalling the description the bispecial smooth words (\Cref{smoothBispecials}). Then we adapt the steps of the proof of \Cref{smoothBispecials} to the context of this paper, thanks to \Cref{occurrences} and especially thanks to the forbidden and mandatory factors.

\subsubsection{Throwback to bispecial smooth words}

In \cite[Sections 3.1 to 3.3]{CH} we express bispecial smooth words with the following notation.

\begin{definition}
	Over a binary alphabet $\mathpzc{A} = \{a,b\}$, the \textup{finite primitives} are the maps
	\begin{align*}
		\mathcal{P}_{f,a}~:~~~~~~~~\mathpzc{A}^* & \longrightarrow \mathpzc{A}^+ & \mathcal{P}_{f,b}~: \mathpzc{A}^* & \longrightarrow \mathpzc{A}^+ \\
    	u_1...u_n & \longmapsto \begin{cases} a^ab^{u_1}a^{u_2}...a^{u_n} b^a \textit{ if } n \textit{ is even,} \\
   		a^ab^{u_1}a^{u_2}...b^{u_n} a^a \textit{ otherwise.} \end{cases} & u & \longmapsto \overline{\mathcal{P}_{f,a}(u)}
	\end{align*}
\end{definition}

\begin{example}
	Over $\{2,4\}$, we have $\mathcal{P}_{f,2}(4) = 22444422$ and $\mathcal{P}_{f,4}(42) = 4422224422$.
\end{example}

We make two important remarks.

\begin{remark}\label{prefixP}
	Let $u \in \mathpzc{A}^*$ and let $c \in \mathpzc{A}$.
	
	\textup{(i)} We have $\left\lvert \mathcal{P}_{f,c}(u) \right\rvert = a \lvert u \rvert_a + b \lvert u \rvert_b + 2a$.
	
	\textup{(ii)} We have $\mathcal{D}_f\left(\mathcal{P}_{f,c}(u)\right) = u$. We deduce that $u$ is smooth if and only if $\mathcal{P}_{f,c}(u)$ is smooth.
\end{remark}

We recall the description of bispecial smooth words and their type.

\begin{proposition}[{\cite{CH}}]\label{smoothBispecials}
	Over a binary alphabet $\mathpzc{A} = \{a,b\}$, the bispecial smooth words are the words of the form
	$$\mathcal{P}_{f,t_0} \mathcal{P}_{f,t_1} ... \mathcal{P}_{f,t_{i-1}} (c^n)$$
	where $i \geq 0$, $t_j \in \mathpzc{A}$ for all $j \in \llbracket 0,i-1\rrbracket$, $c \in \mathpzc{A}$ and $n \in \llbracket 0,b-1 \rrbracket$. Moreover, the strong (resp. weak) ones are the words such that $n = 0$ or $n = a < b-1$ (resp. such that $n = b-1 > a$).
\end{proposition}

\subsubsection{Short bispecial factors}

\begin{lemma}\label{shortBispecials}
	Let $x$ be a smooth sequence over an even or odd alphabet $\mathpzc{A} = \{a,b\}$ that has property $\mathfrak{P}3$ if $a = 1$. Then every short word (i.e., of the form $c^n$ where $c \in \mathpzc{A}$ and $n \in \llbracket 0,b-1 \rrbracket$) is a bispecial factor of $x$ and the strong (resp. weak) ones are the words such that $n = 0$ or $n = a$ (resp. such that $n = b-1$).
\end{lemma}
\begin{proof}
	Let $u$ be a short word over $\mathpzc{A}$, i.e., $u = c^n$ where $c \in \mathpzc{A}$ and $n \in \llbracket 0,b-1 \rrbracket$. Then \Cref{Sfactors} provides the bi-extensions of $u$ in $x$.
\begin{center}
	\begin{tabular}{|l|c|c|c|c|c|}
  		\hline
  		& $c c^n c$ & $c c^n \overline{c}$ & $\overline{c} c^n c$ & $\overline{c} c^n \overline{c}$ & type of $c^n$ \\
  		\hline
  		\hline
  		$n = 0$ or $n = a$ & short $\mathcal{D}_f$ & short $\mathcal{D}_f$ & short $\mathcal{D}_f$ & short $\mathcal{D}_f$ & strong \\
  		\hline
  		$n \in \llbracket 1,a-1 \rrbracket \cup \llbracket a+1,b-2 \rrbracket$ & short $\mathcal{D}_f$ & short $\mathcal{D}_f$ & short $\mathcal{D}_f$ & $\notin \mathpzc{C}_f^\infty$ & neutral \\
  		\hline
  		$n = b-1$ & $\notin \mathpzc{C}_f^\infty$ & short $\mathcal{D}_f$ & short $\mathcal{D}_f$ & $\notin \mathpzc{C}_f^\infty$ & weak \\
  		\hline
		\end{tabular}
	\end{center}
\end{proof}

\subsubsection{Long bispecial factors over even alphabets}

\begin{lemma}\label{evenBispecialsCascade}
	Let $x \in \mathcal{C}^\infty$ over an even alphabet $\mathpzc{A} = \{a,b\}$ and let $u \sqsubset x$.
	
	\textup{(i)} If $u$ is a long bispecial factor of $x$, then we have $u = \mathcal{P}_{f,u_1}\left(\mathcal{D}_f(u)\right)$. Furthermore, if $\mathcal{D}_f(u) \neq \varepsilon$ then $u_1 = \overline{x_0}$ and $\left\lvert \mathcal{D}_f(u) \right\rvert$ is even.
	
	\textup{(ii)} If $u$ is long and $c \in \mathpzc{A}$, then $u$ is a bispecial factor of $x$ if and only if $\mathcal{P}_{f,\overline{c}}(u)$ is a bispecial factor of $\mathcal{P}_c(x)$, and if they are bispecial then they have the same type.
\end{lemma}
\begin{proof}
	(i) Let ${a_1}^{p_1} ... {a_n}^{p_n}$ be the canonical factorization of $u$, where $n \geq 1$. If $n = 1$, $u$ is long so we necessarily have $u = {a_1}^b$ which cannot be bispecial since ${a_1}^b a_1 \notin \mathcal{C}_f^\infty$, therefore we have $n \geq 2$. Now, we have ${a_1}^{p_1+1} ... {a_n}^{p_n} = a_1 u \sqsubset x$ so $p_1 + 1 \leq b$, and $\overline{a_1} {a_1}^{p_1} ... {a_n}^{p_n} = \overline{a_1} u \sqsubset x$ so $p_1 \in \{a,b\}$, therefore $p_1 = a$; and symmetrically we obtain $p_n = a$. We deduce that $u = {a_1}^a {a_2}^{p_2} ... {a_{n-1}}^{p_{n-1}} a_n^a$ and $\mathcal{D}_f(u) = p_2 ... p_{n-1}$, which yields $u = \mathcal{P}_{f,u_1}\left(\mathcal{D}_f(u)\right)$ where $u_1 = a_1$.
	
	Now suppose that $\mathcal{D}_f(u) \neq \varepsilon$, which means that $n \geq 3$. If $u_1 = x_0$ and $p_2 = a$ (resp. $p_2 = b$), we observe that $f_1 = {x_0}^{a+1} \overline{x_0}^a x_0 \sqsubset_p x_0 u$ (resp. $f_2 = \overline{x_0} {x_0}^a \overline{x_0}^{a+1} \sqsubset_p \overline{x_0} u$) so \Cref{evenForbidden} yields $x_0 u \nsqsubset x$ (resp. $\overline{x_0} u \nsqsubset x$). This contradicts the fact that $u$ is bispecial, therefore we have $u_1 = \overline{x_0}$. Symmetrically, if $\left\lvert \mathcal{D}_f(u) \right\rvert$ is odd and $p_{n-1} = a$ (resp. $p_{n-1} = b$), we observe that $u = \overline{x_0}^a {x_0}^{p_2} ... {x_0}^{p_{n-1}} \overline{x_0}^a$ and $f_2 = \overline{x_0} {x_0}^a \overline{x_0}^{a+1} \sqsubset_s u \overline{x_0}$ (resp. $f_1 = {x_0}^{a+1} \overline{x_0}^a x_0 \sqsubset_s u x_0$) so \Cref{evenForbidden} yields $u \overline{x_0} \nsqsubset x$ (resp. $u x_0 \nsqsubset x$). This contradicts the fact that $u$ is bispecial, therefore $\left\lvert \mathcal{D}_f(u) \right\rvert$ is even.

	(ii) Set $n := \lvert u \rvert$. We have $\mathcal{D}_f\left(\mathcal{P}_{f,a}(u)\right) = u \sqsubset x$ so \Cref{DDuu} states that $\mathcal{P}_{f,a}(u) \sqsubset \mathcal{P}_c(x)$ or $\mathcal{P}_{f,b}(u) \sqsubset \mathcal{P}_c(x)$, i.e., there exists $d \in \mathpzc{A}$ such that $v := \mathcal{P}_{f,d}(u) = {a_0}^a {a_1}^{u_1} ... {a_n}^{u_n} {a_{n+1}}^a \sqsubset \mathcal{P}_c(x)$. Also, $\mathcal{D}_f(v) = u$ is long so \Cref{antiInvariance2} states that $\overline{v} \nsqsubset \mathcal{P}_c(x)$. Now we are going to prove the following equivalences:
	\begin{multicols}{2}
		\begin{enumerate}[1)]
			\item $aua \sqsubset x \iff \overline{a_0} v \overline{a_{n+1}} \sqsubset \mathcal{P}_c(x)$
			\item $aub \sqsubset x \iff \overline{a_0} v a_{n+1} \sqsubset \mathcal{P}_c(x)$
			
			\item $bua \sqsubset x \iff a_0 v \overline{a_{n+1}} \sqsubset \mathcal{P}_c(x)$
		
			\item $bub \sqsubset x \iff a_0 v a_{n+1} \sqsubset \mathcal{P}_c(x)$
			
			\item $au \sqsubset x \iff \overline{a_0} v \sqsubset \mathcal{P}_c(x)$
		
			\item $bu \sqsubset x \iff a_0 v \sqsubset \mathcal{P}_c(x)$
		
			\item $ua \sqsubset x \iff v \overline{a_{n+1}} \sqsubset \mathcal{P}_c(x)$
			
			\item $ub \sqsubset x \iff v a_{n+1} \sqsubset \mathcal{P}_c(x)$
		\end{enumerate}
	\end{multicols}
	We prove 1): if $aua \sqsubset x$, we observe that $\mathcal{D}_f(\overline{a_0} v \overline{a_{n+1}}) = aua$ and $\mathcal{D}(\mathcal{P}_c(x)) = x$ so \Cref{DDuu} states that $\overline{a_0} v \overline{a_{n+1}} \sqsubset \mathcal{P}_c(x)$ or $a_0 \overline{v} a_{n+1} \sqsubset \mathcal{P}_c(x)$, and $\overline{v} \nsqsubset \mathcal{P}_c(x)$ yields $\overline{a_0} v \overline{a_{n+1}} \sqsubset \mathcal{P}_c(x)$. Conversely, if $\overline{a_0} v \overline{a_{n+1}} \sqsubset \mathcal{P}_c(x)$, then the facts that $\mathcal{D}_f(\overline{a_0} v \overline{a_{n+1}}) = aua$ and that $\mathcal{D}(\mathcal{P}_c(x)) = x$ immediatly yield $aua \sqsubset x$. The proofs of 2), 3) and 4) are identical since $\mathcal{D}_f(\overline{a_0} v a_{n+1}) = aub$, $\mathcal{D}_f(a_0 v \overline{a_{n+1}}) = bua$ and $\mathcal{D}_f(a_0 v a_{n+1}) = bub$. Also, 5), 6), 7) and 8) directly follow from 1), 2), 3) and 4) (for 7) and 8) we need $x$ to be recurrent, which is provided by \Cref{URLR} (i)).
	
	Now suppose that $u$ is a bispecial factor of $x$. Then 5), 6), 7) and 8) imply that $v = \mathcal{P}_{f,d}(u)$ is a bispecial factor of $\mathcal{P}_c(x)$, and we have $\mathcal{D}_f(v) = u \neq \varepsilon$ so applying (i) to $v \sqsubset \mathcal{P}_c(x)$ ensures that $d = \overline{c}$ and $v = \mathcal{P}_{f,\overline{c}}(u)$. Moreover, 1), 2), 3) and 4) ensure that $u$ and $v$ have the same type.
	
	Conversely, if $\mathcal{P}_{f,\overline{c}}(u)$ is a bispecial factor of $\mathcal{P}_c(x)$, the fact that $\overline{v} = \mathcal{P}_{f,\overline{d}}(u) \nsqsubset \mathcal{P}_c(x)$ yields $\overline{d} \neq \overline{c}$ and $v = \mathcal{P}_{f,\overline{c}}(u)$. Then 5), 6), 7) and 8) imply that $u$ is a bispecial factor of $x$ and 1), 2), 3) and 4) ensure that $u$ and $v$ have the same type.
\end{proof}

\Cref{evenBispecialsCascade} yields a partial description of the bispecial factors of a given smooth sequence.

\begin{definition}
	Let $\mathpzc{A}$ be an even alphabet. If a word $u \in \mathpzc{A}^*$ is long but $\mathcal{D}_f(u)$ is short, we say that $u$ is a \textup{first long} word. Let $x \in \mathpzc{C}^\infty$ and let $u$ be a long bispecial factor of $x$. Then the \textup{degree} of $u$ is the first integer $i \geq 0$ such that $\left(\mathcal{D}_f\right)^i(u)$ is a first long factor of $\mathcal{D}^i(x)$, and we call $\mathcal{D}_f^i(u)$ the \textup{root} of $u$.
\end{definition}

\begin{lemma}\label{evenLongBispecialsForm}
	Let $x$ be a smooth sequence over an even alphabet $\mathpzc{A}$ and let $t := \mathcal{T}(x)$. Then the long bispecial factors of degree $i$ of $x$ are the words of the form
	$$\mathcal{P}_{f,\overline{t_0}} \mathcal{P}_{f,\overline{t_1}} ... \mathcal{P}_{f,\overline{t_{i-1}}} (u)$$
	where $u$ is a first long bispecial factor of $\mathcal{D}^i(x)$ and $t_j \in \mathpzc{A}$ for all $j \in \llbracket 0,i-1\rrbracket$. Moreover, each such bispecial factor in $x$ has the same type as its root $u$ in $\mathcal{D}^i(x)$.
\end{lemma}
\begin{proof}
	We proceed by induction on $i$. If $i = 0$, it is trivial. Now let $i \geq 1$ be such that \Cref{evenLongBispecialsForm} holds for $i$ and let $v$ be a long bispecial factor of degree $i+1$ of $x$. Then $v$ is long and $\mathcal{D}_f(v) \neq \varepsilon$ so \Cref{evenBispecialsCascade} (i) yields $v = \mathcal{P}_{f,\overline{x_0}}\left(\mathcal{D}_f(v)\right) = \mathcal{P}_{f,\overline{t_0}}\left(\mathcal{D}_f(v)\right)$ and \Cref{evenBispecialsCascade} (ii) applied to $\mathcal{D}_f(v) \sqsubset \mathcal{D}(x)$ implies that $\mathcal{D}_f(v)$ is a bispecial factor of $\mathcal{D}(x)$ and has the same type as $v$. Since $\mathcal{D}_f(v)$ has degree $i$ and $\mathcal{T}(\mathcal{D}(x)) = S(t)$, the induction hypothesis provides a first long bispecial factor $u$ of $\mathcal{D}^{i+1}(x)$ such that $\mathcal{D}_f(v) = \mathcal{P}_{f,\overline{t_1}} \mathcal{P}_{f,\overline{t_2}} ... \mathcal{P}_{f,\overline{t_i}} (u)$ has the same type as $u$. We deduce that $v = \mathcal{P}_{f,\overline{t_0}} \mathcal{P}_{f,\overline{t_1}} ... \mathcal{P}_{f,\overline{t_i}} (u)$ has the same type as $u$.
\end{proof}

It remains to find the first long strong bispecial factors of a given smooth sequence and their type.

\begin{lemma}\label{evenFirstLongBispecials}
	Let $x$ be a smooth sequence over an even alphabet $\mathpzc{A} = \{a,b\}$. Then the first long bispecial factors of $x$ are the words of the form
	$$\mathcal{P}_{f,d}(c^n)$$
	where $c \in \mathpzc{A}$, $n \in \llbracket 0,b-2\rrbracket$ is even, $d \in \mathpzc{A}$ and $d = \overline{x_0}$ if $n > 0$. Moreover, the strong (resp. weak) ones are the words such that $n = 0$ and $d = \overline{x_0}$, or $n = a < b-2$ (resp. such that $n = 0$ and $d = x_0$, or $n = b-2 > a$).
\end{lemma}
\begin{proof}
	\Cref{evenBispecialsCascade} (i) states that the first long bispecial factors of $x$ are words of the form $\mathcal{P}_{f,d}(c^n)$ where $c \in \mathpzc{A}$, $n \in \llbracket 0,b-2\rrbracket$ is even, $d \in \mathpzc{A}$ and $d = \overline{x_0}$ if $n > 0$. Conversely, let us inspect the possible bi-extensions of a word $u := \mathcal{P}_{f,d}(c^n)$ of this form in $x$.
	
	If $n = 0$ and $d = \overline{x_0}$, we get $u = \overline{x_0}^a {x_0}^a$. Then we observe that $x_0 u x_0 \sqsubset_s m_1$, $x_0 u \overline{x_0} \sqsubset m_1$, $\overline{x_0} u x_0 \sqsubset m_3$ and $\overline{x_0} u \overline{x_0} \sqsubset_p m_1$, so with \Cref{evenForbidden,evenMandatory} we deduce that $u$ is a strong bispecial factor of $x$.
	
	If $n = 0$ and $d = x_0$, we get $u = {x_0}^a \overline{x_0}^a$. Then we observe that $x_0 u x_0 = f_1$, $x_0 u \overline{x_0} \sqsubset m_3$, $\overline{x_0} u x_0 \sqsubset m_1$ and $\overline{x_0} u \overline{x_0} = f_2$, so with \Cref{evenForbidden,evenMandatory} we deduce that $u$ is a weak bispecial factor of $x$.

	If $n \in \llbracket 2, b-2 \rrbracket$ is even and $c = a$, we get $d = \overline{x_0}$ and $u = \overline{x_0}^a \left[{x_0}^a \overline{x_0}^a\right]^\frac{n}{2} {x_0}^a$. Then, with \Cref{evenForbidden,evenMandatory}, we obtain the bi-extensions and the type of $u$ in $x$:
	\begin{center}
		\begin{tabular}{|l|c|c|c|c|c|}
		\hline
		  & $x_0 u x_0$ & $x_0 u \overline{x_0}$ & $\overline{x_0} u x_0$ & $\overline{x_0} u \overline{x_0}$ & type of $u$ \\
  		\hline
  		\hline
  		$n \in \llbracket 2,a-2 \rrbracket \cup \llbracket a+2,b-4 \rrbracket$ & $\sqsubset_s m_1$ & $\sqsubset m_1$ & $\notin \mathpzc{C}_f^\infty$ & $\sqsubset_p m_1$ & neutral \\
  		\hline
  		$n = a = b-2$ & $\sqsubset_s m_1$ & $= f_3$ & $= m_2$ & $\sqsubset_p m_1$ & neutral \\
  		\hline
  		$n = a < b-2$ & $\sqsubset_s m_1$ & $\sqsubset m_1$ & $= m_2$ & $\sqsubset_p m_1$ & strong \\
  		\hline
  		$n = b-2 > a$ & $\sqsubset_s m_1$ & $= f_3$ & $\notin \mathpzc{C}_f^\infty$ & $\sqsubset_p m_1$ & weak \\
  		\hline
		\end{tabular}
	\end{center}
	If $n \in \llbracket 2, b-2 \rrbracket$ is even and $c = b$, we get $d = \overline{x_0}$ and $u = \overline{x_0}^a \left[{x_0}^b \overline{x_0}^b\right]^\frac{n}{2} {x_0}^a$. Then, with \Cref{evenForbidden,evenMandatory} we obtain the bi-extensions and the type of $u$ in $x$:
	\begin{center}
		\begin{tabular}{|l|c|c|c|c|c|}
		\hline
		  & $x_0 u x_0$ & $x_0 u \overline{x_0}$ & $\overline{x_0} u x_0$ & $\overline{x_0} u \overline{x_0}$ & type of $u$ \\
  		\hline
		\hline
  		$n \in \llbracket 2,a-2 \rrbracket \cup \llbracket a+2,b-4 \rrbracket$ & $\sqsubset_s m_3$ & $\notin \mathpzc{C}_f^\infty$ & $\sqsubset m_3$ & $\sqsubset_p m_3$ & neutral \\
  		\hline
  		$n = a = b-2$ & $\sqsubset_s m_3$ & $= m_4$ & $= f_4$ & $\sqsubset_p m_3$ & neutral \\
  		\hline
  		$n = a < b-2$ & $\sqsubset_s m_3$ & $= m_4$ & $\sqsubset m_3$ & $\sqsubset_p m_3$ & strong \\
  		\hline
  		$n = b-2 > a$ & $\sqsubset_s m_3$ & $\notin \mathpzc{C}_f^\infty$ & $= f_4$ & $\sqsubset_p m_3$ & weak \\
  		\hline
		\end{tabular}
	\end{center}
\end{proof}

Putting the pieces together, \Cref{evenLongBispecialsForm,evenFirstLongBispecials} yield the description of the strong and weak bispecial factors of any smooth sequence.

\begin{corollary}\label{evenLongBispecials}
	Let $x$ be a smooth sequence over an even alphabet $\mathpzc{A} = \{a,b\}$ and let $t := \mathcal{T}(x)$. Then the long bispecial factors of $x$ are the words of the form
	$$\mathcal{P}_{f,\overline{t_0}} \mathcal{P}_{f,\overline{t_1}} ... \mathcal{P}_{f,\overline{t_{i-1}}} \mathcal{P}_{f,d} (c^n)$$
	where $i \geq 0$, $c \in \mathpzc{A}$, $n \in \llbracket 0,b-2\rrbracket$ is even, $d \in \mathpzc{A}$ and $d = \overline{t_i}$ if $n > 0$. Moreover, the strong (resp. weak) ones are the words such that $n = 0$ and $d = \overline{t_i}$, or $n = a < b-2$ (resp. such that $n = 0$ and $d = t_i$, or $n = b-2 > a$).
\end{corollary}

\subsubsection{Long bispecial factors over odd alphabets}

\begin{lemma}\label{oddBispecialsCascade}
	Let $x \in \mathcal{C}^\infty$ over an odd alphabet $\mathpzc{A} = \{a,b\}$ and let $u \sqsubset x$.
	
	\textup{(i)} If $u$ is a long bispecial factor of $x$, then we have $u = \mathcal{P}_{f,u_1}\left(\mathcal{D}_f(u)\right)$. Furthermore, if $\mathcal{D}_f(u) \neq \varepsilon$ then
	$$u_1 = \begin{cases} \overline{x_0} \textit{ if } \mathcal{D}_f(u)_1 = \mathcal{D}(x)_0,\\ x_0 \textit{ otherwise} \end{cases}.$$
	
	\textup{(ii)} If $u$ is long and $c \in \mathpzc{A}$, then $u$ is a bispecial factor of $x$ if and only if $\mathcal{P}_{f,d}(u)$ is a bispecial factor of $\mathcal{P}_c(x)$ where $d := \begin{cases} \overline{c} \textit{ if } u_1 = x_0,\\ c \textit{ otherwise} \end{cases}$, and if they are bispecial then they have the same type.
\end{lemma}
\begin{proof}
	(i) Let ${a_1}^{p_1} ... {a_n}^{p_n}$ be the canonical factorization of $u$, where $n \geq 1$. We show that $u = \mathcal{P}_{f,u_1}\left(\mathcal{D}_f(u)\right)$ exactly as in the proof of \Cref{evenBispecialsCascade} (i), in particular we have $u_1 = a_1$, $u = {a_1}^a {a_2}^{p_2} ... {a_{n-1}}^{p_{n-1}} a_n^a$ and $\mathcal{D}_f(u) = p_2 ... p_{n-1}$. Now suppose that $\mathcal{D}_f(u) \neq \varepsilon$, which means that $n \geq 3$. Thanks to \Cref{oddForbidden}, we obtain the value of $u_1$ by showing that the other possibility contradicts the fact that $u$ is a bispecial factor of $x$.
	\begin{center}
		\begin{tabular}{|l|c|c|c|c|}
		\hline
		$p_2 = \mathcal{D}_f(u)_1$ & $\mathcal{D}(x)_0$ & $u_1$ & contradiction \\
  		\hline
  		\hline
  		$a$ & $a$ & $x_0$ & $f_2 = {x_0}^{a+1} \overline{x_0}^a x_0 \sqsubset_p x_0 u \nsqsubset x$ \\
  		\hline
  		$a$ & $b$ & $\overline{x_0}$ & $\overline{f_2} = \overline{x_0}^{a+1} {x_0}^a \overline{x_0} \sqsubset_p \overline{x_0} u \nsqsubset x$ \\
		\hline
		$b$ & $a$ & $\overline{x_0}$ & $f_1 = x_0 \overline{x_0}^a {x_0}^{a+1} \sqsubset_p x_0 u \nsqsubset x$ \\
		\hline
		$b$ & $b$ & $x_0$ & $\overline{f_1} = \overline{x_0} {x_0}^a \overline{x_0}^{a+1} \sqsubset_p \overline{x_0} u \nsqsubset x$ \\
		\hline
		\end{tabular}
	\end{center}

	(ii) The proof it identical to \Cref{evenBispecialsCascade} (ii) by using (i).
\end{proof}

\Cref{oddBispecialsCascade} yields a partial description of the bispecial factors of a given smooth sequence.

\begin{definition}
	Let $\mathpzc{A}$ be an odd alphabet. If a word $u \in \mathpzc{A}^*$ is long but $\mathcal{D}_f(u)$ is short, we say that $u$ is a \textup{first long} word. Let $x \in \mathpzc{C}^\infty$ and let $u$ be a long bispecial factor of $x$. Then the \textup{degree} of $u$ is the first integer $i \geq 0$ such that $\mathcal{D}_f^i(u)$ is a first long factor of $\mathcal{D}^i(x)$, and we call $\mathcal{D}_f^i(u)$ the \textup{root} of $u$.
\end{definition}

\begin{lemma}\label{oddLongBispecialsForm}
	Let $x \in \mathpzc{C}^\infty$ over an odd alphabet $\mathpzc{A}$ and let $t := \mathcal{T}(x)$. Then the long bispecial factors of degree $i$ of $x$ are the words of the form
	$$\mathcal{P}_{f,t^i_0} \mathcal{P}_{f,t^i_1} ... \mathcal{P}_{f,t^i_{i-1}} (u)$$
	where $u$ is a first long bispecial factor of $\mathcal{D}^i(x)$ and $t^i_j := \begin{cases} \overline{t_j} \textit{ if } u_1 = t_i \textit{ and } i-j \textit{ is odd, or if } u_1 = \overline{t_i} \textit{ and } i-j \textit{ is even,} \\ t_j \textit{ otherwise.} \end{cases}$ Moreover, each such bispecial factor in $x$ has the same type as its root $u$ in $\mathcal{D}^i(x)$.
\end{lemma}
\begin{proof}
	We proceed by induction on $i$. If $i = 0$, it is trivial. Now let $i \geq 1$ be such that \Cref{oddLongBispecialsForm} holds for $i$ and let $v$ be a long bispecial factor of degree $i+1$ of $x$. Then $v$ is long and $\mathcal{D}_f(v) \neq \varepsilon$ so \Cref{oddBispecialsCascade} (i) yields $v = \mathcal{P}_{f,d}\left(\mathcal{D}_f(v)\right)$ where $d = \begin{cases} \overline{x_0} \textit{ if } \mathcal{D}_f(v)_1 = \mathcal{D}(x)_0,\\ x_0 \textup{ otherwise} \end{cases}$ and \Cref{oddBispecialsCascade} (ii) applied to $\mathcal{D}_f(v) \sqsubset \mathcal{D}(x)$ implies that $\mathcal{D}_f(v)$ is a bispecial factor of $\mathcal{D}(x)$ and has the same type as $v$. Since $\mathcal{D}_f(v)$ has degree $i$ and $\mathcal{T}(\mathcal{D}(x)) = S(t)$, the induction hypothesis provides a first long bispecial factor $u$ of $\mathcal{D}^{i+1}(x)$ such that $\mathcal{D}_f(v) = \mathcal{P}_{f,t^i_0} \mathcal{P}_{f,t^i_1} ... \mathcal{P}_{f,t^i_{i-1}} (u)$ has the same type as $u$ where
	$$t^i_j = \begin{cases} \overline{t_{j+1}} \textit{ if } u_1 = t_{i+1} \textit{ and } i-j \textit{ is odd, or if } u_1 = \overline{t_{i+1}} \textit{ and } i-j \textit{ is even,} \\ t_{j+1} \textit{ otherwise} \end{cases}.$$
	We deduce that $v = \mathcal{P}_{f,d} \mathcal{P}_{f,t^i_0} \mathcal{P}_{f,t^i_1} ... \mathcal{P}_{f,t^i_{i-1}} (u)$ has the same type as $u$. By setting
	$$t^{i+1}_j := \begin{cases} \overline{t_j} \textit{ if } u_1 = t_{i+1} \textit{ and } i+1-j \textit{ is odd, or if } u_1 = \overline{t_{i+1}} \textit{ and } i+1-j \textit{ is even,} \\ t_j \textit{ otherwise} \end{cases},$$
	we immediatly get $t^{i+1}_j = t^i_{j-1}$ for all $j \geq 1$ and we have
	\begin{align*}
		t^{i+1}_0 = &~ \begin{cases} \overline{t_0} \textit{ if } u_1 = t_{i+1} \textit{ and } i+1 \textit{ is odd, or if } u_1 = \overline{t_{i+1}} \textit{ and } i+1 \textit{ is even,} \\ t_0 \textit{ otherwise} \end{cases} \\
		= &~ \begin{cases} \overline{x_0} \textit{ if } u_1 = t_{i+1} \textit{ and } i \textit{ is even, or if } u_1 = \overline{t_{i+1}} \textit{ and } i \textit{ is odd,} \\ x_0 \textit{ otherwise} \end{cases} \\
		= &~ \begin{cases} \overline{x_0} \textit{ if } t^i_0 = t_1, \\ x_0 \textit{ otherwise} \end{cases} = \begin{cases} \overline{x_0} \textit{ if } \mathcal{D}_f(v)_1 = \mathcal{D}(x)_0,\\ x_0 \textit{ otherwise} \end{cases} = d.
	\end{align*}
	We deduce that $v = \mathcal{P}_{f,t^{i+1}_0} \mathcal{P}_{f,t^{i+1}_1} ... \mathcal{P}_{f,t^{i+1}_i} (u)$.
\end{proof}

It remains to find the first long strong (resp. weak) bispecial factors of a given smooth sequence.

\begin{lemma}\label{oddFirstLongBispecials}
	Let $x$ be a smooth sequence over an odd alphabet $\mathpzc{A} = \{a,b\}$ that has property $\mathfrak{P}3$ if $a = 1$. Then the first long bispecial factors of $x$ are the words of the form
	$$\mathcal{P}_{f,d}(c^n)$$
	where $c \in \mathpzc{A}$, $n \in \llbracket 0,b-1 \rrbracket$, $d \in \mathpzc{A}$, and 
	$n$ is odd and $d = \begin{cases} \overline{x_0} \textit{ if } c = \mathcal{D}(x)_0, \\ \textit{otherwise} \end{cases}$ if $n > 0$. Moreover, the strong (resp. weak) ones are the words such that $n = a < b-2$ (resp. such that $n = b-2 > a$).
\end{lemma}
\begin{proof}
	\Cref{oddBispecialsCascade} (i) states that the first long bispecial factors of $x$ are words of the form $\mathcal{P}_{f,d}(c^n)$ where $c \in \mathpzc{A}$, $n \in \llbracket 0,b-1\rrbracket$, $d \in \mathpzc{A}$ and $d = \begin{cases} \overline{x_0} \textit{ if } c = \mathcal{D}(x)_0,\\ x_0 \textit{ otherwise} \end{cases}$ if $n > 0$. If $n > 0$, we also show that $n$ must be odd: if $n > 0$ is even, we have $\mathcal{P}_{f,d}(c^n) = d^a \left[\overline{d}^c d^c\right]^\frac{n}{2} \overline{d}^a$ and we show that it leads to a contradiction thanks to \Cref{oddForbidden}:
	\begin{center}
		\begin{tabular}{|l|c|c|c|c|}
		\hline
		$c$ & $\mathcal{D}(x)_0$ & $d$ & contradiction \\
  		\hline
  		\hline
  		$a$ & $a$ & $\overline{x_0}$ & $f_1 = x_0 \overline{x_0}^a {x_0}^{a+1} \sqsubset_s u x_0 \nsqsubset x$ \\
  		\hline
  		$a$ & $b$ & $x_0$ & $\overline{f_1} = \overline{x_0} {x_0}^a \overline{x_0}^{a+1} \sqsubset_s u \overline{x_0} \nsqsubset x$ \\
		\hline
		$b$ & $a$ & $x_0$ & $f_2 = {x_0}^{a+1} \overline{x_0}^a x_0 \sqsubset_s u x_0 \nsqsubset x$ \\
		\hline
		$b$ & $b$ & $\overline{x_0}$ & $\overline{f_2} = \overline{x_0}^{a+1} {x_0}^a \overline{x_0} \sqsubset_s u \overline{x_0} \nsqsubset x$ \\
		\hline
		\end{tabular}
	\end{center}
	Conversely, let us inspect the possible bi-extensions of a word $u := \mathcal{P}_{f,d}(c^n)$ of this form in $x$.
	
	If $n = 0$, $d = x_0$ and $\mathcal{D}(x)_0 = a$ (resp. $\mathcal{D}(x)_0 = b$), we get $u = {x_0}^a \overline{x_0}^a$. Then we observe that $x_0 u x_0 = f_2$, $x_0 u \overline{x_0} \sqsubset m_3$, $\overline{x_0} u x_0 \sqsubset m_1$ and $\overline{x_0} u \overline{x_0} \sqsubset_s m_1$ (resp. $x_0 u x_0 \sqsubset_p \overline{m_1}$, $x_0 u \overline{x_0} \sqsubset \overline{m_3}$, $\overline{x_0} u x_0 \sqsubset \overline{m_1}$ and $\overline{x_0} u \overline{x_0} = \overline{f_1}$), so with \Cref{oddForbidden,oddMandatory} we deduce that $u$ is a neutral bispecial factor of $x$.
	
	If $n = 0$, $d = \overline{x_0}$ and $\mathcal{D}(x)_0 = a$ (resp. $\mathcal{D}(x)_0 = b$), we get $u = {x_0}^a \overline{x_0}^a$. Then we observe that	$x_0 u x_0 = f_1$, $x_0 u \overline{x_0} \sqsubset m_1$, $\overline{x_0} u x_0 \sqsubset m_3$ and $\overline{x_0} u \overline{x_0} \sqsubset_p m_1$ (resp. $x_0 u x_0 \sqsubset_s \overline{m_1}$, $x_0 u \overline{x_0} \sqsubset \overline{m_1}$, $\overline{x_0} u x_0 \sqsubset \overline{m_3}$ and $\overline{x_0} u \overline{x_0} = \overline{f_2}$) so with \Cref{oddForbidden,oddMandatory} we deduce that $u$ is a neutral bispecial factor of $x$.
	
	If $n \in \llbracket 1, b-2 \rrbracket$ is odd, $c = a$ and $\mathcal{D}(x)_0 = a$, we get $d = \overline{x_0}$ and $u = \overline{x_0}^a \left[{x_0}^a \overline{x_0}^a\right]^\frac{n-1}{2} {x_0}^a \overline{x_0}^a$. Then, with \Cref{oddForbidden,oddMandatory}, we obtain the bi-extensions and the type of $u$ in $x$:
	\begin{center}
		\begin{tabular}{|l|c|c|c|c|c|}
		\hline
		 & $x_0 u x_0$ & $x_0 u \overline{x_0}$ & $\overline{x_0} u x_0$ & $\overline{x_0} u \overline{x_0}$ & type of $u$ \\
  		\hline
  		$n \in \llbracket 1,a-2 \rrbracket \cup \llbracket a+2,b-4 \rrbracket$ & $\sqsubset m_1$ & $\sqsubset_s m_1$ & $\sqsubset_p m_1$ & $\notin \mathpzc{C}_f^\infty$ & neutral \\
  		\hline
  		$n = a = b-2$ & $= f_3$ & $\sqsubset_s m_1$ & $\sqsubset_p m_1$ & $= m_2$ & neutral \\
  		\hline
  		$n = a < b-2$ & $\sqsubset m_1$ & $\sqsubset_s m_1$ & $\sqsubset_p m_1$ & $= m_2$ & strong \\
  		\hline
  		$n = b-2 > a$ & $= f_3$ & $\sqsubset_s m_1$ & $\sqsubset_p m_1$ & $\notin \mathcal{C}_f^\infty$ & weak \\
  		\hline
		\end{tabular}
	\end{center}
	If $n \in \llbracket 1, b-2 \rrbracket$ is odd, $c = a$ and $\mathcal{D}(x)_0 = b$, we get $d = x_0$ and $u = {x_0}^a \left[\overline{x_0}^a {x_0}^a\right]^\frac{n-1}{2} \overline{x_0}^a {x_0}^a$. Then, with \Cref{oddForbidden,oddMandatory}, we obtain the bi-extensions and the type of $u$ in $x$:
	\begin{center}
		\begin{tabular}{|l|c|c|c|c|c|}
		\hline
		 & $x_0 u x_0$ & $x_0 u \overline{x_0}$ & $\overline{x_0} u x_0$ & $\overline{x_0} u \overline{x_0}$ & type of $u$ \\
  		\hline
		$n \in \llbracket 1,a-2 \rrbracket \cup \llbracket a+2,b-4 \rrbracket$ & $\notin \mathcal{C}_f^\infty$ & $\sqsubset_p \overline{m_1}$ & $\sqsubset_s \overline{m_1}$ & $\sqsubset \overline{m_1}$ & neutral \\
  		\hline
  		$n = a = b-2$ & $= \overline{m_2}$ & $\sqsubset_p \overline{m_1}$ & $\sqsubset_s \overline{m_1}$ & $= \overline{f_3}$ & neutral \\
  		\hline
  		$n = a < b-2$ & $= \overline{m_2}$ & $\sqsubset_p \overline{m_1}$ & $\sqsubset_s \overline{m_1}$ & $\sqsubset \overline{m_1}$ & strong \\
  		\hline
  		$n = b-2 > a$ & $\notin \mathcal{C}_f^\infty$ & $\sqsubset_p \overline{m_1}$ & $\sqsubset_s \overline{m_1}$ & $= \overline{f_3}$ & weak \\
  		\hline
		\end{tabular}
	\end{center}
	If $n \in \llbracket 1, b-2 \rrbracket$ is odd, $c = b$ and $\mathcal{D}(x)_0 = a$, we get $d = x_0$ and $u = {x_0}^a \left[\overline{x_0}^b {x_0}^b\right]^\frac{n-1}{2} \overline{x_0}^b {x_0}^a$. Then, with \Cref{oddForbidden,oddMandatory}, we obtain the bi-extensions and the type of $u$ in $x$:
	\begin{center}
		\begin{tabular}{|l|c|c|c|c|c|}
		\hline
		 & $x_0 u x_0$ & $x_0 u \overline{x_0}$ & $\overline{x_0} u x_0$ & $\overline{x_0} u \overline{x_0}$ & type of $u$ \\
  		\hline
  		$n \in \llbracket 1,a-2 \rrbracket \cup \llbracket a+2,b-4 \rrbracket$ & $\sqsubset m_3$ & $\sqsubset_s m_3$ & $\sqsubset_p m_3$ & $\notin \mathcal{C}_f^\infty$ & neutral \\
  		\hline
  		$n = a = b-2$ & $= f_4$ & $\sqsubset_s m_3$ & $\sqsubset_p m_3$ & $= m_4$ & neutral \\
  		\hline
  		$n = a < b-2$ & $\sqsubset m_3$ & $\sqsubset_s m_3$ & $\sqsubset_p m_3$ & $= m_4$ & strong \\
  		\hline
  		$n = b-2 > a$ & $= f_4$ & $\sqsubset_s m_3$ & $\sqsubset_p m_3$ & $\notin \mathcal{C}_f^\infty$ & weak \\
  		\hline
		\end{tabular}
	\end{center}
	If $n \in \llbracket 1, b-2 \rrbracket$ is odd, $c = b$ and $\mathcal{D}(x)_0 = b$, we get $d = \overline{x_0}$ and $u = \overline{x_0}^a \left[{x_0}^b \overline{x_0}^b\right]^\frac{n-1}{2} {x_0}^b \overline{x_0}^a$. Then, with \Cref{oddForbidden,oddMandatory}, we obtain the bi-extensions and the type of $u$ in $x$:
	\begin{center}
		\begin{tabular}{|l|c|c|c|c|c|}
		\hline
		 & $x_0 u x_0$ & $x_0 u \overline{x_0}$ & $\overline{x_0} u x_0$ & $\overline{x_0} u \overline{x_0}$ & type of $u$ \\
  		\hline
  		$n \in \llbracket 1,a-2 \rrbracket \cup \llbracket a+2,b-4 \rrbracket$ & $\notin \mathcal{C}_f^\infty$ & $\sqsubset_p \overline{m_3}$ & $\sqsubset_s \overline{m_3}$ & $\sqsubset \overline{m_3}$ & neutral \\
  		\hline
  		$n = a = b-2$ & $= \overline{m_4}$ & $\sqsubset_p \overline{m_3}$ & $\sqsubset_s \overline{m_3}$ & $= \overline{f_4}$ & neutral \\
  		\hline
  		$n = a < b-2$ & $= \overline{m_4}$ & $\sqsubset_p \overline{m_3}$ & $\sqsubset_s \overline{m_3}$ & $\sqsubset \overline{m_3}$ & strong \\
  		\hline
  		$n = b-2 > a$ & $\notin \mathcal{C}_f^\infty$ & $\sqsubset_p \overline{m_3}$ & $\sqsubset_s \overline{m_3}$ & $= \overline{f_4}$ & weak \\
  		\hline
		\end{tabular}
	\end{center}
\end{proof}

Putting the pieces together, \Cref{oddLongBispecialsForm,oddFirstLongBispecials} yield the description of the strong and weak bispecial factors of any smooth sequence.

\begin{corollary}\label{oddLongBispecials}
	Let $x$ be a smooth sequence over an odd alphabet $\mathpzc{A} = \{a,b\}$ that does not have property $\mathfrak{P}2$ and let $t := \mathcal{T}(x)$. Then the long bispecial factors of $x$ are the words of one of the following forms:
	$$\mathcal{P}_{f,t^i_0} \mathcal{P}_{f,t^i_1} ... \mathcal{P}_{f,t^i_{i-1}} \mathcal{P}_{f,d} (\varepsilon)$$
	where $i \geq 0$, $d \in \mathpzc{A}$ and $t^i_j := \begin{cases} \overline{t_j} \textit{ if } d = t_i \textit{ and } i-j \textit{ is odd, or if } d = \overline{t_i} \textit{ and } i-j \textit{ is even,} \\ t_j \textit{ otherwise.} \end{cases}$;
	or
	$$\mathcal{P}_{f,t^i_0} \mathcal{P}_{f,t^i_1} ... \mathcal{P}_{f,t^i_i} (c^n)$$
	where $i \geq 0$, $c \in \mathpzc{A}$, $n \in \llbracket 1,b-1 \rrbracket$ is odd and $t^i_j := \begin{cases} \overline{t_j} \textit{ if } t_{i+1} = \overline{c} \textit{ and } i-j \textit{ is odd, or if } t_{i+1} = c \textit{ and } i-j \textit{ is even,} \\ t_j \textit{ otherwise} \end{cases}$.
	Moreover, the strong (resp. weak) ones are the words of the second form such that $n = a < b-2$ (resp. $n = b-2 > a$).
\end{corollary}

\subsection{Computation of the factor complexity}

In this subsection we prove \Cref{sameComplexity}.

\subsubsection{Smooth sequences with the same factor complexity}

In order to prove \Cref{sameComplexity} (i), we describe the length of the primitives over even alphabets.

\begin{lemma}\label{evenLength}
	Let $\mathpzc{A} = \{a,b\}$ be an even alphabet. If $u \in \mathpzc{A}^*$ has even length and $t \in \mathpzc{A}^\mathbb{N}$, then the sequence $(\ell_n)_{n \geq 0} := \left(\left\lvert \mathcal{P}_{f,t_0} \mathcal{P}_{f,t_1} ... \mathcal{P}_{f,t_n} (u) \right\rvert\right)_{n \geq 0}$ satisfies
	\begin{equation*}
		\begin{cases} \ell_0 = a \lvert u \rvert_a + b \lvert u \rvert_b + 2a \\ \ell_{n+1} = \frac{a+b}{2} \ell_n + 2a \textit{ for all } n \geq 0.
		\end{cases}
	\end{equation*}
\end{lemma}
\begin{proof}
	We first observe that, if $v \in \mathpzc{A}^*$ and $c \in \mathpzc{A}$, then
	\begin{equation}\label{Puaub}
		\left\lvert \mathcal{P}_{f,c}(v) \right\rvert = a \lvert v \rvert_a + b \lvert v \rvert_b + 2a
	\end{equation}
	which already provides the formula for $\ell_0$. Also, if $v$ has even length, then $\left\lvert \mathcal{P}_{f,c}(v) \right\rvert_a = \left\lvert \mathcal{P}_{f,c}(v) \right\rvert_b$. We then obtain by immediate induction that
	$$\left\lvert \mathcal{P}_{f,t_0} \mathcal{P}_{f,t_1} ... \mathcal{P}_{f,t_n} (u) \right\rvert_a = \left\lvert \mathcal{P}_{f,t_0} \mathcal{P}_{f,t_1} ... \mathcal{P}_{f,t_n} (u) \right\rvert_b$$
	for all $t \in \mathpzc{A}^\mathbb{N}$ and all $n \geq 0$, which yields
	\begin{equation}\label{uaubu}
		\left\lvert \mathcal{P}_{f,t_0} \mathcal{P}_{f,t_1} ... \mathcal{P}_{f,t_n} (u) \right\rvert_a = \left\lvert \mathcal{P}_{f,t_0} \mathcal{P}_{f,t_1} ... \mathcal{P}_{f,t_n} (u) \right\rvert_b = \frac{\ell_n}{2}
	\end{equation}
	for all $t \in \mathpzc{A}^\mathbb{N}$ and all $n \geq 0$. Finally, if $t \in \mathpzc{A}^\mathbb{N}$ and $n \geq 0$, we get
	\begin{align*}
		\ell_{n+1} = &~ \left\lvert \mathcal{P}_{f,t_0} \left(\mathcal{P}_{f,t_1} ... \mathcal{P}_{f,t_{n+1}} (u)\right) \right\rvert \\
		= &~ a\left\lvert \mathcal{P}_{f,t_1} \mathcal{P}_{f,t_2} ... \mathcal{P}_{f,t_{n+1}} (u) \right\rvert_a + b\left\lvert \mathcal{P}_{f,t_1} \mathcal{P}_{f,t_2} ... \mathcal{P}_{f,t_{n+1}} (u) \right\rvert_b + 2a \textrm{ with \Cref{Puaub}} \\
		= &~ \frac{a+b}{2} \ell_n + 2a \textrm{ with \Cref{uaubu}}
	\end{align*}
	for all $n \geq 0$.
\end{proof}

This allows to show that all smooth sequences over an even alphabet have the same factor complexity, this is \Cref{sameComplexity} (i).

\begin{proof}[Proof of \Cref*{sameComplexity} \textup{(i)}]
	Let $\mathpzc{A} = \{a,b\}$ be an even alphabet (resp. an odd alphabet such that $b = a+2$) and let $x$ be a smooth sequence (resp. let $x$ be a smooth sequence that does not have property $\mathfrak{P}2$ if $a = 1$). Then \Cref{shortBispecials} implies that the short strong and weak bispecial factors of $x$ do not depend on $x$. Also, \Cref{evenLength} implies that the length of the long bispecial factors of $x$ (cf. \Cref{evenLongBispecials}) does not depend on $\mathcal{T}(x)$, i.e., it does not depend on $x$ (resp. \Cref{oddLongBispecials} states that $x$ has no long strong (resp. weak) bispecial factor). With \Cref{pfrombsw}, we deduce that the second finite difference $b(n)$ does not depend on $x$, which implies that every such smooth sequence has the same factor complexity $p(n)$.
\end{proof}

\subsubsection{Factor complexity over tight alphabets}

By \textit{tight} alphabet we mean an alphabet of the form $\{a,a+2\}$. This is where we prove here \Cref{sameComplexity} (ii).

\begin{proof}[Proof of \Cref*{sameComplexity} \textup{(ii)}]
	Let $\mathpzc{A} = \{a,b\}$ be an even or odd alphabet such that $b = a+2$. Since the smooth sequence $x := \mathpzc{E}(a^\omega)$ does not have property $\mathfrak{P}2$ if $a = 1$, \Cref{sameComplexity} (i) states $p(n) = p_x(n)$. If $\mathpzc{A}$ is even, \Cref{evenLongBispecials} states that the long strong (resp. weak) bispecial factors of $x$ are $\left(\mathcal{P}_{f,b}\right)^{n+1}(\varepsilon)$ (resp. $\left(\mathcal{P}_{f,b}\right)^n \mathcal{P}_{f,a}(\varepsilon)$) for all $n \geq 0$, and \Cref{evenLength} implies that the words $\mathcal{P}_{f,b}^{n+1}(\varepsilon)$ and $\mathcal{P}_{f,b}^n \mathcal{P}_{f,a}(\varepsilon)$ have the same length for all $n \geq 0$. Since one is strong and the other is weak, their cumulated contribution to the factor complexity (cf. \Cref{pfrombsw}) is zero. Similarly, if $\mathpzc{A}$ is odd, \Cref{oddLongBispecials} states that $x$ has no strong (resp. weak) bispecial factor. This leaves only the contribution of the short strong (resp. weak) bispecial factors $\varepsilon$, $a^a$ and $b^a$ ($a^{a+1}$ and $b^{a+1}$) (cf. \Cref{shortBispecials}), and we compute
	\begin{align*}
		b(n) = \begin{cases} 1 \textrm{ if } n = 0 \\ 0 \textrm{ if } n \in \llbracket 1,a-1 \rrbracket \\ 2 \textrm{ if } n = a \\ -2 \textrm{ if } n = a+1 \\ 0 \textrm{ if } n \geq a+2 \end{cases} \textrm{,~~~} & s(n) = \begin{cases} 1 \textrm{ if } n = 0 \\ 2 \textrm{ if } n \in \llbracket 1,a \rrbracket \\ 4 \textrm{ if } n = a+1 \\ 2 \textrm{ if } n \geq a+2 \end{cases}
	\end{align*}
	and
	\begin{equation*}
		p(n) = p_x(n) = \begin{cases} 1 \textrm{ if } n = 0 \\ 2n \textrm{ if } n \in \llbracket 1,a+1 \rrbracket \\ 2n+2 \textrm{ if } n \geq a+2 \end{cases}.
	\end{equation*}
\end{proof}

\subsubsection{Factor complexity over loose even alphabets}

By \textit{loose} alphabet we mean an alphabet of the form $\mathpzc{A} = \{a,b\}$ where $b > a+2$. This is where we prove \Cref{sameComplexity} (iii) to (v), relying on the following lemma.

\begin{lemma}\label{inequalities}
	Let $\mathpzc{A} = \{a,b\}$ be an even alphabet.
	
	\textup{(i)} If $b \leq \frac{1}{2}a^2 + 2a$, then we have $\frac{a+b}{2}a(a+2) > b(b-2)$.
	
	\textup{(ii)} If $b \geq \frac{1}{2}a^2 + 2a + 2$, then we have $\frac{a+b}{2}a(a+2) + 2a < b(b-2)$.
	
	\textup{(iii)} If $b \leq \frac{1}{2}a^2 + 3a$, then we have $\frac{a+b}{2}a(a+4) > b(b-2)$.
	
	\textup{(iv)} If $b \geq \frac{1}{2}a^2 + 3a + 2$, then we have $\frac{a+b}{2}a(a+4) < b(b-2)$.
\end{lemma}
\begin{proof}
	(i) Set $f : x \mapsto x^2 - \left[\frac{1}{2}a^2+a+2\right]x - \left[\frac{1}{2}a^3+a^2\right]$ over $[4,+\infty)$, so that
	$$\frac{a+b}{2}a(a+2) > b(b-2) \iff f(b) < 0 \iff b < g(a) := \frac{1}{4}a^2+\frac{1}{2}a+1+\sqrt{\frac{1}{16}a^4 + \frac{3}{4}a^3 + \frac{7}{4}a^2 + a + 1} = \frac{1}{2}a^2+2a + h(a)$$
	where
	$$h(a) := \sqrt{\frac{1}{16}a^4 + \frac{3}{4}a^3 + \frac{7}{4}a^2 + a + 1} - \left[\frac{1}{4}a^2+\frac{3}{2}a-1\right] = \frac{4a}{\sqrt{\frac{1}{16}a^4 + \frac{3}{4}a^3 + \frac{7}{4}a^2 + a + 1} + \left[\frac{1}{4}a^2+\frac{3}{2}a-1\right]}.$$
	We observe that $h(a) > 0$, so if $b \leq \frac{1}{2}a^2 + 2a$ we have $b < \frac{1}{2}a^2 + 2a + h(a) = g(a)$ and $\frac{a+b}{2}a(a+2) > b(b-2)$.
	
	(ii) Set $f : x \mapsto x^2 - \left[\frac{1}{2}a^2+a+2\right]x - \left[\frac{1}{2}a^3+a^2+2a\right]$ over $[4,+\infty)$, so that
	$$\frac{a+b}{2}a(a+2) + 2a < b(b-2) \iff f(b) > 0 \iff b > g(a) := \frac{1}{4}a^2+\frac{1}{2}a+1+\sqrt{\frac{1}{16}a^4 + \frac{3}{4}a^3 + \frac{7}{4}a^2 + 3a + 1} = \frac{1}{2}a^2+2a + h(a)$$
	where
	$$h(a) := \sqrt{\frac{1}{16}a^4 + \frac{3}{4}a^3 + \frac{7}{4}a^2 + 3a + 1} - \left[\frac{1}{4}a^2+\frac{3}{2}a-1\right] = \frac{6a}{\sqrt{\frac{1}{16}a^4 + \frac{3}{4}a^3 + \frac{7}{4}a^2 + 3a + 1} + \left[\frac{1}{4}a^2+\frac{3}{2}a-1\right]}.$$
	We observe that the map $a \mapsto h(a)$ is decreasing and that $h(2) \approx 1.583$, so we always have $h(a) < 2$. Then, if $b \geq \frac{1}{2}a^2 + 2a + 2$ we have $b > \frac{1}{2}a^2 + 2a + h(a) = g(a)$ and $\frac{a+b}{2}a(a+2) + 2a < b(b-2)$.
	
	(iii) (iv) Set $f : x \mapsto x^2 - \left[\frac{1}{2}a^2+2a+2\right]x - \left[\frac{1}{2}a^3+2a^2\right]$ over $[4,+\infty)$, so that
	$$\frac{a+b}{2}a(a+4) > b(b-2) \iff f(b) < 0 \iff b < g(a)$$
	and
	$$\frac{a+b}{2}a(a+4) < b(b-2) \iff f(b) > 0 \iff b > g(a)$$
	where
	$$g(a) := \frac{1}{4}a^2 + a + 1 + \sqrt{\frac{1}{16}a^4 + a^3 + \frac{7}{2}a^2 + 2a + 1} = \frac{1}{2}a^2 + 3a + h(a)$$
	and
	$$h(a) := \sqrt{\frac{1}{16}a^4 + a^3 + \frac{7}{2}a^2 + 2a + 1} - \left[\frac{1}{4}a^2+2a-1\right] = \frac{6a}{\sqrt{\frac{1}{16}a^4 + a^3 + \frac{7}{2}a^2 + 2a + 1} + \left[\frac{1}{4}a^2+2a-1\right]}.$$
	We observe that $h(a) > 0$, and that $a \mapsto h(a)$ is decreasing and $h(2) \approx 1.292$ so we always have ${h(a) < 2}$. Finally, if $b \leq \frac{1}{2}a^2 + 3a$ (resp. $b \geq \frac{1}{2}a^2 + 3a + 2$) then we have $b < \frac{1}{2}a^2 + 2a + h(a) = g(a)$ (resp. ${b > \frac{1}{2}a^2 + 2a + h(a) = g(a)}$) and $\frac{a+b}{2}a(a+4) > b(b-2)$ (resp. $\frac{a+b}{2}a(a+4) < b(b-2)$).
\end{proof}

\begin{proof}[Proof of \Cref*{sameComplexity} \textup{(iii)} to \textup{(v)}]
	Let $\mathpzc{A} = \{a,b\}$ be an even alphabet such that $b > a+2$. Then \Cref{sameComplexity} (i) states $p(n) = p_x(n)$ where $x := \mathpzc{E}(a^\omega)$. We are going to compute the values of $s(n) := p(n+1)-p(n)$ by ordering the lengths of the bispecial factors of $x$.
	
	\Cref{shortBispecials} states that the short strong (resp. weak) bispecial factors of $x$ are $\varepsilon$, $a^a$ and $b^a$ (resp. $a^{b-1}$ and $b^{b-1}$). This yields the ordering
	\begin{equation}\label{shortOrdering}
		\lvert \varepsilon \rvert = 0 < \left\lvert a^a \right\rvert = \left\lvert b^a \right\rvert = a < \left\lvert a^{b-1} \right\rvert = \left\lvert b^{b-1} \right\rvert = b-1.
	\end{equation}
	Also, \Cref{evenLongBispecials} states that the long strong (resp. weak) bispecial factors of $x$ are the words $\left(\mathcal{P}_{f,b}\right)^{n+1}(\varepsilon)$, $\left(\mathcal{P}_{f,b}\right)^{n+1} (a^a)$ and $\left(\mathcal{P}_{f,b}\right)^{n+1} (b^a)$ (resp. $\left(\mathcal{P}_{f,b}\right)^n \mathcal{P}_{f,a}(\varepsilon)$, $\left(\mathcal{P}_{f,b}\right)^{n+1} (a^{b-1})$ and $\left(\mathcal{P}_{f,b}\right)^{n+1} (b^{b-1})$) for all $n \geq 0$. However, \Cref{evenLength} implies that the words $\left(\mathcal{P}_{f,b}\right)^{n+1}(\varepsilon)$ and $\left(\mathcal{P}_{f,b}\right)^n \mathcal{P}_{f,a}(\varepsilon)$ have the same length for all $n \geq 0$, therefore their contribution to $s(n)$ cancel one another. It remains to order the following lengths that we define for $n \geq 0$:
	\begin{align*}
		& p_n := \left\lvert \left(\mathcal{P}_{f,b}\right)^{n+1} (a^a) \right\rvert, & & q_n := \left\lvert \left(\mathcal{P}_{f,b}\right)^{n+1} (a^{b-2}) \right\rvert, & & r_n := \left\lvert \left(\mathcal{P}_{f,b}\right)^{n+1} (b^a) \right\rvert, & & s_n := \left\lvert \left(\mathcal{P}_{f,b}\right)^{n+1} (b^{b-2}) \right\rvert.
	\end{align*}
	\Cref{evenLength} yields
	\begin{align*}
		& p_0 = a(a+2), & & q_0 = ab, & & r_0 = a(b+2), & & s_0 = b(b-2)+2a, & & p_1 = \frac{a+b}{2}a(a+2) + 2a,
	\end{align*}
	and implies that the four sequences $(p_n)_{n \geq 0}$, $(q_n)_{n \geq 0}$, $(r_n)_{n \geq 0}$ and $(s_n)_{n \geq 0}$ satisfy the recurrence relation $\ell_{n+1} = \frac{a+b}{2} \ell_n + 2a$.
	
	(iii) If $b \in \left\llbracket a+4,\frac{1}{2}a^2+2a \right\rrbracket$, \Cref{inequalities} (i) yields $s_0 < p_1$ and we deduce the ordering
	\begin{equation*}
		p_0 < q_0 < r_0 < s_0 < p_1
	\end{equation*}
	so an immediate induction on $n$ with the recurrence relation yields 
	\begin{equation}\label{longOrdering3}
		p_n < q_n < r_n < s_n < p_{n+1}
	\end{equation}
	for all $n \geq 0$. Finally, to link \Cref{shortOrdering,longOrdering3}, we have $a \leq b-4$ so $\frac{a+b}{2} \leq b-2$, and \Cref{inequalities} (i) yields $p_0 = a(a+2) > \frac{2}{a+b} b(b-2) \geq b > b-1 = \left\lvert a^{b-1} \right\rvert = \left\lvert b^{b-1} \right\rvert$. We deduce the complete ordering
	\begin{equation*}
		\lvert \varepsilon \rvert < \left\lvert a^a \right\rvert = \left\lvert b^a \right\rvert < \left\lvert a^{b-1} \right\rvert = \left\lvert b^{b-1} \right\rvert < p_0 < q_0 < r_0 < s_0 < p_1 < ... < p_n < q_n < r_n < s_n < p_{n+1} < ...
	\end{equation*}
	Since $p_n$ and $r_n$ (resp. $q_n$ and $s_n$) are each the length of a single strong (resp. weak) bispecial factor of $x$, we compute the values of $s(n)$:
	\begin{align*}
		s(n) = \begin{cases} 1 \textrm{ if } n = 0 \\ 2 \textrm{ if } n \in \llbracket 1,a \rrbracket \\ 4 \textrm{ if } n \in ~\llbracket a+1,b-1 \rrbracket \\ 2 \textrm{ if } n \in \llbracket b,p_0 \rrbracket \end{cases} \textrm{ and~~~} & s(n) = \begin{cases} 3 \textrm{ if } n \in ~\llbracket p_k+1,q_k \rrbracket \\ 2 \textrm{ if } n \in ~\llbracket q_k+1,r_k \rrbracket \\ 3 \textrm{ if } n \in ~\llbracket r_k+1,s_k \rrbracket \\ 2 \textrm{ if } n \in ~\llbracket s_k+1,p_{k+1} \rrbracket \end{cases} \textrm{ for some } k \geq 0.
	\end{align*}
	Finally, because the four sequences $(p_n)_{n \geq 0}$, $(q_n)_{n \geq 0}$, $(r_n)_{n \geq 0}$ and $(s_n)_{n \geq 0}$ satisfy the recurrence relation $\ell_{n+1} = \frac{a+b}{2} \ell_n + 2a$, the lengths of the intervals $\llbracket p_k+1,q_k \rrbracket$, $\llbracket q_k+1,r_k \rrbracket$, $\llbracket r_k+1,s_k \rrbracket$, $\llbracket s_k+1,p_{k+1} \rrbracket$ satisfy the recurrence relation $\ell_{k+1} = \frac{a+b}{2} \ell_k$. In particular, $s(n)$ is equal to $2$ (resp. $3$) with non-zero proportion, which implies that $\displaystyle 2 < \liminf_{n \rightarrow \infty} \frac{p_x(n)}{n} < \limsup_{n \rightarrow \infty} \frac{p_x(n)}{n} < 3$.
	
	(iv) (v) If $b \in \left\llbracket \frac{1}{2}a^2+2a+2,+\infty \right)$, \Cref{inequalities} (ii) yields $s_0 > p_1$. We also have $(a+b)(a+2) > 2b$ so $p_1 = \frac{a+b}{2}a(a+2) + 2a > a(b+2) = r_0$, and $(a+b)a > 2(b-2)$ so $q_1 = \frac{a+b}{2}ab + 2a > b(b-2) + 2a = s_0$. We deduce the ordering
	\begin{equation*}
		p_0 < q_0 < r_0 < p_1 < s_0 < q_1
	\end{equation*}
	so an immediate induction on $n$ with the recurrence relation yields 
	\begin{equation}\label{longOrdering4}
		q_n < r_n < p_{n+1} < s_n < q_{n+1}
	\end{equation}
	for all $n \geq 0$. We also have $(a+b)a > 2(b-1)$ so $p_0 = \frac{a+b}{2}a(a+2) + 2a > b-1$. We deduce the complete ordering
	\begin{equation*}
		\lvert \varepsilon \rvert < \left\lvert a^a \right\rvert = \left\lvert b^a \right\rvert < \left\lvert a^{b-1} \right\rvert = \left\lvert b^{b-1} \right\rvert < p_0 < q_0 < r_0 < p_1 < s_0 < q_1 <... < q_n < r_n < p_{n+1} < s_n < q_{n+1} < ...
	\end{equation*}
	and we deduce the values of $s(n)$:
	\begin{align*}
		s(n) = \begin{cases} 1 \textrm{ if } n = 0 \\ 2 \textrm{ if } n \in \llbracket 1,a \rrbracket \\ 4 \textrm{ if } n \in ~\llbracket a+1,b-1 \rrbracket \\ 2 \textrm{ if } n \in \llbracket b,p_0 \rrbracket \\ 3 \textrm{ if } n \in \llbracket p_0+1,q_0 \rrbracket \end{cases} \textrm{ and~~~} & s(n) = \begin{cases} 2 \textrm{ if } n \in ~\llbracket q_k+1,r_k \rrbracket \\ 3 \textrm{ if } n \in ~\llbracket r_k+1,p_{k+1} \rrbracket \\ 4 \textrm{ if } n \in ~\llbracket p_{k+1}+1,s_k \rrbracket \\ 3 \textrm{ if } n \in ~\llbracket s_k+1,q_{k+1} \rrbracket \end{cases} \textrm{ for some } k \geq 0.
	\end{align*}
	Because the four sequences $(p_n)_{n \geq 0}$, $(q_n)_{n \geq 0}$, $(r_n)_{n \geq 0}$ and $(s_n)_{n \geq 0}$ satisfy the recurrence relation $\ell_{n+1} = \frac{a+b}{2} \ell_n + 2a$, the lengths of the intervals $\llbracket p_k+1,q_k \rrbracket$, $\llbracket q_k+1,r_k \rrbracket$, $\llbracket r_k+1,s_k \rrbracket$, $\llbracket s_k+1,p_{k+1} \rrbracket$ satisfy the recurrence relation $\ell_{k+1} = \frac{a+b}{2} \ell_k$. In particular, $s(n)$ is equal to $2$ (resp. $4$) with non-zero proportion, which implies that $\displaystyle 2 < \liminf_{n \rightarrow \infty} \frac{p_x(n)}{n} < \limsup_{n \rightarrow \infty} \frac{p_x(n)}{n} < 4$.
	
	Also, the length of the intervals $\llbracket q_k,q_{k+1} \rrbracket$, $\llbracket r_k,r_{k+1} \rrbracket$, $\llbracket p_{k+1},p_{k+2} \rrbracket$ and $\llbracket s_k,s_{k+1} \rrbracket$ tends to infinity when $k$ grows. We deduce that ${\displaystyle \limsup_{n \rightarrow \infty} \frac{p_x(n)}{n} = \max\left(\frac{p(q_{k+1}+1)-p(q_k+1)}{q_{k+1}-q_k}, \frac{p(s_{k+1}+1)-p(s_k+1)}{s_{k+1}-s_k} \right)}$ where
	\begin{align*}
		\frac{p(q_{k+1}+1)-p(q_k+1)}{q_{k+1}-q_k} = &~ \frac{\displaystyle \sum_{i = q_k+1}^{q_{k+1}} s(i)}{q_{k+1}-q_k} = 3 + \frac{\displaystyle \sum_{i = q_k+1}^{q_{k+1}} \left[s(i) - 3\right]}{q_{k+1}-q_k} = 3 + \frac{(s_k-p_{k+1})-(r_k-q_k)}{q_{k+1}-q_k} \\
		= &~ 3 + \frac{(s_0-p_1) - (r_0 - q_0)}{q_1-q_0} = 3 + \frac{b(b-2)-\frac{a+b}{2}a(a+2)-2a}{q_1-q_0}
	\end{align*}
	and
	$$\frac{p(s_{k+1}+1)-p(s_k+1)}{s_{k+1}-s_k} = 3 + \frac{(s_1-p_2)-(r_1-q_1)}{s_1-s_0} = 3 + \frac{a+b}{2}\frac{b(b-2)-\frac{a+b}{2}a(a+2)-2a}{s_1-s_0},$$
	and \Cref{inequalities} (ii) yields $\displaystyle \limsup_{n \rightarrow \infty} \frac{p_x(n)}{n} > 3$.
	
	Similarly, we deduce that $\displaystyle \liminf_{n \rightarrow \infty} \frac{p_x(n)}{n} = \max\left( \frac{p(r_{k+1}+1)-p(r_k+1)}{r_{k+1}-r_k}, \frac{p(p_{k+2}+1)-p(p_{k+1}+1)}{p_{k+2}-p_{k+1}} \right)$ where
	$$\frac{p(r_{k+1}+1)-p(r_k+1)}{r_{k+1}-r_k} = 3 + \frac{(s_0-p_1)-(r_1-q_1)}{r_1-r_0} = 3 + \frac{b(b-2)-\frac{a+b}{2}a(a+4)}{r_1-r_0}$$
	and
	$$\frac{p(p_{k+2}+1)-p(p_{k+1}+1)}{p_{k+2}-p_{k+1}} = 3 + \frac{(s_0-p_1)-(r_1-q_1)}{p_2-p_1} = 3 + \frac{b(b-2)-\frac{a+b}{2}a(a+4)}{p_2-p_1}.$$
	Now, if $b \in \left\llbracket \frac{1}{2}a^2+2a+2,\frac{1}{2}a^2+3a \right\rrbracket$ (resp. if $b \in \left\llbracket \frac{1}{2}a^2+3a+2,+\infty \right)$), \Cref{inequalities} (iii) (resp. \Cref{inequalities} (iv)) yields $\displaystyle \liminf_{n \rightarrow \infty} \frac{p_x(n)}{n} < 3$ (resp. $\displaystyle \liminf_{n \rightarrow \infty} \frac{p_x(n)}{n} > 3$).
\end{proof}

\subsection{Linear factor complexity}

\Cref{linearCorollary} naturally follows from \Cref{URLR,sameComplexity}.

\begin{proof}[Proof of \Cref*{linearCorollary}]
	If $x$ is a smooth sequence in $\mathscr{M}$, then \Cref{aperiodic} and \Cref{URLR} \textup{(ii)} state that $x$ is aperiodic and linearly recurrent. Now the following statements apply to $x$: an aperiodic sequence has at least linear factor complexity (\cite{MH}) and a linearly recurrent sequence has at most linear factor complexity (\cite[Theorem 24 \textit{i)}]{DHS}). We deduce that $p_x(n) = \Theta(n)$.
	
	If $x$ is a smooth sequence in $\mathscr{M}' \backslash \mathscr{M}$, i.e. $x$ is defined over the alphabet $\{1,3\}$ and does not have property $\mathfrak{P}1$ or $\mathfrak{P}2$, then \Cref{sameComplexity} (i) and (ii) imply that $p_x(n) = p(n) = \Theta(n)$.
\end{proof}

\section{Frequencies in smooth sequences}\label{frequencies}

\subsection{Uniform frequencies}

Similarly to \Cref{linearCorollary}, \Cref{frequenciesCorollary} naturally follows from \Cref{URLR,sameComplexity}.

\begin{proof}[Proof of \Cref*{frequenciesCorollary}]
	If $x$ is a smooth sequence in $\mathscr{M}$, then \Cref{URLR} \textup{(ii)} states that $x$ is linearly recurrent. Now the following statement applies to $x$: a linearly recurrent sequence has uniform frequencies (\cite[Theorem 15]{Durand00}).
	
	If $x$ is a smooth sequence in $\mathscr{M}' \backslash \mathscr{M}$, i.e. $x$ is defined over the alphabet $\{1,3\}$ and does not have property $\mathfrak{P}2$, so \Cref{sameComplexity} (i) and (ii) imply that $\displaystyle\limsup_{n \rightarrow \infty} \frac{p_x(n)}{n} =  \displaystyle\limsup_{n \rightarrow \infty} \frac{p(n)}{n} = 2$. Now the following statement applies to $x$: a sequence $y$ that satisfies $\displaystyle\limsup_{n \rightarrow \infty} \frac{p_y(n)}{n} < 3$ has uniform frequencies (\cite[Theorem 1.5]{Boshernitzan}).
\end{proof}

\begin{remark}
	Over any even alphabet $\mathpzc{A} = \{a,b\}$ such that $b \in \llbracket a+2, \frac{1}{2}a^2 + 2a \rrbracket$, we can also prove \Cref{frequenciesCorollary} in the same way as for the sequences in $\mathscr{M}' \backslash \mathscr{M}$: \Cref{sameComplexity} \textup{(i)} to \textup{(iii)} imply that a smooth sequence $x$ over such an alphabet satisfies $\displaystyle\limsup_{n \rightarrow \infty} \frac{p_x(n)}{n} < 3$ so \cite[Theorem 1.5]{Boshernitzan} yields uniform frequencies.
\end{remark}

\subsection{Letter frequencies}

In this subection we prove \Cref{smoothFreq}. We begin by considering smooth sequences over even alphabets, this is \Cref{smoothFreq} (i).

\begin{proof}[Proof of \Cref*{smoothFreq} \textup{(i)}]
	Let $x$ be a smooth sequence over an even alphabet $\mathpzc{A} = \{a,b\}$, and let $t = \mathcal{T}(x)$, so that $x = \pi_{x_0} \circ \mathcal{F} \circ S(t)$. We observe that, for any sequence $y \in \Sigma^\mathbb{N}$, the sequence $\sigma_{t_1}(y)$ has equal frequencies of $A$'s and $B$'s. In particular this holds for $\mathcal{F} \circ S(t) = \sigma_{t_1} \circ \mathcal{F} \circ S^2(t)$, and we deduce that $x = \pi_{x_0} \circ \mathcal{F} \circ S(t)$ has equal frequencies of $a$'s and $b$'s.
\end{proof}

Then we consider smooth sequences over odd alphabets, this is \Cref{smoothFreq} (ii) and (iii).

\begin{proof}[Proof of \Cref*{smoothFreq} \textup{(ii)}]
	Let $x$ be a smooth sequence over an odd alphabet $\mathpzc{A} = \{a,b\}$, and suppose that $\mathcal{T}(x) = cu^\omega$ for some $c \in \mathpzc{A}$ and some $u \in \mathpzc{A}^n$ ($n \geq 1$), and that and that $u$ is not of the form $(1b)^n$ or $(b1)^n$ if $a = 1$. Then \Cref{smoothMorphic} states that $\mathcal{F}(u^\omega) = \sigma^\infty(D)$ and $x = \pi_c \circ \mathcal{F}(u^\omega)$ where $\sigma := \sigma_{u_1} ... \sigma_{u_n}$ and $D \in \Sigma$. If $a = 1$, the condition on $u$ implies that there exists $n \geq 1$ such that $(1b)^n \nsqsubset u^\omega$, therefore \Cref{FLR} states that $\mathcal{F}(u^\omega) = \sigma^\infty(D)$ is a proper primitive $S$-adic sequence. We deduce that the substitution $\sigma$ is primitive, so, following \cite{Queffelec}, the letter frequencies of the purely morphic sequence $\sigma^\infty(D)$ are given by the normalized right Perron eigenvector of $M_\sigma$, that we denote $\left(\begin{matrix} \alpha \\ \beta \\ \gamma \end{matrix}\right)$. Since $\pi_c(A) = cc$, $\pi_c(B) = c \overline{c}$ and $\pi_{c}(D) = \overline{c} \overline{c}$, the frequency of $c's$ in $x = \pi_c \left(\sigma^\infty(D)\right)$ is $\frac{2\alpha+\beta}{2} = \alpha+\frac{1}{2}\beta$.
\end{proof}

\begin{proof}[Proof of \Cref*{smoothFreq} \textup{(iii)}]
	Let $x$ be the smooth sequence over an odd alphabet $\mathpzc{A} = \{1,b\}$.
	
	If $\mathcal{T}(x) = (1b)^\omega$, then \Cref{1bb1morphic} \textup{(i)} states that $x = \mathcal{E}((1b)^\omega) = \pi_1 (\sigma^\infty(A))$ where $\sigma : A \mapsto A^\frac{b-1}{2} B$, ${B \mapsto A^\frac{b-1}{2} B^\frac{b+1}{2}}$. The substitution $\sigma$ is primitive with incidence matrix $\left(\begin{matrix} \frac{b-1}{2} & \frac{b-1}{2} \\ 1 & \frac{b+1}{2} \end{matrix}\right)$, it has dominant eigenvalue $\frac{b+\sqrt{2b-1}}{2}$ with right eigenvector $\left(\begin{matrix} 1-\frac{2}{\sqrt{2b-1}+1} \\ \frac{2}{\sqrt{2b-1}+1}\end{matrix}\right)$ which provides the frequencies of $A$'s and $B$'s in $\sigma^\infty(A)$. Finally, since $\pi_1(A) = 11$ and $\pi_1(B) = 1b$, the frequency of $b$'s in $x = \pi_1(\sigma^\infty(A))$ is $\frac{1}{2} \frac{2}{\sqrt{2b-1}+1} = \frac{1}{\sqrt{2b-1}+1}$.
	
	If $\mathcal{T}(x) = (b1)^\omega$, the proof is identical thanks to \Cref{1bb1morphic} \textup{(ii)}.
\end{proof}

\section*{Acknowledgements}

This work was supported by ANR-22-CE40-0011 project Inside Zero Entropy Systems.

\bibliographystyle{plainurl}
\bibliography{refs}

\end{document}